\documentclass[reqno]{amsart}
\pdfoutput=1\relax\pdfpagewidth=8.26in\pdfpageheight=11.69in\pdfcompresslevel=9
\usepackage{mathptmx,amsmath,mathbbol,lineno}
\usepackage{hyperref}
\usepackage{amsfonts}
\usepackage{amssymb}
\usepackage{amscd}
\usepackage{latexsym}
\usepackage{graphics}
\usepackage{subfig}
\usepackage{float}
\usepackage[all,knot,poly,arc,web]{xy}
\usepackage{ifthen}
\usepackage{epsfig}
\usepackage{tikz}

\newtheorem{theorem}{Theorem}
\newtheorem{Hypotheses (H)}{Hypotheses (H)}
\newtheorem{definition}{Definition}

\newtheorem{lemma}{Lemma}
\newtheorem{remark}{Remark}

\newtheorem{prop}{Proposition}

\def\neweq#1{\begin{equation}\label{#1}}
\def\endeq{\end{equation}}

\def\phi{\varphi}

\date{}

\def\pp{\overrightarrow{p}(\cdot) }
\def\qq{\overrightarrow{q}(\cdot) }

\def\fun(#1,#2,#3){\mathcal{E}_{_{#1}}(#2, #3)}
\def\sob(#1,#2){W^{1,#1}(#2)}
\def\forma(#1){\fun({\pp,\qq},\cdot,\cdot)}

\def\l2{L^2(\Omega)}

\def\funn(#1,#2,#3,#4,#5){\langle J_{#1} #2, #3 \rangle_{_{#4,#5}}}
\def\nn(#1,#2){\|#1\|_{_{#2}}}

\newcounter{appendix}
\newenvironment{apequation}{
\addtocounter{equation}{-1}
\refstepcounter{appendix}

\begin{equation}}
{\end{equation}}
\newtheorem{adefinicion}{Definition A}

\newtheorem{ateorema}{Theorem A}

\AtBeginDocument{{\noindent\small
{\em {\color{blue}Calculus of Variations and Partial Differential Equations} {\color{orange}(to appear)}}}
\vspace{9mm}}

\begin{document}

\title[Diffusion over ramified domains]{Diffusion over ramified domains: solvability and fine regularity}

\author[K. Silva\,-\,P\'erez,\,\,\,A. V\' elez\,-\,Santiago]{Kevin Silva\,-\,P\'erez,\,\,\,Alejandro V\' elez\,-\,Santiago}

\address{Kevin Silva\,-\,P\'erez\hfill\break
Department of Mathematics\\
University of Puerto Rico at R\'io Piedras\\
San Juan, PR \,00925}
\email{kevin.silva1@upr.edu}

\address{Alejandro V\'elez\,-\,Santiago\hfill\break
Department of Mathematics\\
University of Puerto Rico at R\'io Piedras\\
San Juan, PR \,00925}
\email{alejandro.velez2@upr.edu,\,\,\,dr.velez.santiago@gmail.com}


\numberwithin{equation}{section}

\begin{abstract}
We consider a domain $\Omega\subseteq\mathbb{R\!}^{\,2}$ with branched fractal boundary $\Gamma^{\infty}$ and parameter $\tau\in[1/2,\tau^{\ast}]$ introduced by Achdou and Tchou \cite{ACH08}, for $\tau^{\ast}\simeq 0.593465$, which acts as an idealization of the bronchial trees in the lungs systems. For each $\tau\in[1/2,\tau^{\ast}]$, the corresponding region $\Omega$ is a non-Lipschitz domain, which attains its roughest structure at the critical value $\tau=\tau^{\ast}$ in such way that in this endpoint parameter the region $\Omega$ fails to be an extension domain, and its ramified boundary $\Gamma^{\infty}$ is not post-critically finite. Then, we investigate a model equation related to the diffusion of oxygen through the bronchial trees by considering the realization of a generalized diffusion equation $$\frac{\partial u}{\partial t}-\mathcal{A} u+\mathcal{B}u\,=\,f(t,x)\,\,\,\,\,\,\,\,\textrm{in}\,\,\,(0,\infty)\times\Omega$$ with inhomogeneous mixed-type boundary conditions $$\displaystyle\frac{\partial u}{\partial\nu_{_{\mathcal{A}}}}+\beta u\,=\,g(x,t)\,\,\,\,\,\,\,\textrm{on}\,\,(0,\infty)\times\Gamma^{\infty},\,\,\,\,\,\,\,\,\,\,\,\,\,\,\,\,u=0\,\,\,\,\,\,\,\textrm{in}\,\,\,(0,\infty)\times(\partial\Omega\setminus\Gamma^{\infty}),$$ and $u(x,0)=u_0\in C(\overline{\Omega})$,  where $\mathcal{A}$ is an uniformly elliptic second-order (non-symmetric) differential operator with bounded measurable coefficients, $\mathcal{B}$ is as a lower-order (non-symmetric) differential operator with unbounded measurable coefficients,\, $\displaystyle\frac{\partial u}{\partial\nu_{_{\mathcal{A}}}}$ stands as a generalized notion of a normal derivative over rough surfaces (in the sense of Definition \ref{Def-gen-normal}), and $\beta\in L^s_{\mu}(\Gamma^{\infty})^+$ with $\displaystyle{\textrm{ess}\inf_{x\in\Gamma^{\infty}}}|\beta(x)|\geq\beta_0$ for a sufficiently large constant $\beta_0>0$, and $s>1$. Under minimal assumptions, we first show that the stationary version of the above  diffusion equation in uniquely solvable, and that the corresponding weak solution in globally H\"older continuous on $\overline{\Omega}$. Since we are including the critical case $\tau=\tau^{\ast}$, this is the first time in which global uniform continuity of weak solutions of a Robin-type boundary value problem is attained over a non-extension domain. Furthermore, after two transitioning procedures, we prove the unique solvability of the inhomogeneous time-dependent diffusion equation, and we show that the corresponding weak solution is globally uniformly continuous over $[0,T]\times\overline{\Omega}$ for each fixed parameter $T>0$.
\end{abstract}
\maketitle

\section{Introduction and statements of main results}\label{sec1}

\subsection{Bronchial trees and ramified domains}\label{subsec1.1}

\indent In the respiratory system of mammals, the collective term \textit{bronchial tree} refers to the bronchi and all their branches. The bronchi are the airways of the lower respiratory tract. At the level of the third or fourth thoracic vertebra, the trachea branches into the left and right main bronchus. The right main bronchus is shorter and runs more vertically than the left. Both bronchi divide further into secondary or lobar bronchi, which branch progressively to distribute respiratory air adequately to the left and right lungs. The terminal segment of each bronchus contains millions of alveoli (air sacs) in which gas exchange takes place. For more information related to the human bronchial tree and its fractal-type structure, we refer to the works in \cite{GLE11,T-S-S-H20,WEI-GOM62} (among others).

\begin{figure}[h]
    \centering
    \includegraphics[scale=0.2]{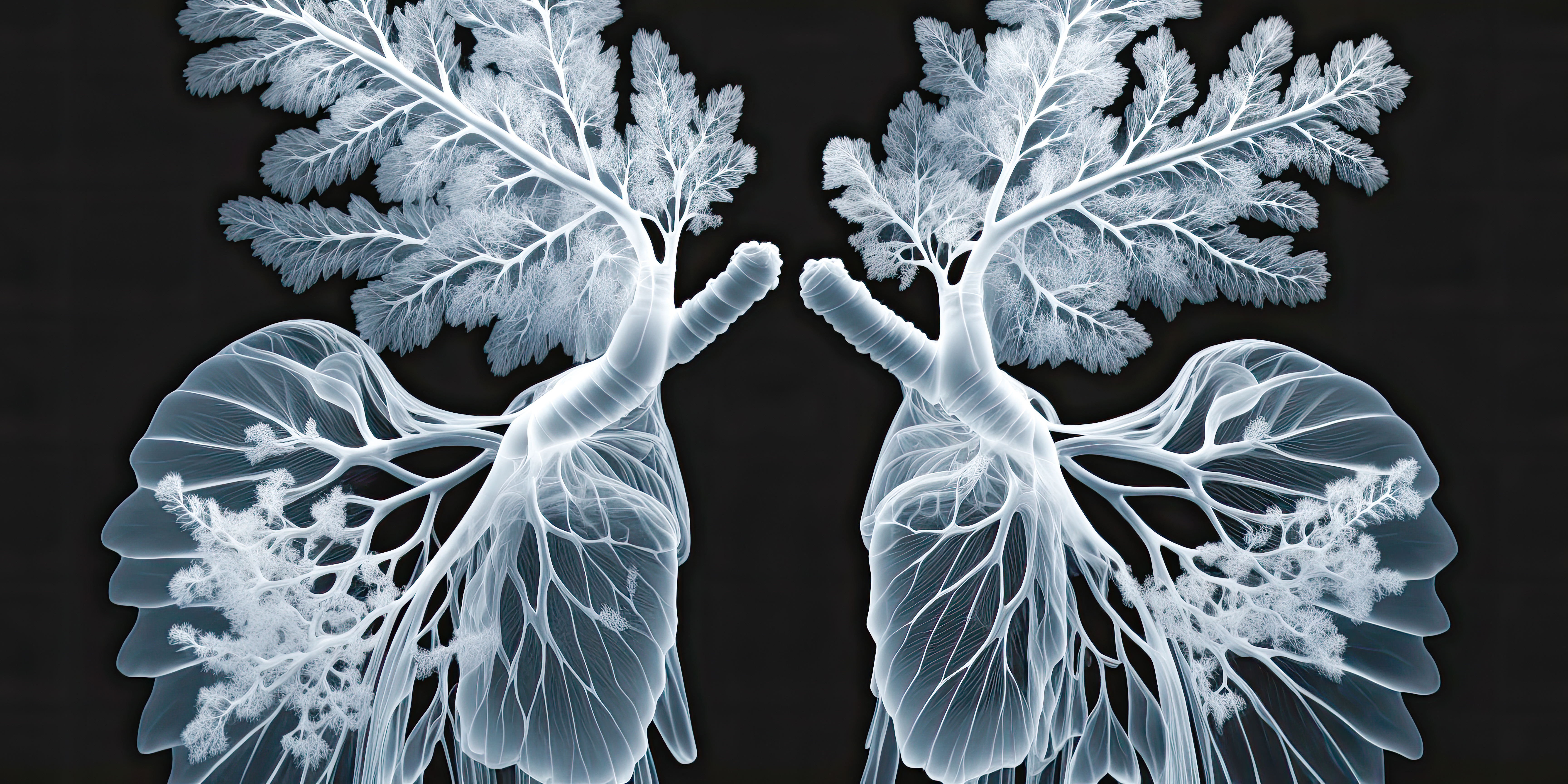}
    \caption{X-ray of human bronchial trees.}
    \label{x-ray}
\end{figure}

\indent The mammalian bronchial tree is designed to provide adequate gas flow to the alveoli, minimal entropy production in respiratory mechanics, and minimal expenditure of matter and energy. However, the airways represent only a part of the respiratory system and their geometry must be adapted to the function of the entire system by solving the problem of distributing an inspired volume of air over a large surface area in a confined volume. Therefore, the topology of bronchi exhibits the features of space utilisation with its progressive branching and a reduction in bronchial diameter associated with a fractal geometry; see for instance \cite{MEASURING-BRONCHIAL}. The effect of scale is also studied to verify self-similarity in the fractal geometry of the bronchial tree.\\
\indent Providing a good measure (and fractal dimension) to the bronchial trees guarantees advantages in the biology of these, allowing optimization of gas exchange, and it will also be helpful to investigate more precisely the ability of the airways to take up space by branching, and quickly reach the alveolar surface. We can also use these information to investigate the heterogeneity of alveolar flow and ventilation and the ability to buffer developmental defects (e.g. \cite{Frac-pulmon2}). Additionally, in \cite{Frac-pulmon}, it was suggested that fractal patterns of endothelial junction organization determine the functional selectivity of solute transport through the pulmonary capillary network.

\subsection{Domain with ramified fractal boundary}\label{subsec1.2}

\indent Bronchial trees exhibit a ramified-like structure, so in order to study further the geometry and diffusion of oxygen, we turn our attention to the domain $\Omega\subseteq\mathbb{R\!}^{\,2}$ with ramified boundary $\Gamma^{\infty}$ introduced by Achdou and Tchou \cite{ACH08}. This set is precisely an idealization of the bronchial trees in the pulmonary system (in the sense of a projection into the $2$-dimensional setting), and it is dependent upon a parameter $\tau\in[1/2,\tau^{\ast}]$, where $\tau^{\ast}\simeq 0.593465$ is the solution of the equation (\ref{tau-eq}), which is needed in order to avoid overlapping in the construction of the ramified set $\Omega$ (see \cite{ACH08}).\\
\indent When $1/2\leq\tau<\tau^{\ast}$, Achdou and Tchou proved that $\Omega$ is is a locally uniform domain (in the sense of Rogers \cite{ROGERS06}), but this situation fails when $\tau=\tau^{\ast}$. Moreover, when $\tau=\tau^{\ast}$, one has that the ramified boundary $\Gamma^{\infty}$ is not post-critically finite (e.g. \cite[Remark 4]{ACH08}). Therefore, although we will deal with ramified domains for arbitrary $\tau\in[1/2,\tau^{\ast}]$, the critical situation $\tau=\tau^{\ast}$ will be the main case considered in this paper, since in this case the geometry of the domain is more irregular, and as pointed out in \cite[Remark 4]{ACH08}, most of the known results related to function spaces on self-similar fractals are valid only for post-critically finite sets, which is not the case when $\tau=\tau^{\ast}$.

\begin{figure}[h!]
    \centering
\includegraphics[scale=0.35]{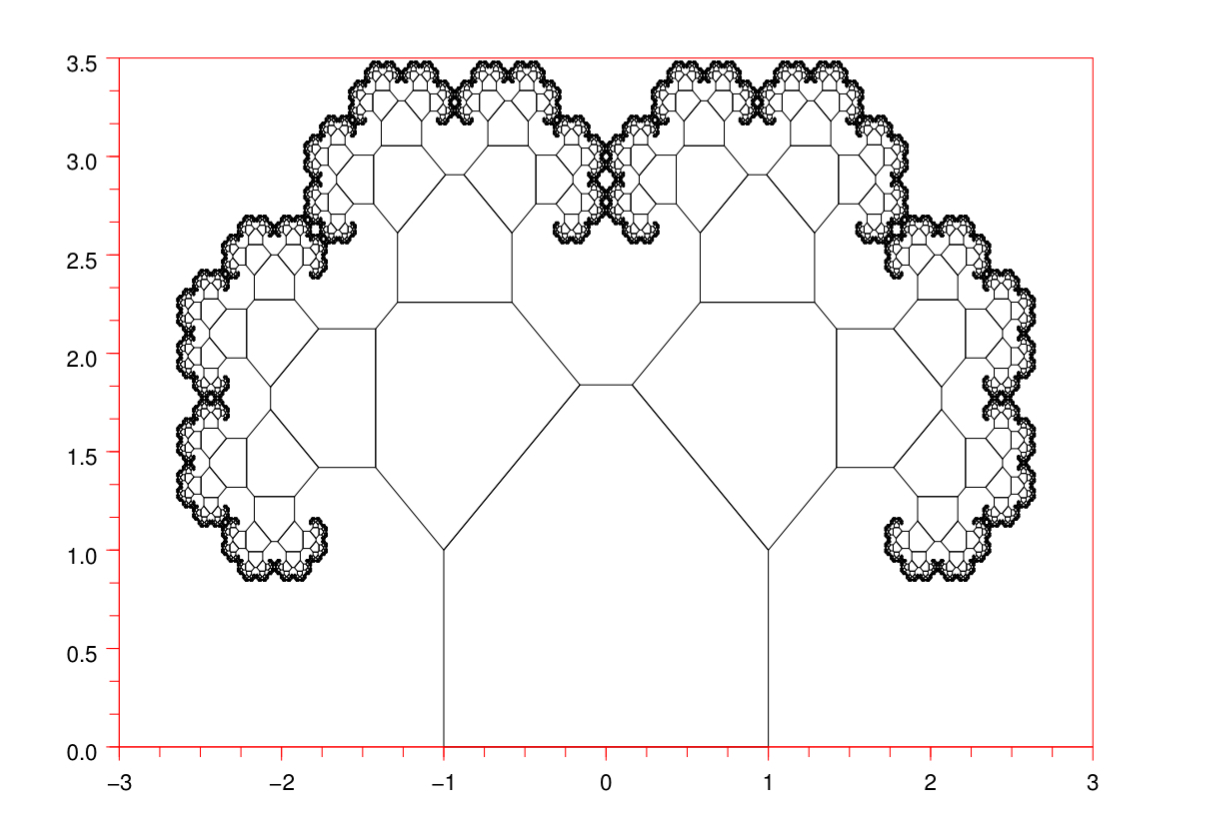}
    \caption{The domain $\Omega$ with ramified fractal boundary $\Gamma^{\infty}$ in the case $\tau=\tau^{\ast}$ (as in \cite{ACH08}).}
    \label{fig:my_label}
\end{figure}

\indent However, in \cite{ACH10,ACH08}, the authors were able to define function spaces and Sobolev spaces over $\Omega$, and were able to establish traces of Sobolev function into $L^p$-spaces and Besov spaces over the ramified boundary $\Gamma^{\infty}$, even in the critical case $\tau=\tau^{\ast}$. In fact, these results constitute a key framework in which one can define Robin-type boundary value problems, which play a crucial role in the theme of this paper.

\subsection{The diffusion equation}\label{subsec1.3}

\indent One of the main parts in this paper consists on studying the diffusion of medical sprays in the lungs system. Since the gas exchanges between the lungs and the circulatory system take place only in the last generations of the bronchial tree (the smallest structures), reasonable models for the diffusion of gases (such as oxygen) may involve inhomogeneous Neumann or Robin boundary conditions over the ramified boundary $\Gamma^{\infty}$. Similarly, the lungs are mechanically coupled to the diaphragm, which also implies inhomogeneous conditions over the boundary $\Gamma^{\infty}$ (in the case one is interested in a coupled fluid-structure model).\\
\indent Motivated by this, we consider the generalized
diffusion equation formally given by
\begin{equation}\label{1.01}
\displaystyle\frac{\partial u}{\partial t}-\mathcal{A}u+\mathcal{B} u\,=\,f(t,x)\,\,\,\,\,\,\,\,\,\,\textrm{in}\,\,\,(0,\infty)\times\Omega
\end{equation}
with inhomogeneous mixed Dirichlet-to-Robin boundary conditions
\begin{equation}\label{1.02}
\displaystyle\frac{\partial u}{\partial\nu_{_{\mathcal{A}}}}+\beta u\,=\,g(t,x)\,\,\,\,\textrm{on}\,\,\,(0,\infty)\times\Gamma^{\infty};\,\,\,\,\, u=0\,\,\,\textrm{over}\,\,\,(0,\infty)\times(\Gamma\setminus\Gamma^{\infty});\,\,\,\,\,u(0,x)=u_0\in L^2(\Omega).
\end{equation}
where $\mathcal{A}$ is an uniformly elliptic second-order (non-symmetric) differential operator with bounded measurable coefficients (for a concrete definition of $\mathcal{A}$, refer to section \ref{sec5}), $\mathcal{B}$ stands as a lower-order differential non-symmetric operator with unbounded measurable coefficients, $\displaystyle\frac{\partial u}{\partial\nu_{_{\mathcal{A}}}}$ stands as a generalized notion of a normal derivative in the sense of Definition \ref{Def-gen-normal}, and $\beta$ is a (possibly unbounded) measurable function defined over $\Gamma^{\infty}$. In the case of the Laplacian, Achdou and Tchou \cite{ACH07} established solvability and approximation results for an inhomogeneous elliptic Neumann problem (the case $\beta\equiv0$) over the ramified domain $\Omega$, but the realization of the non-stationary problem over $\Omega$ has not been fully investigated, up to the present time. Furthermore, our model equation (\ref{1.01}) with boundary conditions (\ref{1.02}) involve unbounded lower-order coefficient under minimal assumptions.\\
\indent Concerning the well-posedness of inhomogeneous parabolic problems with unbounded coefficients, there is little literature, almost non-existent in the case of non-smooth domains. One can recall first the results in \cite{LAD-SOL-URAL68,LAD-URAL68} over Lipschitz-type domains (for the Dirichlet problem), and their generalization to a broader class of domains recently achieved in \cite{KIM-RYU-WOO22} (for the Dirichlet problem, but under more general situations on the coefficients and function spaces). However, the results obtained in
\cite{KIM-RYU-WOO22,LAD-SOL-URAL68,LAD-URAL68} deal with inhomogeneous interior differential equation with homogeneous boundary conditions. For Neumann and Robin boundary value problems, Nittka \cite{NITTKA2014} developed a very useful procedure whose motivation resides in \cite{LAD-SOL-URAL68}, from which he was able to establish the solvability and global regularity for boundary value problems of type (\ref{1.01}), Robin boundary conditions, and unbounded lower-order terms over bounded Lipschitz domains. In particular, global regularity results were established by Nittka \cite{NITTKA2014} in the case of bounded Lipschitz domains, and his methods represent the main motivations for the establishment of our main results in this paper. Recently in \cite{APU-NAZ-PAL-SOF22}, a generalized solvability result and Sobolev norm estimate was developed for a class of parabolic Wentzell-type boundary value problem with unbounded coefficients over domains of class $C^{1,1}$, and in \cite{LAN-AVS-23-1}, a priori estimates were given for multiple boundary value problems involving various boundary conditions, unbounded measurable coefficients, over a general class of irregular domains possessing the extension property (in the sense of Definition \ref{ext-domain}).\\
\indent In the case of our problem, if $1/2\leq\tau<\tau^{\ast}$, by virtue of \cite{JON}, $\Omega$ is an extension domain, but this situation fails in the case $\tau=\tau^{\ast}$. Therefore, all results in the present paper are new and original, especially in the critical case $\tau=\tau^{\ast}$.

\subsection{Main results}\label{subsec1.4}

\indent The establishment of global regularity results for the diffusion equation (\ref{1.01}) with boundary conditions (\ref{1.02}) is divided into three main stages. At first, we consider the stationary version of problem (\ref{1.01}) with boundary conditions (\ref{1.02}), where we have the following global regularity result.

\begin{theorem}\label{Holder-continuity}
Given $f_i\in L^{P_i}(\Omega)$ and $g\in L^q_{\mu}(\Gamma^{\infty})$ for $p_i>\chi_{_{\{i=0\}}}+2\chi_{_{\{1\leq i\leq 2\}}}$ and $q>1$, the elliptic mixed boundary value problem
\begin{equation}
\label{E01}\left\{
\begin{array}{lcl}
\mathcal{A}u+\mathcal{B}u\,=\,f_0(x)-\displaystyle\sum^2_{j=1}\partial_{x_j}f_j(x)\,\,\,\,\,\textrm{in}\,\,\Omega;\\
\,\,\,\,u\,=\,0\,\,\,\,\,\,\indent\,\,\,\,\,\textrm{on}\,\,\Gamma\setminus\Gamma^{\infty};\\
\frac{\partial u}{\partial\nu_{_{\mathcal{A}}}}+\beta u\,=\,g(x)+\displaystyle\sum^2_{j=1}f_j(x)\nu_{\mu_j}\,\,\,\,\,\,\textrm{on}\,\,\Gamma^{\infty},\\
\end{array}
\right.
\end{equation}
possesses a unique weak solution $u\in\mathcal{V}_2(\Omega)$ (in the sense of (\ref{Existence-Elliptic})). Moreover,  there exists a constant $\delta_0\in(0,1)$ (independent of $u$) such that $u\in C^{0,\delta_0}(\overline{\Omega})$, that is, $u$ is H\"older continuous over $\overline{\Omega}$, and there exists a constant $C>0$ (independent of $u$) such that
\begin{equation}\label{Holder-norm}
\|u\|_{_{C^{0,\delta_0}(\overline{\Omega})}}\,\leq\,C\left(\displaystyle\sum^2_{i=0}\|f_i\|_{_{p_i,\Omega}}+\|g\|_{_{q,\Gamma^{\infty}}}\right).
\end{equation}
\end{theorem}

\indent To achieve this latter strong result, one needs to establish several a priori global and local estimates, which together with an oscilation result enables us to reach the global regularity for the weak solution. For bounded Lipschitz domains and bounded coefficients, in \cite{NITTKA}, Nittka established this conclusion using a different method which relies strongly on the Lipschitz structure of the boundary (and thus is not valid for our case). For the case of the (symmetric) Laplacian with bounded lower-order coefficients, by using a completely different approach, in \cite{VELEZ2013-1}, global uniform continuity for weak solutions to the Robin problem was extended to hold for bounded extension domains whose boundaries are $d$-sets (see Definition \ref{Ahlfors}). We will use the motivations of this latter method to deduce fine regularity results for our case. In fact, when $\tau\in[1/2,\tau^{\ast})$, since $\Omega$ is locally uniform and $\Gamma^{\infty}$ is a $d$ with respect to $\mu(\cdot):=\mathcal{H}^d(\cdot)$, many of the results in \cite{VELEZ2013-1} can be easily adapted to hold for our problem, but this is no longer a valid argument in the critical case $\tau=\tau^{\ast}$. Therefore, a new idea needs to be developed in order to include this key case. Thankfully, the trace results in \cite{ACH08} together with the fact that $\Omega$ is a $W^{1,p}$-extension domain for lower values of $p$ allow us to find an alternative road to deduce the needed inequalities and estimates, leading to the fulfillment of fine regularity for the elliptic case.\\
\indent Turning now to the parabolic problem (\ref{1.01}) with initial condition and boundary conditions of type (\ref{1.02}), we first use the semigroup theory to transition into the homogeneous version of problem (\ref{1.01}). Given $A_{\mu}$ the unique operator associated with the bilinear form related to the variational formulation of the elliptic problem (as in (\ref{Existence-Elliptic})), let $\{T_{\mu}(t)\}_{t\geq 0}$ denote the corresponding $C_0$-semigroup of operators generated by $-A_{\mu}$ in $L^2(\Omega)$. Our second main result, which serves as a bridge between the elliptic problem and the desired inhomogeneous diffusion equation, reads as follows.

\begin{theorem}\label{feller-semi}
The part of the operator $A_{\mu}$ in $C(\overline{\Omega})$ generates a compact analytic $C_0$-semigroup $\{T^c_{\mu}(t)\}_{t\geq 0}$ of operators with angle $\pi/2$ over $C(\overline{\Omega})$. In particular, for each $u_0\in C(\overline{\Omega})$ and $t\geq0$, the function $u(t):=T_{\mu}(t)u_0\in C(\overline{\Omega})$ is the unique mild solution of the abstract Cauchy problem
$$\left\{
     \begin{array}{ll}
       u_t= A^c_{\mu}u\,\,\,\,\,\,\,\,\,\,\,\,\,\,\,\,\,\,\,\,\,\,\textrm{for}\,\,t\in(0,\infty);\\
       u(0)=u_0\,\,\,\,\,\,\,\,\,\,\,\,\,\,\,\,\,\,\,\,\,\,\textrm{in}\,\,\Omega
     \end{array}
   \right.$$
\end{theorem}

\indent There are two key ingredients used in order to reach to the conclusion of Theorem \ref{feller-semi}. First, the elliptic regularity theory obtained implies that $\{T_{\mu}(t)\}_{t\geq0}$ enjoys the {\it Feller property} in the sense that it maps $L^{\infty}(\Omega)$ into $C(\overline{\Omega})$ for each $t>0$ (in fact we deduce a stronger conclusion for this semigroup). Secondly, we show that for each positive time $t$, the set $T_{\mu}(t)(C(\overline{\Omega}))$ is dense over $C(\overline{\Omega})$. This is a standard step employed in order to obtain a $C_0$-semigroup over $C(\overline{\Omega})$, as in \cite{FUK-TOM96,NITTKA,WAR06}. However, this result has not been established for non-Lipschitz domains, up to the present time, due to the fact that one frequently needs to deal an approximation of the co-normal components of a smooth function, which in the case of irregular boundaries may not exist in the classical way, and may belong to a larger space. In fact, one does not know (up to the present time) how to deal with this strange element lying in a very big space, so one needs to find an alternative way to go around it. In the present paper, we deal with this situation for the first time, and use another approach motivated by some results in \cite{NITTKA2010,ZIE89} in order avoid dealing with a normal-type term.\\
\indent Finally, we adapt the approach introduced by Nittka \cite{NITTKA2014} to pass into our inhomogeneous boundary value problem (\ref{1.01}) with boundary conditions (\ref{1.02}). After establishing the unique solvability of this inhomogeneous mixed-type heat equation, we go through an extensive procedure of technical and complicated a priori estimates, adopting to our general non-smooth situation the ideas originally presented by Ladyzhenskaya, Solonnikov, and Ural'tseva \cite{LAD-SOL-URAL68}, and refined by Nittka \cite{NITTKA2014} (for bounded Lipschitz domains). We arrive then to the following general statement, which corresponds to the main result of this paper.

\begin{theorem}\label{Main-T2}
Given $T>0$ arbitrarily fixed, if $f\in L^2(0,T;L^2(\Omega))$, $g\in L^2(0,T;L^2_{\mu}(\Gamma^{\infty}))$, and $u_0\in L^2(\Omega)$, then the diffusion equation (\ref{1.01}) with mixed boundary conditions (\ref{1.02}) admits a unique weak solution $u\in C([0,T];L^2(\Omega))\cap L^2((0,T); \mathcal{V}_2(\Omega))$. Furthermore, under the condition (\ref{sharp-coercivity}), if $u_0\in C(\overline{\Omega})$,\, $f\in L^{\kappa_p}(0,T;L^p(\Omega))$, and $g\in L^{\kappa_q}(0,T;L^q_{\mu}(\Gamma^{\infty}))$, where
$\kappa_p,\,\kappa_q,\,p,\,q\in[2,\infty)$ are such that
$$\displaystyle\frac{1}{\kappa_p}+\displaystyle\frac{1}{p}<1\,\,\,\,\,\,\,\,\textrm{and}
\,\,\,\,\,\,\,\,\displaystyle\frac{1}{\kappa_q}+\displaystyle\frac{1}{2q}<\frac{1}{2},$$
then $u\in C([0,T];C(\overline{\Omega}))$ with
\begin{equation}\label{Main-Eq2}
\|u\|_{_{C([0,T];C(\overline{\Omega}))}}
\,\leq\,C^{\ast}\left(\|u_0\|_{_{C(\overline{\Omega})}}+\|f\|_{_{L^{\kappa_p}(0,T;L^{p}(\Omega))}}+
\|g\|_{_{L^{\kappa_q}(0,T;L^{q}_{\mu}(\Gamma^{\infty}))}}\right),
\end{equation}
for some constant $C^{\ast}>0$ (independent of $u$).
\end{theorem}

\begin{figure}[h!]
\begin{center}
  \subfloat{
    \includegraphics[width=0.45\textwidth]{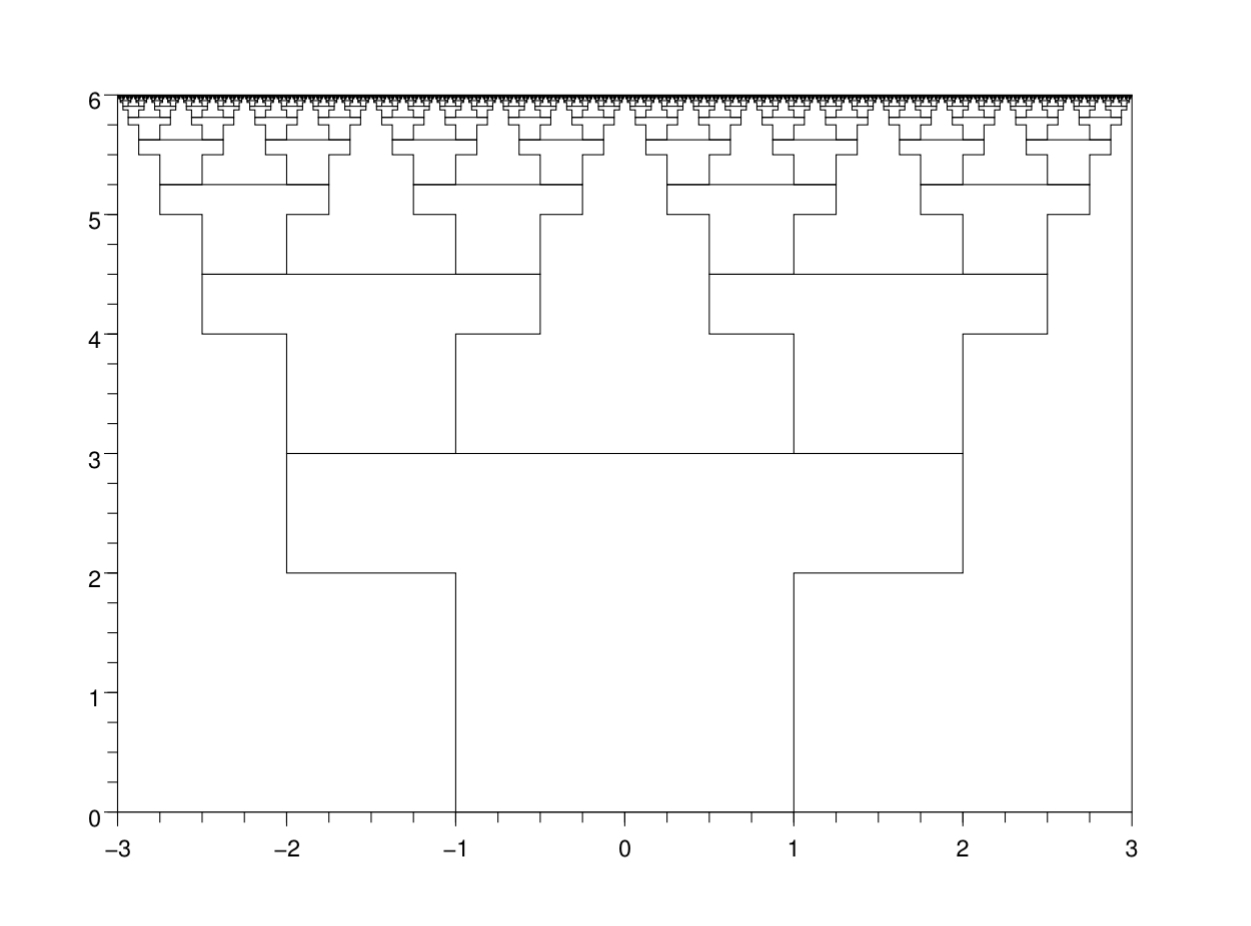}}
  \subfloat{
    \includegraphics[width=0.4\textwidth]{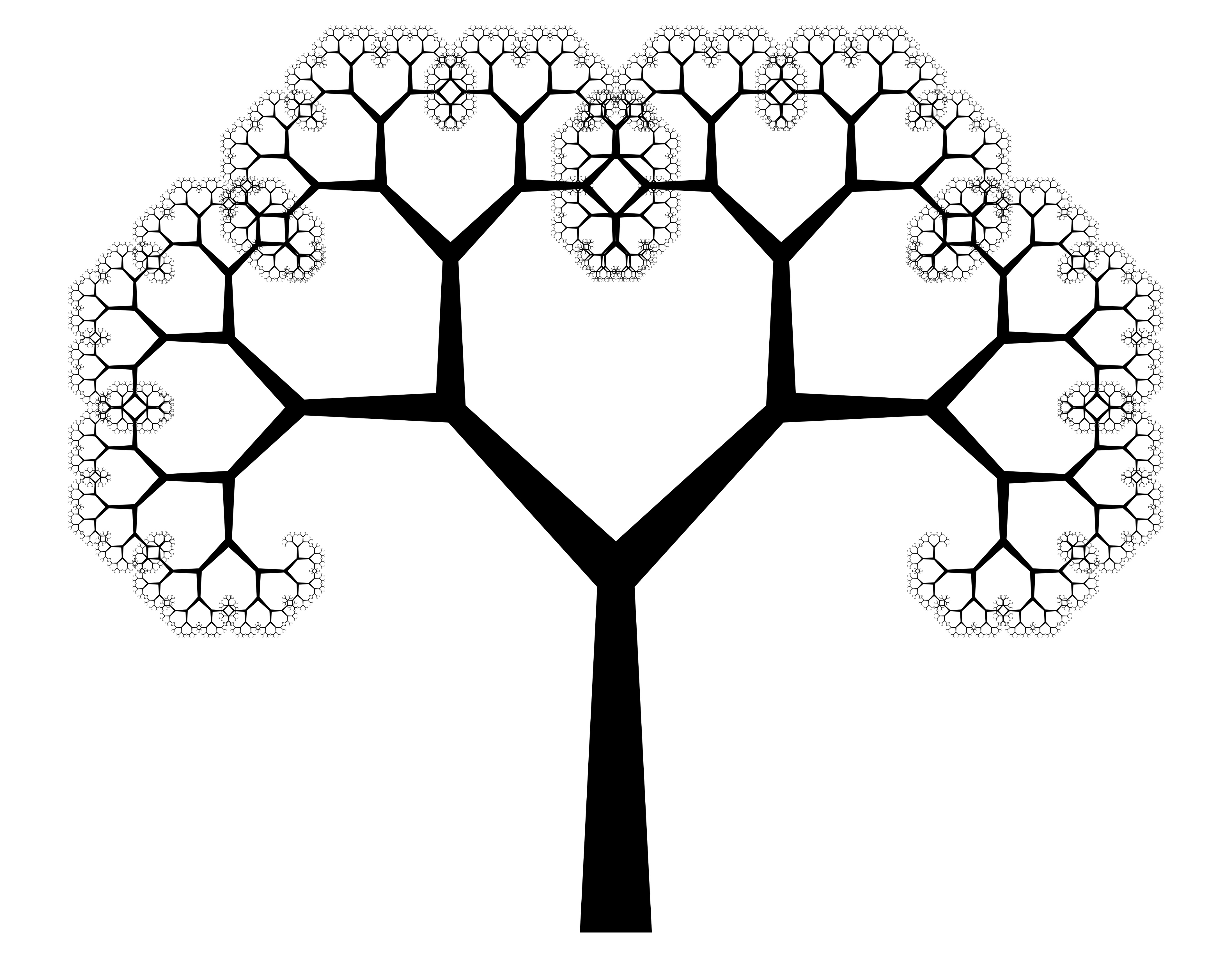}}
 \caption{Left figure: a T-shaped domain. \,Right figure: a variant of the ramified domain $\Omega$ involving different pre-fractal hexagons and a different angle in the critical case $\tau=\tau^{\ast}$.}
 \label{T-shaped-trees}
 \end{center}
\end{figure}

\begin{remark}\label{extension_results_appl}
The validity and applications of Theorem \ref{Holder-continuity}, Theorem \ref{feller-semi}, Theorem \ref{Main-T2}, extend to the following more general situations.
\begin{enumerate}
\item[(a)]\, The main results in this paper extend to more general classes of ramified and tree-shaped domains, as in \cite{ACH-DEH14,ACH-DEH12,ACH-SAB06-01,ACH-SAB06-02,ACH10,DEH16} (among others). Examples of these ramified sets can be seen in the Figure \ref{T-shaped-trees}. Both examples in Figure \ref{T-shaped-trees} are not extension domains (and also not locally uniform domains), but are $W^{1,p}$-extension domains for lower values of $p<2$.
\item[(b)]\, In addition to the mentioned application to the bronchial trees in the respiratory system of mammals, the model equation in this paper can be applied to other applications involving diffusion, such as diffusion processes on the circulatory system, gas exchange of plants with ramified-like structure in their processes of photosynthesis, as well as other areas of science and engineering (see for instance \cite{G-R-K05,K-K-M-P-G98,P-G-K15,S-Y-C-L20,X-Y-Y-Z06} and the references therein).
\end{enumerate}
\end{remark}

\subsection{Organization of the paper}\label{subsec1.5}

\indent Section \ref{sec2} concerns some basic notions and definitions about functions spaces, locally uniform domains, $d$-sets, and some known analytical intermediate results that will be applied throughout subsequent sections in the paper. Section \ref{sec3} introduces the concrete definition and main properties of the ramified domain. An extensive collection of properties regarding the ramified domain $\Omega$ and ramified boundary $\Gamma^{\infty}$ can be found in \cite{ACH10,ACH08,ACH07}, so we only state the main and relevant properties. We then establish sharp embedding, traces, and inequalities that are fundamental for the latter sections regarding the diffusion boundary value problem. It is important to point out here that some of these sharp inequalities are only known for extension domains (e.g. \cite[Theorem 3.3 and Theorem 3.6]{VELEZ2013-1}), which mean that in our domain, such results may not hold in the critical case $\tau=\tau^{\ast}$. Nevertheless, one can compensate the absence of the extension property by looking on the other properties of the domain, obtaining an analogous conclusion for our ramified domain, including the case $\tau=\tau^{\ast}$. Furthermore, we include a new approximation result for the Hausdorff measure of the ramified set $\Gamma^{\infty}$, namely, Theorem \ref{Main-T1}, whose proof is given at the end, in an appendix.\\
\indent The last two sections of the paper are devoted to the diffusion equation (\ref{1.01}) with Robin-type boundary conditions (\ref{1.02}). Section \ref{sec5} deals entirely with the elliptic version of this problem, and it is divided into three main parts. In the first part, we show the existence and uniqueness of a weak solution, under minimal assumptions. In the next part, we first derive global a priori estimate and $L^{\infty}$-estimates for weak solutions. Also, given $x_0\in\overline{\Omega}$ and $0<\rho<1$, the same estimates are established for a class of solutions (lying in a suitable space) solving problem (\ref{1.01}) with boundary conditions (\ref{1.02}) over $\Omega\cap B(x_0;\rho)$. Additional local measure and norm estimates over $\Omega\cap B(x_0;\rho)$ are given, with additional measure-type consequences. All these results are applied in the third part to produce an oscillation-type estimate, which, together with the global boundedness of weak solutions, allow us to obtain global H\"older continuity of weak solutions under minimal assumptions. Particulary, Theorem \ref{Holder-continuity} is established in this section, which is particularly new and original in the critical case $\tau=\tau^{\ast}$. In section \ref{sec6}, we consider the solvability and global regularity of the inhomogeneous diffusion equation (\ref{1.01}) with boundary conditions (\ref{1.02}), which we divide into four subsections. The regularity results in the previous section provide the tools to transition into the time-dependent homogeneous equation via Feller semigroups and the existence of a $C_0$-semigroup over $C(\overline{\Omega})$ related to problem (\ref{1.01}) with boundary conditions (\ref{1.02}). Theorem \ref{feller-semi} is established in subsection \ref{proof-main-2}. Subsection \ref{subsec6.2} brings all the classical theory of the inhomogeneous diffusion problem (\ref{1.01}), which runs as in \cite{NITTKA2014} with minor modifications. In particular, existence and uniqueness of weak solutions to problem (\ref{1.01}) fulfilling (\ref{1.02}) is proved here. Subsection \ref{subsec6.3} is all devoted to a priori estimates (under various assumptions on the initial conditions). In subsection \ref{subsec6.4}, all these results are combined in order to establish Theorem \ref{Main-T2}, which is the third and key main result of the paper.\\
\indent At the end, we add an appendix, in which we devote our attention in the establishment of Theorem \ref{Main-T1}. Indeed, after settling the proper language and stating some key known results from Jia \cite{JIA07-1,JIA07-2}, Theorem \ref{Main-T1} is fully proved. The proof is motivated by the procedures used in Jia \cite{JIA07-1,JIA07-2} together with some refinements made in \cite{FERR-VELEZ18-1}. This result is somewhat unrelated with respect to the resolution of the diffusion problem (\ref{1.01}) with boundary conditions (\ref{1.02}), but since the conclusion in Theorem \ref{Main-T1} is new and original, we have decided to include its proof in the appendix.

\section{Basic preliminaries}\label{sec2}

\subsection{Function spaces}\label{subsec2.1}

\indent In this section we review some fundamental properties of function spaces, and present some key definitions that will be valuable later on.\\
\indent As usual, $L^p_{\mu}(E):=L^p(E;d\mu)$ denotes the $L^p$-space on $E$ with respect to the measure $\mu$; if $\mu$ is the $N$-dimensional Lebesgue measure, then we write $L^p_{\mu}(E)=L^p(E)$. Also, given $\Omega\subseteq\mathbb{R\!}^{\,2}$ a bounded domain with boundary $\Gamma:=\partial\Omega$, by $W^{1,p}(\Omega)$ we mean the well-known $L^p$-based Sobolev space, and we write $H^1(\Omega):=W^{1,2}(\Omega)$. At times, we will deal with Besov spaces $H^{d/2}_{\mu}(\Gamma):=\left\{u\in L^2_{\mu}(\Gamma)\mid\mathcal{N}^2_{_{d/2}}(u,\Gamma,\mu)<\infty\right\}$, for $0<d<2$, where
$$\mathcal{N}^2_{_{d/2}}(u,\Gamma,\mu):=\left(\displaystyle\int_{\Gamma}\int_{\Gamma}\left[\frac{|u(x)-u(y)|}
{|x-y|^{d}}\right]^2\,d\mu_xd\mu_y\right)^{\frac{1}{2}},$$
and its corresponding dual will be denoted by $H^{d/2}_{\mu}(\Gamma)^{\ast}$. Also, for $r,\,s\in[1,\infty)$, or $r=s=\infty$, we will sometimes refer to the Banach Space $\mathbb{X\!}^{\,r,s}(\Omega;\Gamma):=L^r(\Omega)\times L^s_{\mu}(\Gamma)$ endowed with norm
$$|\|(f,g)\||_{_{\mathbb{X\!}^{\,r,s}(\Omega;\Gamma)}}:=\|f\|_{_{r,\Omega}}+\|g\|_{_{s,\Gamma}}\,,
\,\,\,\,\,\,\textrm{if}\,\,\,r,\,s\in[1,\infty),\,\,\,\,\,\,\textrm{and}\,\,\,\,\,
|\|(f,g)\||_{_{\mathbb{X\!}^{\,\infty,\infty}(\Omega;\Gamma)}}:=
\max\left\{\|f\|_{_{\infty,\Omega}},\,\|g\|_{_{\infty,\Gamma}}\right\}.$$
If $r=s$, then we will write $\mathbb{X\!}^{\,r}(\Omega;\Gamma)=\mathbb{X\!}^{\,r,r}(\Omega;\Gamma)$. Also if $u\in L^1(\Omega)$ is such
that $u|_{_{\Gamma}}\in L^1_{\mu}(\Gamma)$, we write $\mathbf{u}:=(u,u|_{_{\Gamma}})$. For more information regarding these spaces, refer to \cite{ADA,MAZ,NECAS,TAR07,ZIE89} (among many others).\\
\indent Next, we give the following three standard definitions (e.g. \cite{BIE09,DLL06,HAJLASZ-KOS-TUO08,JON,JO-WAL}).

\begin{definition}\label{ext-domain}
Let $p\in[1,\infty]$. A domain
$\Omega\subseteq\mathbb{R\!}^N$ is called a $W^{1,p}$-{\bf extension domain}, if there exists a bounded linear operator
$S:W^{1,p}(\Omega)\rightarrow W^{1,p}(\mathbb{R\!}^N)$ such that $Su=u$\, a.e. on $\Omega$. If $p=2$ and the latter property holds for
$\Omega$, we will call it simply an {\it extension domain}.
\end{definition}

\begin{definition}\label{ep-del-dom}
A domain $\Omega\subseteq\mathbb{R\!}^N$ is said to be an\, {\bf ($\mathbf{\epsilon,\,\delta}$)-domain},
if there exists $\delta\in (0,+\infty]$ and there exists $\epsilon\in (0,1]$, such that for each two points $x,\,y\in\Omega$
with $|x-y|\leq\delta$, there exists a continuous rectifiable curve $\gamma :[0,t]\rightarrow\Omega$ such that $\gamma(0)=x$ and
$\gamma(t)=y$, with the following properties:
\begin{enumerate}
\item[(a)]\,\,\,\,$\ell(\{\gamma\})\leq\frac{1}{\epsilon}|x-y|$.
\item[(b)]\,\,\, $\textrm{dist}(z,\partial\Omega)\geq\,\epsilon\,\min
\{|x-z|,|y-z|\},\,\,\,\textrm{for all}\,\,z\in\{\gamma\}$.
\end{enumerate}
These sets are also referred as \textbf{locally uniform domains} (e.g. \cite{ROGERS06}).
\end{definition}

\indent In a classical paper by Jones \cite{JON}, it was shown that every $(\epsilon,\delta)$-domain is a $W^{1,p}$-extension domain for all $p\in[1,\infty]$, and the converse of this result holds for any finitely connected open set in $\mathbb{R\!}^{\,2}$.

\begin{definition}\label{Ahlfors}
Let $d\in(0,N)$ and $\mu$ a Borel measure supported on a bounded set $F\subseteq\mathbb{R\!}^N$.
Then $\mu$ is said to be a {\bf $d$-Ahlfors measure}, if there exist constants $M_1,\,M_2,\,R_0>0$ such that
\begin{equation}\label{2.01}
M_1r^2\,\leq\,\mu(F\cap B(x;r))\,\leq\,M_2r^d,\indent\,\,\,\textrm{for all}\,\,\,0<r<R_0\,\,\,\textrm{and}\,\,\,x\in F,
\end{equation}
where $B(x;r)$ denotes the ball of radius $r$ centered at $x\in F$. In this case the set $F\subseteq\mathbb{R\!}^N$ will be called a
{\bf $d$-set} (with respect to the measure $\mu$; see for instance \cite{JO-WAL}).
\end{definition}

\begin{remark}\label{facts_extension_domains}
When $\Omega\subseteq\mathbb{R\!}^N$ is bounded, it is obvious that the inequality
$|B(x,r)\cap\Omega|\leq cr^N$ is valid for some constant $c>0$ and for every $x\in\Omega$ and $r>0$.
If in addition $\Omega$ is a bounded $W^{1,p}$-extension domain, then it is well-known that $\Omega$ satisfies
the so called {\bf measure density condition}, that is, the inequality
$$|B(x,r)\cap\Omega|\,\geq\,c' \,r^N$$
holds some constant $c' >0$, and for all $x\in\Omega$ and $r\in(0,1]$ (e.g. \cite{SHV07}), where $|\cdot|$ denotes the $N$-dimensional Lebesgue measure. Thus, one sees that in this
case the set $\Omega$ is a $N$-set with respect tho the Lebesgue measure $|\cdot|$. Moreover, if $\Omega$ is a
$W^{1,p}$-extension domain, then $|\partial\Omega|=0$.
\end{remark}

\subsection{Known analytical results}\label{subsec2.2}

\indent The following results of analytical nature will be important in the development of regularity results.

\begin{lemma}\label{lemma1}\,(see \cite{LAD-SOL-URAL68})\,
Let $\{y_n\}_{n\geq0},\,\{z_n\}_{n\geq0}$ be sequences of nonnegative real numbers, and let $c$,\, $b$,\, $\varepsilon$, and $\delta$ be positive constants
with $b\geq1$, such that
$$y_{n+1}\,\leq\,cb^n\left(y^{1+\delta}_n+z^{1+\varepsilon}_ny^{\delta}_n\right)\,\,\,\,\,\,\,\,\textrm{and}
\,\,\,\,\,\,\,\,z_{n+1}\,\leq\,cb^n\left(y_n+z^{1+\varepsilon}_n\right),\,\,\,\,\,\,\,\,\textrm{for all}\,\,n\in\mathbb{N\!}\,.$$
Define $$l:=\min\left\{\delta,\,\displaystyle\frac{\varepsilon}{1+\varepsilon}\right\}\,\,\,\,\,\,\,\,\textrm{and}
\,\,\,\,\,\,\,\,\omega:=\min\left\{(2c)^{^{-\frac{1}{\delta}}}b^{^{-\frac{1}{\delta l}}},\,
(2c)^{^{-\frac{1+\varepsilon}{\varepsilon}}}b^{^{-\frac{1}{\varepsilon l}}}\right\},$$ and assume that
$$y_0\,\leq\,\omega\,\,\,\,\,\,\,\,\textrm{and}\,\,\,\,\,\,\,\,z_0\,\leq\,\omega^{1/(1+\varepsilon)}.$$
Then $$y_n\,\leq\,\omega b^{^{-\frac{n}{b}}}\,\,\,\,\,\,\,\,\textrm{and}\,\,\,\,\,\,\,\,z_n\,\leq\,\left(\omega b^{^{-\frac{n}{b}}}\right)^{^{\frac{1}{1+\varepsilon}}},
\,\,\,\,\,\,\,\textrm{for every integer}\,\,n\geq0.$$
\end{lemma}
\indent\\
\begin{lemma}\label{lemma2}\,(see \cite{STAM68})\,
{\it Let} $\varphi\,=\,\varphi(t,r)$ {\it be a nonnegative
function on a half-strip} $\{t\geq k_0\geq 0\}\times\{0\leq r<R_0\}$ {\it such that:}\\[2ex]
\indent (a)\,\,{\it$\varphi(\cdot,r)$ is non-increasing for every fixed $r$,}\\[2ex]
\indent (b)\,\,{\it$\varphi(h,\cdot)$ is non-decreasing for every fixed $h$,}\\[2ex]
{\it and such that there exist $c,\tau,\gamma>0$, and} $\delta>1$ {\it with}
$$\varphi(h,r)\leq c\,(h-k)^{-\tau}(R-r)^{-\gamma}\varphi(k,R)^{\delta},$$
{\it for all} $h>k\geq k_0,\,r<R<R_0$. {\it Then for $\xi\in(0,1)$ arbitrary we have that}
$$\varphi(k_0+\varsigma,(1-\xi)R_0)=0,$$
{\it where} $$\varsigma^{\tau}=c((1-\xi)R_0)^{-\gamma}\varphi(k_0,R_0)^{\delta-1}2^{^{\delta\left
(\frac{\tau+\gamma}{\delta-1}\right)}}.$$
\end{lemma}
\indent\\
\begin{lemma}\label{lemma3}\,(see \cite{STAM68})\,
{\it Let} $\varphi\,=\,\varphi(t)$ {\it be a nonnegative,
non-increasing function on a closed interval} $\{k_0\leq t\leq M\}$, {\it and assume that there exist}
$c,\,\tau,\,\delta>0$ {\it with}
$$(h-k)^{\tau}\varphi(h)^{\delta}\leq c\,(M-k)^{\tau}(\varphi(k)-\varphi(h)),$$
{\it for all} $M>h>k\geq k_0$. {\it Then}
$\displaystyle\lim_{^{h\rightarrow M}}\varphi(h)=0.$\\
\end{lemma}

\begin{lemma}\label{lemma4}\,(see \cite{GIL-TRU})\,
{\it Let $\Phi\,=\,\Phi(t)$ be a nondecreasing function defined on \,$\{0<t\leq R\}$.
If there are constants $\eta,\,\varsigma\in(0,1),\,\tau>0$,  such that}
$$\Phi(\varsigma\rho)\,\leq\,\eta\,\Phi(\rho)+\psi(\rho)\indent\textit{for all}\,\,\,\rho\in(0,R],$$
{\it where $\psi(\cdot)$ is another nondecreasing function over \,$\{0<t\leq R\}$,
then for every $\delta\in(0,1)$ and $\rho\in(0,R]$, we have}
$$\Phi(\rho)\,\leq\,C\left[\left(\frac{\rho}{R}\right)^{\vartheta}\Phi(R)+\psi(\rho^{\delta}R^{1-\delta})\right],$$
for some constants $C=C(\eta,\varsigma)>0$ and $\vartheta=\vartheta(\eta,\varsigma,\theta)>0$.
\end{lemma}

\section{Ramified domain: construction \& properties}\label{sec3}

\subsection{Basic construction of the ramified domains}\label{subsec3.1}

\indent We consider a class of domains with ramified fractal-like boundaries. Such domains are of particular interest in applications to
bronchial trees (see for instance \cite{ACH08} among others).\\
\indent In fact, we consider the following particular domain $\Omega$ investigated in
\cite{ACH08} (among others).
Let $\{G_1,G_2\}$ denote a set consisting of a pair of similitudes on $\mathbb{R\!}^{\,2}$ defined by
$$G_i(x_1,x_2):=\left(\begin{array}{lcl}
(-1)^i\left[1-\displaystyle\frac{\tau}{\sqrt{2}}\right]+\displaystyle\frac{\tau}{\sqrt{2}}[x_1+(-1)^ix_2]\\
\indent\\
1+\displaystyle\frac{\tau}{\sqrt{2}}+\displaystyle\frac{\tau}{\sqrt{2}}[x_2+(-1)^{i+1}x_1]\\
\end{array}
\right)$$
\noindent ($i\in\{1,2\}$) for $\tau\leq\tau^{\ast}\simeq 0.593465$. In fact, $\tau^{\ast}$ is the solution of the equation
\begin{equation}\label{tau-eq}
2\sqrt{2}\,\tau^5+2\tau^4+2\tau^2+\sqrt{2}\,\tau-2=0
\end{equation}
(see \cite{ACH08}). Let $\mathfrak{A}_n$ be the set containing all the maps from
$\{1,\ldots,n\}$ to $\{1,2\}$. Then for $\sigma\in\mathfrak{A}_n$, we define the affine map
$$\mathcal{M}_{\sigma}(G_1,G_2):=G_{\sigma(1)}\circ\cdot\cdot\cdot\circ G_{\sigma(n)}.$$
Then, letting $V$ be the interior of the convex hull of the
points $P_1=(-1,0)$,\, $P_2=(1,0)$,\, $P_3=G_1(P_1)$,\, $P_4=G_2(P_2)$,\, $P_5=G_1(P_2)$, and $P_6=G_2(P_1)$,
we define
$$\Omega:=\textrm{Interior}\left(\overline{V}\cup\left(\displaystyle\bigcup^{\infty}_{n=1}\displaystyle\bigcup_{\sigma\in\mathfrak{A}_n}
\mathcal{M}_{\sigma}(G_1,G_2)(\overline{V})\right)\right).$$
The control $\tau\leq\tau^{\ast}$ ensures that the sets $\mathcal{M}_{\sigma}(G_1,G_2)$ will not overlap.
It follows that the ramified boundary $\Gamma^{\infty}$ ($i\in\{1,2\}$) corresponds to the unique compact self-similar set in $\mathbb{R\!}^{\,2}$ fulfilling
$$\Gamma^{\infty}=G_1(\Gamma^{\infty})\cup G_2(\Gamma^{\infty}).$$

For more information and properties of these sets, refer to \cite{ACH08} (among others). In particular, one can summarize
several of the main properties in the following result.

\begin{theorem}\label{prop-ramified}\,(see \cite{ACH08})\,The following properties are valid for the set $\Omega$.
\begin{enumerate}
\item[(a)]\,\, $\Omega$ has the open set condition.
\item[(b)]\,\, $\Omega$ is an $(\epsilon,\delta)$-domain whenever $1/2\leq\tau<\tau^{\ast}$.
\item[(c)]\,\, $\Omega$ is not an extension domain when $\tau=\tau^{\ast}$, but it is a $W^{1,q}$- extension domain for all $q\in[1,2)$. In particular, $\Omega$ is not an $(\epsilon,\delta)$-domain when $\tau=\tau^{\ast}$.
\item[(d)]\,\, The ramified boundary $\Gamma^{\infty}$ is a $d$-set with respect to the self-similar measure supported on $\Gamma^{\infty}$, where
$$d:=-\log(2)/\log(\tau).$$
\item[(e)]\,\, $\Omega$ is a $2$-set with respect to the $2$-dimensional Lebesgue measure $|\cdot|$.
\end{enumerate}
\end{theorem}

\begin{figure}[h!]
\begin{center}
  \subfloat{
   \label{1_1}
    \includegraphics[width=0.4\textwidth]{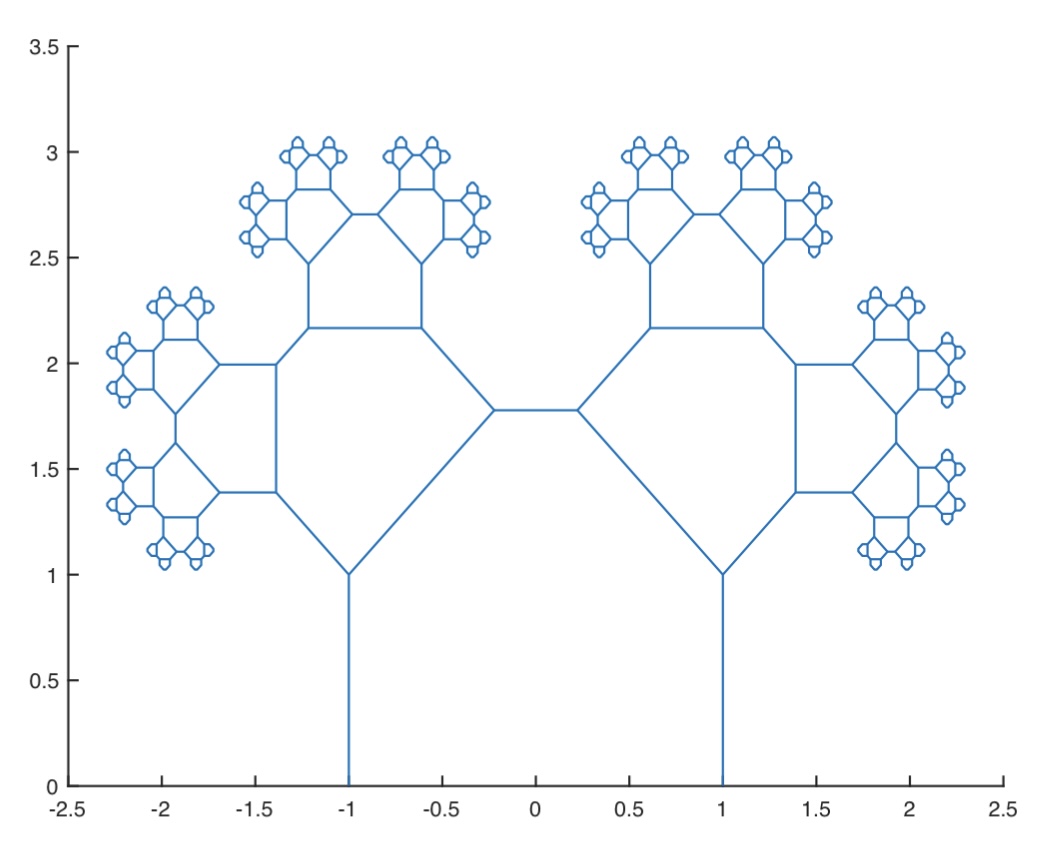}}
  \subfloat{
    \includegraphics[width=0.4\textwidth]{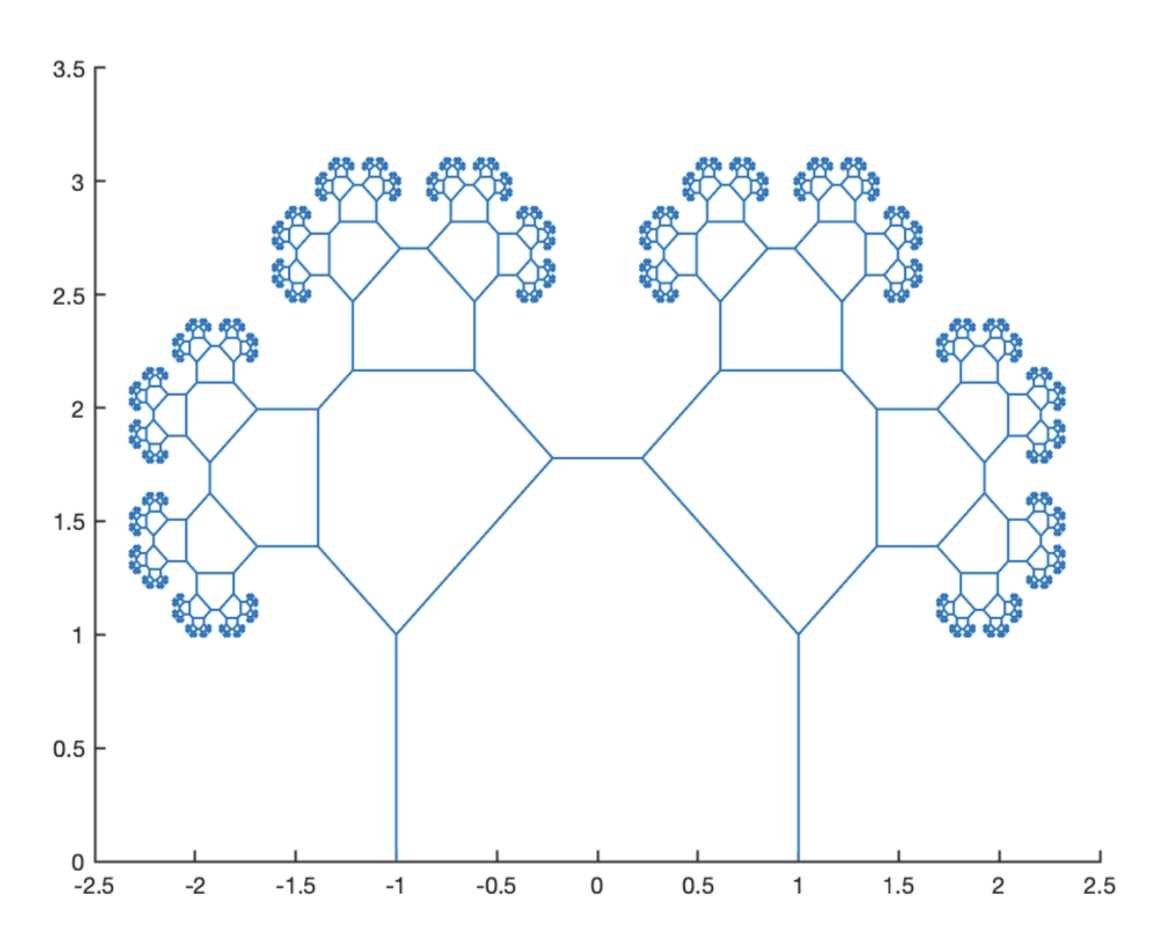}}
 \caption{The pre-fractals: $\Omega^{(5)}$ (left), and $\Omega^{(9)}$ (right).}
 \label{pre-fractals-pictures}
 \end{center}
\end{figure}

\indent For each $m\in\mathbb{N\!}_{\,0}:=\{0\}\cup\mathbb{N\!}$\,, we denote by
\begin{equation}
    \label{Omega-ramified}
\Omega^{(m)}:=\textrm{Interior}\left(\overline{V}\cup\left(\displaystyle\bigcup^{m}_{n=0}\displaystyle\bigcup_{\sigma\in\mathfrak{A}_n}
\mathcal{M}_{\sigma}(G_1,G_2)(\overline{V})\right)\right)
\end{equation}
the pre-fractal domain obtained by stopping the construction of the ramified domain $\Omega$ at the step $m+1$ (see figure \ref{pre-fractals-pictures}).
It follows that $\{\Omega^{(m)}\}_{_{m\geq0}}$ defines an increasing sequence of bounded Lipschitz domains with $\Omega=\displaystyle\bigcup^{\infty}_{m=0}\Omega^{(m)}$.\\

\indent There is a procedure useful to approximate accurately the ``length" of the ramified boundary $\Gamma^{\infty}$. It is known that the notion of classical length in $\Gamma^{\infty}$ is infinite (and thus meaningless), so one needs a correct way to measure the length of a fractal ramified set such as $\Gamma^{\infty}$. This is done by measuring $\Gamma^{\infty}$ with respect to its corresponding Hausdorff measure $\mu(\cdot):=\mathcal{H}^d(\cdot)$ for $d$ its corresponding Hausdorff dimension. The computation and estimation of the Hausdorff dimension and measure of a fractal set is an important problem in fractal geometry.\\
\indent For this part, we build a rectangle with a base of length 2 and position it on a plane centered at the origin, so that the base aligns with the interval $[-1, 1]$. The height of the rectangle can vary. Then, we take transversal cuts at the upper corners with a common angle, for example, 45 degrees, and these cuts have a size of $2\tau$ (see figure \ref{omega2}). From here, the above statement is contained in the following new original result. Since the proof is technical and kind of unrelated with the main subject of the paper, we will provide a complete derivation of the result at the end, in the Appendix.

\begin{theorem}\label{Main-T1}
There exists a monotone decreasing sequence $\{a_n\}\subseteq(0,1)$ and a monotone increasing sequence $\{b_n\}\subseteq(0,1)$ such that the $d$-dimensional Hausdorff measure of the fractal boundary $\Gamma^{\infty}$ of the ramified domain $\Omega$ satisfies the estimation
\begin{equation}\label{Gen_Est_Hausdorff}
0.7582916255\,\leq\,b_1\,\leq\,b_2\,\leq\cdot\cdot\cdot\leq\,\displaystyle\sup_{n\in\mathbb{N\!}}\{b_n\}=\mathcal{H}^d(\Gamma^{\infty})=\displaystyle\inf_{n\in\mathbb{N\!}}\{a_n\}\,\leq\cdot\cdot\cdot\leq\,a_2\,\leq\,a_1\,\leq\,0.894684280597,
\end{equation}
for  $d:=-\log_{\tau}2$ the Hausdorff dimension of $\Gamma^{\infty}$ (the above estimate taken uniformly with respect to any $\tau\in[1/2,\tau^{\ast}]$). Furthermore,  
$$b_n=a_n\exp\left(-\frac{2\sqrt{2}\,d}{
\hat{\delta}_{\tau}(1-\tau)(\sqrt{2}+3\tau)}\tau^{n}\right), \,\,\,\,\,\textrm{for}\,\,\,\hat{\delta}_{\tau}=\min_{\delta \in \Delta}\{\delta\}$$
where $\Delta := \{\delta \in (0, 1): \text{diam}(\bigcup_{i=1}^{k}\Delta_{i}) \geq \delta \mathfrak{H}\}$, with $\mathfrak{H}:=\sup_{(x_1, x_2)\in \Omega}x_{2}=\frac{1+\frac{3\tau}{\sqrt{2}}}{1-\tau^2}$ denoting the notion of height for $\Omega$. Finally, for each fixed $\tau\in[1/2,\tau^{\ast}]$, the lower and upper bounds in (\ref{Gen_Est_Hausdorff}) can be improved. 
\end{theorem}

\indent One can easily compute and obtain the last assertion in Theorem \ref{Main-T1}. In fact, if one takes $\tau=\tau^{\ast}$, the Hausdorff measure of $\Gamma^{\infty}$ can be computationally evaluated to obtain that
\begin{equation}\label{a_1-computation}
0.83383674565 \approx a_{1}\exp\left(-\frac{2\sqrt{2}d}{\delta_{\tau^*}(1-\tau)(\sqrt{2}+3\tau)}\tau^{n}\right)\,\leq\, \mathcal{H}^{d}(\Gamma^{\infty})\,\leq \,a_{1}\approx 0.894684280597,
\end{equation}
giving forth a fairly good approximation for $\mathcal{H}^{d}(\Gamma^{\infty})$. Also, note that for every $\tau_{1}, \tau_{2} \in \left[\frac{1}{2}, \tau^*\right]$, $\hat{\delta}_{\tau_2}\leq \hat{\delta}_{\tau_1}$ holds, and thus if $\tau=\frac{1}{2}$, then we deduce that
\begin{equation}\label{a_1-1/2}
0.7582916255\approx b_{1}\leq \mathcal{H}^{d}(\Gamma^\infty)\leq a_{1}\approx 0.8535533906
\end{equation}
As the sequence $\{a_{n}\}$ is monotonically decreasing, within the range  of the parameter $\tau\in \left[\frac{1}{2}, \tau^*\right]$, one can verify that the measure approximating interval for $\Gamma^\infty$ runs to the left whenever $\tau$ takes a smaller value. Consequently, for the sequence $\{b_{n}\}$ the value closest to the measure to be found will occur just when $\tau=\tau^*$. More precisely, $$\inf_{n\in\mathbb{N}, \tau \leq \tau^*}\{a_{n}\}=\inf_{n\in\mathbb{N}, \tau=\frac{1}{2}}\{a_{n}\}\,\,\,\,\,\,\,\,\textrm{and}\,\,\,\,\,\,\,\,\sup_{n\in\mathbb{N}, \tau \leq \tau^*}\{b_{n}\}=\sup_{n\in\mathbb{N}, \tau=\tau^*}\{a_{n}\}.$$

\indent Since the ramified boundary $\Gamma^{\infty}$ is self-similar, we can follow a procedure similar to the one developed by Jia \cite{JIA07-2} to obtain accurate upper and lower bounds for $\mathcal{H}^d(\Gamma^{\infty})$ as in Theorem \ref{Main-T1}. The central tools of the method introduced in \cite{JIA07-2} has been used to estimate the Hausdorff measure of some classical self-similar fractals, such as the Sierpinski gasket and the Koch curve (e.g. \cite{JIA07-1,JIA07-2}), and recently to a family of $3$-dimensional Koch-type fractal crystals introduced in \cite{FERR-VELEZ18-1}.

\subsection{Embedding, trace, and sharp inequalities over ramified domains}\label{subsec3.2}

\indent We now present the key embedding, trace results, as well as important Poincar\'e-type inequalities over the ramified set $\Omega$ described in the previous subsection. For a set $E\subseteq\mathbb{R\!}^{\,2}$ and for $B(x_0,\rho)$ a ball centered in $x_0\in\mathbb{R\!}^{\,2}$ with radius $\rho>0$, we write $$E(x_0,\rho):=E\cap B(x_0,\rho).$$

\begin{theorem}\label{embeddings-trace}
Given $\Omega\subseteq\mathbb{R\!}^{\,2}$ the domain with ramified fractal boundary $\Gamma^{\infty}$ described in the previous subsection, let $r,\,s\in[1,\infty)$ be arbitrarily fixed, take any $x_0\in\overline{\Omega}$, and let $\rho>0$. Then the following hold:
\begin{enumerate}
\item[(a)]\,\,\,There exists a continuous and compact embedding
$H^1(\Omega)\hookrightarrow L^{r}(\Omega)$ and a constant $c_1>0$ such that
\begin{equation}\label{2.03}
\|u\|_{_{r,\Omega}}\,\leq\,c_1\|u\|_{_{H^{1}(\Omega)}},
\indent\,\,\,\,\,\textrm{for every}\,\,\,u\in H^{1}(\Omega).
\end{equation}
\item[(b)]\,\,\,There exists a linear compact mapping $H^{1}(\Omega)\hookrightarrow L^{s}_{\mu}(\Gamma^{\infty})$ and a constant $c_2>0$ such that
\begin{equation}\label{2.04}
\|u\|_{_{s,\Gamma^{\infty}}}\,\leq\,c_2\|u\|_{_{H^{1}(\Omega)}},
\indent\,\,\,\,\,\textrm{for all}\,\,\,u\in H^{1}(\Omega).
\end{equation}
\item[(c)]\,\,\,There exist constants $\rho_0>0$ and $c_3>0$ (independent of $\rho$) such that
\begin{equation}\label{2.05}
\|u\|^2_{_{r,\Omega(x_0,\rho)}}\,\leq\,c_3\,\|\nabla u\|^2_{_{2,\Omega(x_0,\rho)}},
\end{equation}
for all $0<\rho<\rho_0$ and for every $u\in H^{1}(\Omega)$ such that $u=0$ over $E$ for some measurable subset $E$ of $\overline{\Omega}(x_0,\rho)$
such that $|E|\geq\kappa\,|\Omega(x_0,\rho)|$.
\item[(d)]\,\,\,There exist constants $\rho_0>0$ and $c_4>0$ (independent of $\rho$) such that
\begin{equation}\label{2.06}
\|u\|^2_{_{s,\Gamma(x_0,\rho)}}\,\leq\,c_4\,\|\nabla u\|^2_{_{2,\Omega(x_0,\rho)}},
\end{equation}
for all $0<\rho<\rho_0$ and for every $u\in H^{1}(\Omega)$ such that $u=0$ over $E$ for some measurable subset $E$ of $\overline{\Omega}(x_0,\rho)$
such that $|E|\geq\kappa\,|\Omega(x_0,\rho)|$.
\item[(e)]\,\,\,There exists a constant $c_5=c(N,|E|)>0$ such that
\begin{equation}\label{2.07}
\|u\|_{_{H^{1}(\Omega)}}\,\leq\,c_5\|\nabla u\|_{_{2,\Omega}},
\end{equation}
for every $u\in H^{1}(\Omega)$ that integrates to zero over a measurable set $E\subseteq\Omega$ with $|E|>0$.\\
\end{enumerate}
\end{theorem}

\begin{proof}
If $1/2\leq\tau<\tau^{\ast}$, then by virtue of Theorem \ref{prop-ramified} together with the results in \cite{BIE09,VELEZ2013-1,ZIE89}, one deduces all the conclusions of the theorem. For the case $\tau=\tau^{\ast}$, this is no longer the case. However, by virtue of \cite{ACH08}, we can establish the desired results over a simplified geometrical fractal-type domain $\widetilde{\Omega}$ described in \cite[proof of Theorem 4]{ACH08}, whose structure behaves in a similar way as in \cite{ACH-SAB06-01}. Therefore without loss of generality, we will still write $\widetilde{\Omega}=\Omega$. Then, if $r\in[2,\infty)$ is fixed, then one can find $r_0\in[1,2)$ such that $r=2r_0(2-r_0)^{-1}$. Then, choosing $\varepsilon_0>0$ small enough such that $1<r_{\varepsilon_0}:=r_0+\varepsilon_0<2$, from \cite[Proposition 3.3]{ACH-SAB06-01}, one get the following continuous mappings:
\begin{equation}\label{mult-embeddings}
H^1(\Omega)\hookrightarrow W^{1,r_{\varepsilon_0}}(\Omega)\stackrel{\textrm{compact}}{\hookrightarrow} L^{^{\frac{2r_0}{2-r_0}}}(\Omega)\,=\,L^r(\Omega),
\end{equation}
which yields statement (a). Point (b) follows directly from \cite[Theorem 4]{ACH08}. To prove (c), we recall that by virtue of Theorem \ref{prop-ramified}(c), it follows that $\Omega$ is a $W^{1,r}$-extension domain for each $r\in[1,2)$. Then, under all the assumptions in statement (c), following \cite[proof of Theorem 3.3]{VELEZ2013-1}, from (\ref{mult-embeddings}) and H\"older's inequality, we obtain that
$$\|u\|^2_{_{r,\Omega(x_0,\rho)}}\,\leq\,c\,\|\nabla u\|^2_{_{r_{\varepsilon_0},\Omega(x_0,\rho)}}\,\leq\,c'\,\|\nabla u\|^2_{_{2,\Omega(x_0,\rho)}},$$
for some constants $c,\,c'>0$. This completes the proof of part (c). Now, invoking \cite[Lemma 10]{ACH08}, we get the existence of a constant $s_0\in[1,2)$ such that the the trace mapping $W^{1,s_0}(\Omega)\hookrightarrow L^s_{\mu}(\Gamma^{\infty})$ is bounded. From here, proceeding as in \cite[proof of Theorem 3.3 and Theorem 3.6]{VELEZ2013-1}, we conclude that
$$\|u\|^2_{_{s,\Gamma^{\infty}(x_0,\rho)}}\,\leq\,\eta\,\|\nabla u\|^2_{_{s_0,\Omega(x_0,\rho)}}\,\leq\,\eta'\,\|\nabla u\|^2_{_{2,\Omega(x_0,\rho)}},$$ for some constants $\eta,\,\eta'>0$, which yields (d). Finally, statement (e) follows at once from a combination of the results in \cite{ACH08,ZIE89}, completing the proof.
\end{proof}

\begin{remark}\label{Sharp-epsilon-ineq}
Given $\Omega$ the domain in Theorem \ref{embeddings-trace}, in views of \cite[Lemma 2.4.7]{NITTKA2010}, for any $\epsilon>0$, there exists constants $C_{\epsilon},\,C'_{\epsilon}>0 $ such that
\begin{equation}\label{epsilon-interior}
\|u\|_{_{r,\Omega}}\,\leq\,\epsilon\|\nabla u\|_{_{2,\Omega}}+C_{\epsilon}\|u\|_{_{2,\Omega}}
\end{equation}
and
\begin{equation}\label{epsilon-trace}
\|u\|_{_{s,\Gamma^{\infty}}}\,\leq\,\epsilon\|\nabla u\|_{_{2,\Omega}}+C'_{\epsilon}\|u\|_{_{2,\Omega}},
\end{equation}
for each $u\in H^{1}(\Omega)$, whenever $1\leq r<\infty$ and $1\leq s<\infty$.
\end{remark}

\begin{remark}\label{Equivalent-norm}
If $\Omega\subseteq\mathbb{R\!}^{\,2}$ is the set in Theorem \ref{embeddings-trace}, 
then it is clear that all the conclusions in Theorem \ref{embeddings-trace} and Remark \ref{Sharp-epsilon-ineq} hold if one replaces $H^1(\Omega)$ with the closed subspace
\begin{equation}\label{V-space}
 \mathcal{V}_2(\Omega):=\left\{u\in H^1(\Omega)\mid u|_{_{\Gamma\setminus\Gamma^{\infty}}}=0\right\},
\end{equation}
with equivalent (and even refined) norms. Moreover,
from (\ref{2.04}) together with
\cite[Theorem 1]{ACH08}, it follows that the norm
$$|\|u\||^2_{_{\widetilde{\mathcal{V}}_2(\Omega)}}:=\|\nabla u\|^2_{_{2,\Omega}}+\left\|u|_{_{\Gamma^{\infty}}}\right\|^2_{_{2,\Gamma^{\infty}}}$$
is an equivalent norm for $\mathcal{V}_2(\Omega)$.
\end{remark}

\section{The elliptic problem}\label{sec5}

\indent In this part we will summarize some auxiliary results regarding an elliptic boundary value problem, from where we will make the transition into the desired diffusion equation. Recall that $\Omega\subseteq\mathbb{R\!}^{\,2}$ will always denote the ramified domain with fractal boundary $\Gamma^{\infty}$ described in the previous section. Before posing the elliptic problem, we need to understand the interpretation of a notion for a co-normal derivative, which is not expected to exist in its classical definition over a fractal-type ramified boundary $\Gamma^{\infty}$.

\begin{definition}\label{Def-gen-normal}
Given $\alpha_{ij}\in L^{\infty}(\Omega)$ (for $i,\,j\in
\{1,2\}$), let $u\in W^{1,1}_{loc}(\Omega)$
be such that $$\displaystyle\sum^2_{i,j=1}\alpha_{ij}\partial_{x_j}u\partial_{x_i}\varphi\in L^1(\Omega),
\indent\textrm{for every}\,\,\,\varphi\in C^1(\overline{\Omega}).$$ Let $\mu$ be a Borel regular measure on the boundary
$\Gamma:=\partial\Omega$. If there exists a function $f\in L^1_{loc}(\Omega)$ such that $$\displaystyle\int_{\Omega}
\displaystyle\sum^2_{i,j=1}\alpha_{ij}\partial_{x_j}u\partial_{x_i}\varphi\,dx=\displaystyle\int_{\Omega}f\varphi\,dx+
\displaystyle\int_{\Gamma}\varphi\,d\mu,\indent\textrm{for all}\,\,\,\varphi\in C^1(\overline{\Omega}),$$
then we say that $\mu$ is the {\bf generalized co-normal derivative} of $u$, and we denote
$$\displaystyle\frac{\partial u}{\partial\nu_{_{\mathcal{A}}}}:=\displaystyle\sum^2_{i,j=1}
\alpha_{ij}\partial_{x_j}u\,\nu_{\mu_i}:=\mu,$$
and in this case, we say that $u$ is the {\bf weak solution} of the Eq.
$$\left\{
\begin{array}{lcl}
-\displaystyle\sum^2_{i,j=1}\partial_{x_i}(\alpha_{ij}\partial_{x_j}u)\,=\,f\,\,\,\textrm{in}\,\,\Omega\\[2ex]
\,\,\frac{\partial u}{\partial\nu_{_{\mathcal{A}}}}\,=\,1\,\,\,\textrm{on}\,\,\,\Gamma\\
\end{array}
\right.$$
\end{definition}

\indent When $\tau\in[1/2,\tau^{\ast})$, since $\Omega$ is an extension domain, one can follow the approach given by
Hinz, Rozanova-Pierrat and Teplyaev \cite{HINZ-ROZ-TEP22,HINZ-ROZ-TEP21} to give sense to a more structured notion of a generalized normal-type derivative. However, in the critical case $\tau=\tau^{\ast}$, such approach cannot be used since $\Omega$ fails to have the extension property, as shown in \cite{ACH08}.

\subsection{Existence of weak solutions}\label{subsec5.1}

We are now concerned with the elliptic mixed boundary value problem (\ref{E01}),
where $(f_0,g)\in L^p(\Omega)\times L^q_{\mu}(\Gamma)$ for given $p,\,q\in[1,\infty]$,\, $f_1,\,f_2\in L^2(\Omega)$, and $\nu_{\mu_j}$ stands as an interpretation of the $j$-component of the unit outer normal vector over rough boundaries (see for instance \cite[Definition 1]{LAN-AVS-23-1}, or \cite[Definition 1]{HEN-AVS19-1}). Here
\begin{equation}
\label{E02}\mathcal{A}u:=-\displaystyle\sum^2_{^{i,j=1}}\partial_{x_j}(\alpha_{ij}(x)\partial_{x_i}u)\,\,\,\,\,\,\,\,\textrm{and}\,\,\,\,\,\,\,\,\mathcal{B}u:=\displaystyle\sum^2_{i=1}\eta_i\partial_{x_i}u+\lambda u,
\end{equation}
where $\alpha_{ij}\in L^{\infty}(\Omega)$ (for each $i,\,j\in\{1,2\}$), \,$\eta_i\in L^{\zeta_i}(\Omega)$ with $\hat{\zeta}:=\zeta_1\wedge\zeta_2:=\min\{\zeta_1,\zeta_2\}>2$, \,$(\lambda,\beta)\in\mathbb{X\!}^{\,r,s}(\Omega;\Gamma^{\infty})$ for\, $r\wedge s>1$, and where $\frac{\partial u}{\partial\nu_{_{\mathcal{A}}}}$ is defined as in Definition \ref{Def-gen-normal}. We will denote $\bar{\alpha}:=(\alpha_{11},\alpha_{12},\alpha_{21},\alpha_{22})$ and $\bar{\eta}:=(\eta_1,\eta_2)$.\\
\indent We also assume that $\mathcal{A}$ is {\it uniformly elliptic} in the sense that there exists a positive constant $\alpha_0>0$ such that
\begin{equation}
\label{E03}\displaystyle\sum^2_{^{i,j=1}}\alpha_{ij}(x)\xi_i\xi_j\,\geq\,\alpha_0\,|\xi|^2,
\indent\,\forall\,x\in\Omega,\,\,\,\forall\,\xi=(\xi_1,\xi_2)\in\mathbb{R\!}^{\,2},
\end{equation}
\indent Consider now the bilinear closed (non-symmetric) energy form $(\mathcal{E},D(\mathcal{E})$ associated to problem (\ref{E01}), where $D(\mathcal{E}):=\mathcal{V}_2(\Omega)\times\mathcal{V}_2(\Omega)$ (for $\mathcal{V}_2(\Omega)$ given by (\ref{V-space})), and
\begin{equation}\label{E-form}
\mathcal{E}(u,v):=\displaystyle\int_{\Omega}\displaystyle\sum^2_{i,j=1}\alpha_{ij}\partial_{x_i}u\partial_{x_j}v\,dx+\displaystyle\int_{\Omega}\displaystyle\sum^2_{i=1}\eta_{i}(\partial_{x_i}u)v\,dx+\displaystyle\int_{\Omega}\lambda uv\,dx+\displaystyle\int_{\Gamma^{\infty}}\beta uv\,d\mu,\,\,\,\,\,\,\,\,\,\,\,\,\forall\,u,\,v\in\mathcal{V}_2(\Omega).
\end{equation}
Using H\"older's inequality together with (\ref{2.03}) and (\ref{2.04}), we clearly see that $\mathcal{E}(\cdot,\cdot)$ is continuous with $$|\mathcal{E}(u,v)|\leq\max\left\{\|\bar{\alpha}\|_{_{\infty,\Omega}},\,c_1\left(\|\bar{\eta}\|^2_{_{\hat{\zeta},\Omega}}+\|\lambda\|_{_{r,\Omega}}\right),\,c_2\|\beta\|_{_{s,\Gamma^{\infty}}}\right\}\|u\|_{_{\mathcal{V}_2(\Omega)}}\|v\|_{_{\mathcal{V}_2(\Omega)}}\,\,\,\,\,\,\,\,\,\,\,\,\forall\,u,\,v\in\mathcal{V}_2(\Omega).$$ Furthermore, recalling (\ref{E03}), using Remark \ref{Sharp-epsilon-ineq} for $\epsilon:=\left[4\left(\|\bar{\eta}\|^2_{_{\hat{\zeta},\Omega}}+\|\lambda\|_{_{r,\Omega}}\right)\right]^{-1}\alpha_0>0$, and recalling Remark \ref{Equivalent-norm}, we obtain that
$$\frac{\alpha_0}{2}\|\nabla u\|^2_{_{2,\Omega}}\,\leq\,\mathcal{E}(u,u)+\int_{\Gamma^{\infty}}\left(C_{\epsilon}\left(\|\bar{\eta}\|^2_{_{\hat{\zeta},\Omega}}+\|\lambda\|_{_{r,\Omega}}\right)+\beta\right)|u|^2\,d\mu,\,\,\,\,\,\,\,\,\,\,\,\,\forall\,u\in\mathcal{V}_2(\Omega).$$
Therefore, if one assumes that $\displaystyle\textrm{ess}\inf_{x\in\Gamma^{\infty}}\beta(x)\geq\beta_0$ for some constant $\beta_0\geq C_{\epsilon}\left(\|\bar{\eta}\|^2_{_{\hat{\zeta},\Omega}}+\|\lambda\|_{_{r,\Omega}}\right)$, recalling \cite[Theorem 1]{ACH08}, one deduces that the form $\mathcal{E}(\cdot,\cdot)$ is coercive. If $(f_0,g)\in\mathbb{X\!}^{\,p_0,q}(\Omega;\Gamma^{\infty})$ with $\min\{p_0,q\}>1$ and $f_1,\,f_2\in L^2(\Omega)$, then from (\ref{2.03}) and (\ref{2.04}) we immediately get that the functional $\Phi:\mathcal{V}_2(\Omega)\rightarrow\mathbb{R\!}\,$ defined by: 
\begin{equation}\label{Functional}
\Phi(v):=\displaystyle\int_{\Omega}\left(f_0v+\displaystyle\sum^2_{i=1}f_i\partial_{x_i}v\right)\,dx+\displaystyle\int_{\Gamma^{\infty}}gv\,d\mu
\end{equation}
lies in $\mathcal{V}_2(\Omega)^{\ast}$. Therefore, combining all the above together with Lax Milgram's lemma, we have shown the following result.

\begin{theorem}\label{solvability}
Assume that $\displaystyle\textrm{ess}\inf_{x\in\Gamma^{\infty}}\beta(x)\geq\beta_0$ for some constant $\beta_0\geq C_{\epsilon}\left(\|\bar{\eta}\|^2_{_{\hat{\zeta},\Omega}}+\|\lambda\|_{_{r,\Omega}}\right)$, for $\epsilon:=\left[4\left(\|\bar{\eta}\|^2_{_{\hat{\zeta},\Omega}}+\|\lambda\|_{_{r,\Omega}}\right)\right]^{-1}\alpha_0$, let $(f_0,g)\in\mathbb{X\!}^{\,p_0,q}(\Omega;\Gamma^{\infty})$ with $\min\{p_0,q\}>1$, and let $f_1,\,f_2\in L^2(\Omega)$. Then problem (\ref{E01}) is uniquely solvable over $\mathcal{V}_2(\Omega)$ in the sense that there exists a unique function $u\in\mathcal{V}_2(\Omega)$ such that
\begin{equation}\label{Existence-Elliptic}
\mathcal{E}(u,\varphi)=\displaystyle\int_{\Omega}\left(f_0\varphi+\displaystyle\sum^2_{i=1}f_i\partial_{x_i}\varphi\right)\,dx+\displaystyle\int_{\Gamma^{\infty}}g\varphi\,d\mu,\,\,\,\,\,\,\,\,\,\,\,\,\forall\,\varphi\in\mathcal{V}_2(\Omega).
\end{equation}
\end{theorem}

\indent To conclude this subsection, for an open set $D$ (not necessary contained in $\Omega$, but such that $\Omega\cap D\neq\emptyset$),
consider the space $\mathcal{V}_c(\Omega;D)$ as the closure in $\mathcal{V}_2(\Omega)$ of the space
$$\mathcal{W}(\Omega;D):=\left\{\varphi u\mid u\in\mathcal{V}_2(\Omega),\,
\,\,\varphi\in C^{\infty}_c(\mathbb{R\!}^N),\,\,\,\textrm{supp}[\varphi]\subseteq D\right\}.$$
Define the bilinear form $\mathcal{E}_{c}$ by
$$\mathcal{E}_{c}(\cdot,\cdot)=\mathcal{E}(\cdot,\cdot)\indent\,\textrm{with}\,\,\,\,D(\mathcal{E}_c):=\mathcal{V}_c(\Omega;D)\times \mathcal{V}_c(\Omega;D).$$
Clearly $\mathcal{V}_c(\Omega;D)$ is a closed subspace of $\mathcal{V}_2(\Omega)$, and thus under the same conditions of Theorem \ref{solvability}, we see that
$\mathcal{E}_{c}$ is continuous and coercive. Furthermore, the functional
$\Phi_c:=\Phi|_{_{\mathcal{V}_c(\Omega;D)}}\in\mathcal{V}_c(\Omega;D)^{\ast}$, and thus it also makes sense the existence of a function $\hat{u}\in\mathcal{V}_c(\Omega;D)$ fulfilling
\begin{equation}\label{c-Existence-Elliptic}
\mathcal{E}_{c}(\hat{u},\varphi)=
\displaystyle\int_{\Omega\cap D}f_0\varphi dx+\displaystyle\int_{\Omega\cap D}\displaystyle\sum^2_{i=1}f_i\partial_{x_i}\varphi\,dx+
\displaystyle\int_{\Gamma\cap D}g\varphi\,d\mu,\,\,\,\,\,\,\,\,\,\forall\,\varphi\in\mathcal{V}_c(\Omega;D).
\end{equation}
Finally, we recall that if $D=B(x_0;\rho)$ for $x_0\in\overline{\Omega}$ and $\rho>0$, and if $E\subseteq\mathbb{R\!}^{\,2}$ is such that $E\cap B(x_0;\rho)\neq\emptyset$, we write
\begin{equation}\label{Omega-cap-Ball}
E(x_0,\rho):=E\cap B(x_0;\rho).
\end{equation}

\subsection{Global and local estimates}\label{subsec5.2}

\indent In this section we develop key a priori estimates together with other norm estimates, which will play a very important role into the derivation of fine regularity results.\\
\indent To begin, given $u\in H^1(\Omega)$ and $k\geq0$ a fixed real number, we put
\begin{equation}\label{k-sets}
E_k:=\{x\in E\mid u\geq k\}\,\,\,\,\,\,\,\,\,\,\,\,\textrm{for}\,\,\,E\,\,\,\textrm{either}\,\,\,\Omega,\,\,\,\textrm{or}\,\,\,\Gamma:=\partial\Omega,
\end{equation}
and we define
\begin{equation}\label{k-function}
u_k:=(u-k)^+\in H^1(\Omega).
\end{equation}
Clearly $\partial_{x_i}u_k=\partial_{x_i}u\chi_{_{\Omega_k}}$ for each $i\in\{1,2\}$.\\
\indent The following lemma will be needed in order to establish the main result of this subsection.

\begin{lemma}\label{Lem0}
Given $m_1,\,m_2\geq1$ integers and
$\xi\in\{1,\ldots,m_1\}$, and $\zeta\in\{1,\ldots,m_2\}$, let $s_{1,\xi}\in[2,\infty)$ and $s_{2,\zeta}\in[2,\infty)$.
\begin{enumerate}
\item[(a)]\,\, Let $u\in H^1(\Omega)$. If there exists constants $k_0,\,\gamma_0\geq0$, and positive constants $\theta_{1,\xi},\,\theta_{2,\zeta}$, such that
for all  $k\geq k_0$, we have
\begin{equation}\label{Lem0-02-0}
\|u_{k}\|^2_{_{H^1(\Omega)}}\,\leq\,\gamma_0\left(\|u_{k}\|^2_{_{2,\Omega_k}}
+k^2\displaystyle\sum^{m_1}_{\xi=1}|\Omega_k|
^{^{\frac{2(1+\theta_{1,\xi})}{s_{1,\xi}}}}+k^2\displaystyle\sum^{m_2}_{\zeta=1}\mu(\Gamma_k)^{^{\frac{2(1+\theta_{2,\zeta})}{s_{2,\zeta}}}}\right),
\end{equation}
then there exists a constant $C^{\ast}_0>0$ (independent of $u$ and $k_0$) such that
\begin{equation}\label{Lem0-03}
\textrm{ess}\displaystyle\sup_{\overline{\Omega}}\{u\}
\,\leq\,C^{\ast}_0\left(\|u_{k_0}\|^2_{_{2,\Omega_{k_0}}}+k^2_0\right)^{1/2}.
\end{equation}
\item[(b)]\,\,Let $u\in\mathcal{V}_c(\Omega;D)$ for $D:=B(x_0,\rho)$ (where $x_0\in\overline{\Omega}$,\, $\rho>0$, and $\Omega(x_0,\rho)\neq\emptyset$).
If there exists constants $k_0,\,\gamma_0\geq0$, and positive constants $\theta_{1,\xi},\,\theta_{2,\zeta}$, such that
\begin{equation}\label{Lem0-02}
\|u_{k}\|^2_{_{\mathcal{V}_c(\Omega;D)}}\,\leq\,\gamma_0\|u_{k}\|^2_{_{2,\Omega_{k}(\rho)}}+\gamma_0k^2\left(\displaystyle\sum^{m_1}_{\xi=1}|\Omega_{k}(\rho)|
^{^{\frac{2(1+\theta_{1,\xi})}{s_{1,\xi}}}}+\displaystyle\sum^{m_2}_{\zeta=1}\mu(\Gamma_{k}(\rho))^{^{\frac{2(1+\theta_{2,\zeta})}{s_{2,\zeta}}}}\right),
\end{equation}
then there exists constants $C^{\ast}_0>0$ and $\sigma>0$ (independent of $u$ and $k_0$) such that
\begin{equation}\label{Lem0-03-2}
\textrm{ess}\displaystyle\sup_{\overline{\Omega}(x_0,\rho)}\{u\}
\,\leq\,C^{\ast}_0\rho^{\sigma}\left(\|u_{k_0}\|^2_{_{r,\Omega_{k_0}(\rho)}}
+k^2_0\right)^{1/2},
\end{equation}
where $r\in(2,\infty)$.
\end{enumerate}
\end{lemma}

\begin{proof}
It suffices to establish part (b). Given $\hat{k}>\hat{k}_0\geq0$ fixed real numbers, without loss of generality, we prove the result for
$$\theta_{1,\xi}=\theta_{2,\zeta}=\theta:=\min\left\{\min\{\theta_{1,\xi}\mid1\leq\xi\leq m_1\}\cup\min\{\theta_{2,\zeta}\mid1\leq\zeta\leq m_2\}\right\}$$
for each $\xi\in\{1,\ldots,m_1\}$ and $\zeta\in\{1,\ldots,m_2\}$. Since many of the steps in the proof can be derived in the same way as in the proof of Lemma \ref{Lem1} (in section \ref{sec6}), we will only sketch the main estimates and inequalities. Having said this,
define the sequences $\{k_n\}_{n\geq0}$,\, $\{y_n\}_{n\geq0}$,\, $\{z_n\}_{n\geq0}$ of positive real numbers, by:
\begin{equation}\label{Lem0-06}
k_n:=\hat{k}_0\chi_{_{\{n=0\}}}+(2-2^{-(n-1)})\hat{k}\chi_{_{\{n\geq1\}}},\,\,\,\,\,\,\,\,\,\,\,\,\,y_n:=\hat{k}^{-2}\|u_{k_n}\|^2_{_{2,\Omega_{k}(\rho)}},\,\,\,\,\,\,\,\,\,\,\,\,\,
z_n:=\displaystyle\sum^{m_1}_{\xi=1}|\Omega_{k_n}(\rho)|^{^{\frac{2}{s_{1,\xi}}}}
+\displaystyle\sum^{m_2}_{\zeta=1}\mu(\Gamma_{k_n}(\rho))^{^{\frac{2}{s_{2,\zeta}}}}.
\end{equation}
Clearly $k_n\geq\hat{k}_0$ for each nonnegative integer $n$. Taking into account
the definitions of the sequences defined above together with H\"older's inequality and (\ref{2.03}), we calculate and get that
\begin{equation}\label{Lem0-08}
\hat{k}^2y_{n+1}\,\leq\,c\,4^{n+1}y_n^{^{\frac{r-2}{r}}}
\|u_{k_{n+1}}\|^2_{_{\mathcal{V}_c(\Omega;D)}}\,\,\,\,\,\,\,\,\textrm{and}\,\,\,\,\,\,\,\,4^{-(n+1)}\hat{k}^2z_{n+1}\,\leq\,c\|u_{k_{n}}\|^2_{_{\mathcal{V}_c(\Omega;D)}}.
\end{equation}
where $c>0$ is a constant that varies from line to line (this will be assumed for the remaining of the proof).
Applying assumption (\ref{Lem0-02}) and using (\ref{Lem0-08}) yield that
\begin{equation}\label{Lem0-11}
y_{n+1}\,\leq\,c'\,8^n(y^{1+\delta}_n+z^{1+\theta}_ny^{\delta}_n)\,\,\,\,\,\textrm{and}\,\,\,\,\,
z_{n+1}\,\leq\,c'\,8^n(y_n+z^{1+\theta}_n),\,\,\textrm{for each integer}\,\,n\geq0,
\end{equation}
for some constant $c'>0$, where $\delta:=(r-2)/r$. Now, it is clear that
$y_0\,\leq\,\hat{k}^{-2}\|u_{k}\|^2_{_{2,\Omega_{k}(\rho)}}.$
and by using (\ref{Lem0-02}) and recalling that $|\cdot|$ and $\mu(\cdot)$ are upper $2$-Ahlfors and upper $d$-Ahlfors measures,
respectively, we proceed as above to find that
$$(\hat{k}-k_0)^2z_0\,\leq\,c\gamma_0\left(\rho^2\left\{\|u_{k_0}\|^2_{_{r,\Omega_{k_0}(\rho)}}\right\}+2k^2_0
\displaystyle\sum^{m_1}_{\xi=1}\rho^{^{\frac{4(1+\theta)}{s_{1,\xi}}}}
+2k^2_0\displaystyle\sum^{m_2}_{\zeta=1}\rho^{^{\frac{2d(1+\theta)}{s_{2,\zeta}}}}\right).$$
Thus setting
$$c''_{\rho}:=c\gamma_0\max\left\{\rho^2,\,2\displaystyle\sum^{m_1}_{\xi=1}\rho^{^{\frac{4(1+\theta)}{s_{1,\xi}}}}
+2\displaystyle\sum^{m_2}_{\zeta=1}\rho^{^{\frac{2d(1+\theta)}{s_{2,\zeta}}}}\right\},$$
we obtain that
\begin{equation}\label{Lem0-13}
z_0\,\leq\,\displaystyle\frac{c''_{\rho}}{(\hat{k}-k_0)^2}\left(\|u_{k_0}\|^2_{_{r,\Omega_{k_0}(\rho)}}+k^2_0\right),
\,\,\,\,\,\,\,\,\,\,\,\,\textrm{for all}\,\,\,\hat{k}\geq k_0.
\end{equation}
Selecting $$l:=\min\left\{\delta,\,\frac{\theta}{1+\theta}\right\}\,\,\,\,\,\textrm{and}\,\,\,\,\,
\eta:=\min\left\{\frac{1}{(2c''_{\rho})^{^{\frac{1}{\delta}}}8^{^{\frac{1}{\delta l}}}},\,\frac{1}
{(2c''_{\rho})^{^{\frac{1+\theta}{\theta}}}\,8^{^{\frac{1}{\theta l}}}}\right\},$$
and choosing
\begin{equation}\label{Lem0-14}
\hat{k}:=\max\left\{\displaystyle\frac{\rho^2}{\sqrt{\eta}}\|u_{k_0}\|_{_{r,\Omega_{k_0}(\rho)}},\,
\displaystyle\frac{\sqrt{c''_{\rho}}}{\eta^{^{\frac{1}{2(1+\theta)}}}}
\left(\|u_{k_0}\|_{_{r,\Omega_{k_0}(\rho)}}+k^2_0\right)^{1/2}+k_0\right\},
\end{equation}
we see that
\begin{equation}\label{Lem0-15}
\hat{k}\,\leq\,c^{\star}\rho^{\sigma}\left(\|u_{k_0}\|_{_{r,\Omega_{k_0}(\rho)}}+k^2_0\right)^{1/2}
\end{equation}
for some constants $c^{\star},\,\sigma>0$, and moreover the selection implies that
\begin{equation}\label{Lem0-16}
y_0\,\leq\,\eta\,\,\,\,\,\,\,\,\textrm{and}\,\,\,\,\,\,\,\,z_0\,\leq\,\eta^{^{\frac{1}{1+\theta}}}.
\end{equation}
Therefore, the sequences $\{y_n\}$,\, $\{z_n\}$ satisfy the conditions of Lemma \ref{lemma1}.
Thus, applying Lemma \ref{lemma1} for $\hat{k}$ given by (\ref{Lem0-14})
shows that $\displaystyle\lim_{n\rightarrow\infty}z_n=0$, and consequently
\begin{equation}\label{Lem0-17}
u(t)\,\leq\,\displaystyle\lim_{n\rightarrow\infty}k_n=2\hat{k}\,\,\,\,\,\,\,\,\textrm{a.e. in}\,\,\overline{\Omega}.
\end{equation}
Combining (\ref{Lem0-17}) with (\ref{Lem0-15}) we obtain (\ref{Lem0-03-2}), as desired.
The derivation of (\ref{Lem0-03}) follows in the exact way (with simpler arguments).
\end{proof}

\indent Now we establish global boundedness and a priori estimates for weak solutions to problem (\ref{E01}).

\begin{theorem}\label{Elliptic-Infinity-Case}
Let $(f,g)\in\mathbb{X\!}^{\,p_0,q}(\Omega;\Gamma^{\infty})$ for $\min\{p_0,q\}>1$, assume that $\beta\in L^s_{\mu}(\Gamma^{\infty})$ is as in Theorem \ref{solvability}, and let $f_i\in L^{p_i}(\Omega)$ for $i\in\{1,2\}$ with $\min\{p_1,p_2\}>2$.
\begin{enumerate}
\item[(a)]\,\, If $u\in\mathcal{V}_2(\Omega)$ is a weak solution of problem (\ref{E01}), then there exists a constant $C^{\ast}_1=C^{\ast}_1(p,q,|\Omega,\mu(\Gamma^{\infty}))>0$ (independent of $u$) such that
\begin{equation}\label{Thm0-01}
\max\left\{\|u\|_{_{\infty,\Omega}},\|u\|_{_{\infty,\Gamma}}\right\}\,\leq\,C^{\ast}_1\left(\displaystyle\sum^2_{i=0}\|f_i\|_{_{p_i,\Omega}}+\|g\|_{_{q,\Gamma^{\infty}}}\right).
\end{equation}
\item[(b)]\,\, Let $\rho\in(0,1)$ and $x_0\in\overline{\Omega}$ be such that $\Omega(x_0,\rho)\neq\emptyset$, and
let $D:=B(x_0,\rho)$. If $u\in\mathcal{V}_c(\Omega;D)$ solves (\ref{c-Existence-Elliptic}), then there exists a constant $C^{\star}_1=C^{\star}_1(p,q,|\Omega|,\mu(\Gamma^{\infty}))>0$ (independent of $u$) and a constant $\sigma^{\star}>0$, such that
\begin{equation}\label{Thm0-01b}
\max\left\{\|u\|_{_{\infty,\Omega(x_0,\rho)}},\|u\|_{_{\infty,\Gamma(x_0,\rho)}}\right\}
\,\leq\,C^{\star}_1\rho^{\sigma^{\star}}\left(\displaystyle\sum^2_{i=0}\|f_i\|_{_{p_i,\Omega}}+\|g\|_{_{q,\Gamma^{\infty}}}\right).
\end{equation}
\end{enumerate}
\end{theorem}

\begin{proof}
Again  we will only sketch the main steps of part (b).
Let $u\in\mathcal{V}_c(\Omega;D)$ be a solution of problem (\ref{c-Existence-Elliptic}), and let $u_{k}$ be the function given by (\ref{k-function}). Then we have
\begin{equation}\label{thm0-03}
\mathcal{E}_{c}(u,u_{k})=\displaystyle\int_{\Omega_{k}(\rho)}f_0u_{k}\,dx+
\displaystyle\int_{\Omega_{k}(\rho)}\displaystyle\sum^2_{i=1}f_i\partial_{x_i}u_{k}\,dx+
\displaystyle\int_{\Gamma^{\infty}_{k}(\rho)}gu_{k}\,d\mu.
\end{equation}
After multiple calculation procedures (see the proof of Theorem \ref{Thm-Partial-Bounded}) we arrive at
\begin{equation}\label{thm0-12}
\min\left\{\beta_0,\displaystyle\frac{c_0}{2}\right\}\|u_{k}\|^2_{_{\mathcal{V}_c(\Omega;D)}}
\,\leq\,\displaystyle\int_{\Omega_{k}(\rho)}\left(\frac{1}{k}|f_0|+C_0(\bar{\eta})\right)(|u_{k}|^2+k^2)\,dx
+\frac{1}{c_0}\displaystyle\int_{\Omega_{k}(\rho)}\displaystyle\sum^N_{i=1}|f_i|^2\,dx
+\displaystyle\int_{\Gamma^{\infty}_{k}(\rho)}\frac{1}{k}|g|(|u_{k}|^2+k^2)\,d\mu,
\end{equation}
where $C_0(\bar{\eta}):=\left(8C_{\epsilon^2_0}\|\bar{\eta}\|^2_{_{\hat{\zeta},\Omega}}\right)/c_0+\lambda^-$ for $\epsilon_0:=c_0\left(8\|\bar{\eta}\|_{_{\hat{\zeta},\Omega}}\right)^{-1}$.
We now estimate the terms in (\ref{thm0-12}). Put
\begin{equation}\label{thm0-13}
k^2_0:=\|f_0\|^2_{_{p_0,\Omega}}+\displaystyle\sum^N_{i=1}\|f_i\|^2_{_{p_i,\Omega}}+\|g\|^2_{_{q,\Gamma^{\infty}}}\,\,\,\,\,\textrm{and}\,\,\,\,\,k\geq k_0,
\end{equation}
and apply H\"older's inequality together with Young's inequality and Remark \ref{Equivalent-norm} to deduce that\\[2ex]
\begin{equation}\label{thm0-14}
\displaystyle\int_{\Omega_{k}(\rho)}\frac{1}{k}|f_0||u_{k}|^2\,dx\,\leq\,\displaystyle\frac{k_0}{k}
\|u_{k}\|^2_{_{\frac{2p_0}{p_0-1},\Omega_{k}(\rho)}}
\,\leq\,\epsilon_1\|u_{k}\|^2_{_{\mathcal{V}_c(\Omega;D)}}+C_{\epsilon_1}\|u_{k}\|^2_{_{2,\Omega_{k}(\rho)}},
\end{equation}
for all $\epsilon_1>0$, and for some constant $C_{\epsilon_1}>0$. Proceeding in the same way as in the derivation of (\ref{thm0-14}), we get that
\begin{equation}\label{thm0-16}
\displaystyle\int_{\Omega_{k}(\rho)}\lambda|u_{k}|^2\,dx
\,\leq\,\|\lambda\|_{_{r,\Omega}}\left(\epsilon_2\|u_{k}\|^2_{_{\mathcal{V}_c(\Omega;D)}}+
C_{\epsilon_2}\|u_{k}\|^2_{_{2,\Omega_{k}(\rho)}}\right)
\,\,\,\,\textrm{and}\,\,\,\,
\displaystyle\int_{\Gamma^{\infty}_{k}(\rho)}\frac{1}{k}|g||u_{k}|^2\,d\mu
\,\leq\,\epsilon_3\|u_{k}\|^2_{_{\mathcal{V}_c(\Omega;D)}}+C_{\epsilon_3}\|u_{k}\|^2_{_{2,\Omega_{k}(\rho)}},
\end{equation}
for all $\epsilon_2,\,\epsilon_3>0$, and for some constants $C_{\epsilon_2},\,C_{\epsilon_3}>0$. In a similar way, for $k\geq k_0$ we see that
\begin{equation}\label{thm0-17}
\displaystyle\int_{\Omega_{k}(\rho)}\frac{1}{k}|f_0|\,dx\,\leq\,|\Omega_{k}(\rho)|^{^{1-\frac{1}{p_0}}},\,\,\,\,\,\,\,\,\,\,\,\,\,\,\,\,
\displaystyle\int_{\Gamma^{\infty}_{k}(\rho)}\frac{1}{k}|g|\,d\mu\,\leq\,\mu(\Gamma^{\infty}_{k}(\rho))^{^{1-\frac{1}{q}}},
\end{equation}
and
\begin{equation}\label{thm0-18}
\displaystyle\int_{\Omega_{k}(\rho)}\lambda\,dx
\,\leq\,\|\lambda\|_{_{r,\Omega}}|\Omega_{k}(\rho)|^{^{1-\frac{1}{r}}},
\,\,\,\,\,\,\,\,\,\,\,\,
\displaystyle\int_{\Omega_{k}(\rho)}
\displaystyle\sum^2_{i=1}|f_i|^2\,dx\,\leq\,k^2\displaystyle\sum^2_{i=1}|\Omega_{k}(\rho)|^{^{1-\frac{2}{p_i}}}.
\end{equation}
Since the maps $\mathcal{V}_c(\Omega;D)\hookrightarrow L^{r_0}(\Omega)$ and $\mathcal{V}_c(\Omega;D)\hookrightarrow L^{s_0}_{\mu}(\Gamma^{\infty})$ are compact for each $r_0,\,s_0\in[1,\infty)$, choosing
$$\epsilon_1=\epsilon_3=\frac{1}{6}\min\left\{1,\displaystyle\frac{c_0}{2}\right\},\,\,\,\,\,\,\,\,\,\,
\epsilon_2=\frac{1}{6\|\hat{\mathfrak{h}}\|_{_{r,\Omega}}}\min\left\{1,\displaystyle\frac{c_0}{2}\right\},$$
selecting $$\gamma_0:=\displaystyle\frac{2\max\left\{1,C_{\epsilon_1},C_{\epsilon_3},\|\lambda\|_{_{r,\Omega}}C_{\epsilon_2}\right\}}{\min\{1,c_0/2\}},$$
and putting
$$\theta_{1,1}:=\frac{(p_0-1)\tilde{p}_0-2p_0}{2p_0}-1,\,\,\,\,\,\,\,
\theta_{1,2}:=\frac{(r-1)\tilde{r}-2r}{2r}-1,\,\,\,\,\,\,\,\theta_{1,i+2}:=\frac{(p_i-2)\tilde{p}_i-2p_i}{2p_i}-1,\,\,\,\,\,\,\,\theta_{2,1}:=\frac{(q-1)\tilde{q}-2q}{2q}-1,$$
$$s_{1,1}:=\tilde{p}_0,\,\,\,\,\,\,\,\,\,\,\,\,\,\,\,s_{1,2}:=\tilde{r}\,\,\,\,\,\,\,\,\,\,\,\,\,\,\,s_{1,i+2}:=\tilde{p}_i,\,\,\,\,\,\,\,\,\,\,\,\,\,\,\,s_{2,1}:=\tilde{q}\,\,\,\,\,\,\,\,\,\,\,\,\,\,\,\,\textrm{(}\,i\in\{1,2\}\,\textrm{)},$$
for some constants $\tilde{p}_0,\,\tilde{r},\,\tilde{p}_1,\,\tilde{p}_2,\,\tilde{q}\in(2,\infty)$ ensuring that $\theta_{j,m}>0$ for each $(j,m)\in\{1,2\}\times\{1,\ldots,4\}$,
we insert (\ref{thm0-14}), (\ref{thm0-16}), (\ref{thm0-17}), and (\ref{thm0-18}) into (\ref{thm0-12}), to arrive at
\begin{equation}\label{thm0-25}
\|u_{k}\|^2_{_{\mathcal{V}_c(\Omega;D)}}\,\leq\,\gamma_0\|u_{k}\|^2_{_{2,\Omega_{k}(\rho)}}+\gamma_0k^2\left(\displaystyle\sum^{4}_{\xi=1}|\Omega_{k}(\rho)|
^{^{\frac{2(1+\theta_{1,\xi})}{s_{1,\xi}}}}+\displaystyle\sum^{2}_{\zeta=1}\mu(\Gamma^{\infty}_{k}(\rho))^{^{\frac{2(1+\theta_{2,\zeta})}{s_{2,\zeta}}}}\right).
\end{equation}
Therefore, all the conditions of Lemma \ref{Lem1}(b) are fulfilled, and thus
applying Lemma \ref{Lem1} to $u$ and to $-u$ (where $-u$ can be seen as the unique weak solution to problem (\ref{c-Existence-Elliptic}) with respect to the data $-f_i$ and $-g$ for $i\in\{0,1,2\}$), and using the fact that $u=0$ over $\Gamma\setminus\Gamma^{\infty}$, one has
\begin{equation}\label{thm0-26}
\max\left\{\|u\|_{_{\infty,\Omega(x_0,\rho)}},\|u\|_{_{\infty,\Gamma(x_0,\rho)}}\right\}
\,\leq\,C^{\ast}_0\rho^{\sigma}\left(\|u_{k_0}\|^2_{_{\hat{r},\Omega_{k_0}(\rho)}}+\|f_0\|^2_{_{p_0,\Omega}}+\displaystyle\sum^N_{i=1}\|f_i\|^2_{_{p_i,\Omega}}+\|g\|^2_{_{q,\Gamma^{\infty}}}\right)^{1/2},
\end{equation}
for some constant $C^{\ast}_0>0$ and for some $\hat{r}\in(2,\infty)$. Using the coercivity of the form $\mathcal{E}_c(\cdot,\cdot)$ together with Theorem \ref{embeddings-trace}(a), we clearly deduce that
\begin{equation}\label{thm0-27}
\|u_{k_0}\|^2_{_{\hat{r},\Omega_{k_0}(\rho)}}\,\leq\,c\left(\|f_0\|^2_{_{p_0,\Omega}}+\displaystyle\sum^N_{i=1}\|f_i\|^2_{_{p_i,\Omega}}+\|g\|^2_{_{q,\Gamma^{\infty}}}\right),
\end{equation}
and from here a combination of (\ref{thm0-26}) and (\ref{thm0-27}) yield (\ref{Thm0-01b}), as desired.
\end{proof}

\indent Next we give an useful local gradient estimate for weak solutions, which will be established under a small additional assumption. In fact, for the rest of the paper, we will assume that
\begin{equation}\label{sharp-coercivity}
\lambda\geq0\,\,\,\textrm{a.e. and}\,\,\,\beta-C_{\epsilon}\|\bar{\eta}\|^2_{_{\hat{\zeta},\Omega}}\geq\beta_0\,\,\,\mu\textrm{-a.e},\,\,\,\,\,\,\,\,\textrm{or}\,\,\,\,\,\,\,\,\lambda-C_{\epsilon}\|\bar{\eta}\|^2_{_{\hat{\zeta},\Omega}}\geq\lambda_0\,\,\,\textrm{a.e. and}\,\,\,\beta\geq0\,\,\,\mu\textrm{-a.e.},
\end{equation}
for some positive constants $C_{\epsilon}$,\, $\lambda_0$ and $\beta_0$ (as in subsection \ref{subsec5.1}).
It is straightforward to see that when (\ref{sharp-coercivity}) holds, the bilinear form $\mathcal{E}(\cdot,\cdot)$ is coercive.

\begin{lemma}\label{gradient-local}
Under all the conditions of Theorem \ref{Elliptic-Infinity-Case} together with (\ref{sharp-coercivity}), if $u\in\mathcal{V}_2(\Omega)$ is a weak solution of (\ref{E01}), then for each $x_0\in\overline{\Omega}$,\, $k\geq k_0\geq0$, and $0<\varrho_0<\rho\leq\rho^{\ast}<1$, there exists a constant $c>0$ such that
\begin{equation}\label{G01}
\|\nabla u\|^2_{_{2,\Omega_k(\varrho_0)}}\,\leq\,c\left(\frac{\|u_k\|^2_{_{\Omega_k(\rho)}}}{(\rho-\varrho_0)^2}+\Upsilon(f_0,f_1,f_2,g)\right),
\end{equation}
where $u_k\in\mathcal{V}_2(\Omega)$ is defined as in (\ref{k-function}), and
\begin{equation}\label{Upsilon}
\Upsilon(f_0,f_1,f_2,g):=\|f_0\|^2_{_{p_0,\Omega}}+\displaystyle\sum^N_{i=1}\|f_i\|^2_{_{p_i,\Omega}}+\|g\|^2_{_{q,\Gamma^{\infty}}}.
\end{equation}
\end{lemma}

\begin{proof}
Under all the conditions of the lemma, take $\phi\in C^{\infty}_c(B(x_0;\rho_0)$ with $0\leq\phi\leq1$,\,$\phi|_{_{B(x_0;\varrho_0)}}=1$,\,$\phi|_{_{B(x_0;\rho^{\ast})\setminus B(x_0;\rho)}}=0$, and $\|\nabla\phi\|_{_{\infty,B(x_0;\rho)}}\leq(\rho-\varrho_0)^{-1}$. Then, given $u\in\mathcal{V}_2(\Omega)$ a weak solution of problem (\ref{E01}), we test (\ref{Existence-Elliptic}) with the function $u_{_{k,\phi}}:=\phi^2u_k\in\mathcal{V}_c(\Omega;D)$ (for $D=B(x_0;\rho)$), recalling (\ref{E03}) and (\ref{sharp-coercivity}), we use the product rule to get that\\[2ex]
$c_0\|\phi\nabla u_k\|^2_{_{2,\Omega_k(\varrho_0)}}\,\leq\,2\left|\displaystyle\int_{\Omega_k(\rho)}\displaystyle\sum^2_{i,j=1}\alpha_{ij}\phi(\partial_{x_i}u_k)u_k\partial_{x_j}\phi\,dx\right|+\left|\displaystyle\int_{\Omega_k(\rho)}\displaystyle\sum^2_{i=1}\eta_{i}\phi(\partial_{x_i}u_k)\phi u_k\,dx\right|+$\\
\begin{equation}\label{G02}
+\left|\displaystyle\int_{\Omega_k(\rho)}\displaystyle f_0\phi^2u_k\,dx\right|+
\left|\displaystyle\int_{\Omega_k(\rho)}\displaystyle\sum^2_{i=1}f_{i}\partial_{x_i}(\phi^2u_k)\,dx\right|+
\left|\displaystyle\int_{\Gamma_k(\rho)}\displaystyle g\phi^2u_k\,d\mu\right|.
\end{equation}
Using Young's inequality, it is clear that
\begin{equation}\label{G03}
2\left|\displaystyle\int_{\Omega_k(\rho)}\displaystyle\sum^2_{i,j=1}\alpha_{ij}\phi(\partial_{x_i}u_k)u_k\partial_{x_j}\phi\,dx\right|\,\leq\,\epsilon\|\phi\nabla u_k\|^2_{_{2,\Omega_k(\rho)}}+2\epsilon^{-1}\|\bar{\alpha}\|^2_{_{\infty,\Omega}}\|u_k\nabla\phi\|^2_{_{2,\Omega_k(\rho)}},
\end{equation}
for all $\epsilon>0$. Also, in views of Remark \ref{Sharp-epsilon-ineq}, we have
\begin{equation}\label{G04}
\left|\displaystyle\int_{\Omega_k(\rho)}\displaystyle\sum^2_{i=1}\eta_{i}\phi(\partial_{x_i}u_k)\phi u_k\,dx\right|\,\leq\,
2\epsilon\|\nabla(\phi u_k)\|^2_{_{2,\Omega_k(\rho)}}+
C_{\epsilon^2}\|\bar{\eta}\|^2_{_{\hat{\zeta},\Omega}}\epsilon^{-1}\|\phi u_k\|^2_{_{2,\Omega_k(\rho)}},
\end{equation}
for all $\epsilon>0$ and for some constant $C_{\epsilon^2}>0$, where $\hat{\zeta}:=\min\{\zeta_1,\zeta_2\}$. Using (\ref{G03}) in (\ref{G04}), and proceeding in the same way for the rest of the terms in the right hand side of the inequality (\ref{G02}), we arrive at the following inequality:\\[2ex]
$\|\phi\nabla u\|^2_{_{2,\Omega_k(\rho)}}$
\begin{equation}\label{G05}
\leq\,8\epsilon\|\phi\nabla u\|^2_{_{2,\Omega_k(\rho)}}+2\left(\epsilon^{-1}\|\bar{\alpha}\|^2_{_{\infty,\Omega}}+3\epsilon\right)\|u_k\nabla\phi\|^2_{_{2,\Omega_k(\rho)}}+2\left(\epsilon^{-1}\|\bar{\eta}\|^2_{_{\hat{\zeta},\Omega}}C_{\epsilon^2}+3\epsilon\right)\|\phi u_k\|^2_{_{2,\Omega_k(\rho)}}+C'_{\epsilon}\Upsilon(f_0,f_1,f_2,g),
\end{equation}
for every $\epsilon>0$ and for some constants $C_{\epsilon^2},\,C'_{\epsilon}>0$. Selecting $\epsilon>0$ suitably and using the conditions over the cut-off function $\phi$, we are lead into inequality (\ref{G01}), as desired.
\end{proof}

\indent Taking into account Theorem \ref{embeddings-trace} and Lemma \ref{gradient-local}, we can proceed in a exact way as in \cite[proof of Corollary 5.6]{VELEZ2013-1} to deduce the following consequences.

\begin{lemma}\label{consequences}
Under all the conditions of Lemma \ref{gradient-local},
there exists constants $c'_{i}>0$ ($1\leq i\leq 4$) such that for each $x_0\in\overline{\Omega}$ and for every $0<\varrho_0<\rho<\rho^{\ast}<1$,
a weak solution $u\in\mathcal{V}_2(\Omega)$ of problem (\ref{E01}) fulfills the following conditions.
\begin{enumerate}
\item[(a)]\,\, $\|u_k\|^2_{_{2,\Omega_k(\varrho_0)}}\,\leq\,c'_{1}\,|\Omega_k(\varrho_0)|^{\delta_1}\left[(\rho-\varrho_0)^{-2}\,
\|u_k\|^2_{_{2,\Omega_k(\rho)}}+\Upsilon(f_0,f_1,f_2,g)\right]$;
\item[(b)]\,\, $|\Omega_h(\varrho_0)|^{\delta_2}\,\leq\,c'_{2}\,(h-k)^{-2}
\left[(\rho-\varrho_0)^{-2}\,
\|u_k\|^2_{_{2,\Omega_k(\rho)}}+\Upsilon(f_0,f_1,f_2,g)\right]$;
\item[(c)]\,\, $\mu(\Gamma^{\infty}_h(\varrho_0))^{\delta_3}\,\leq\,c'_{3}\,(h-k)^{-2}
\left[(\rho-\varrho_0)^{-2}\,
\|u_k\|^2_{_{2,\Omega_k(\rho)}}+\Upsilon(f_0,f_1,f_2,g)\right]$;
\item[(d)]\,\, $|\Omega_h(\varrho_0)|^{\delta_4}\,\leq\,c'_{4}\,(h-k)^{-2}(|\Omega_k(\varrho_0)|-
|\Omega_h(\varrho_0)|)\left[(\rho-\varrho_0)^{-2}\,
\|u_k\|^2_{_{2,\Omega_k(\rho)}}+\Upsilon(f_0,f_1,f_2,g)\right]$,
\end{enumerate}
for all $k\geq0$ such that $|\Omega_k(\rho)|\leq\frac{1}{2}|\Omega(x_0;\rho)|$, for every $h>k$, where $\delta_i\in(0,1)$ ($1\leq i\leq 4$) are arbitrarily fixed, and where $\Upsilon(f_0,f_1,f_2,g)$ is given by (\ref{Upsilon}).
\end{lemma}

\subsection{An oscillation result}\label{subsec5.3}

\indent This section is devoted to establish a very important oscillation lemma, which, when combined with the previous results,
will lead to the fulfilment of the first main result of the paper. We will assume all the conditions of Theorem \ref{Elliptic-Infinity-Case} together with condition (\ref{sharp-coercivity}).\\
\indent To begin, for a bounded function $v$ on $\overline{\Omega}$, and for $x_0\in\overline{\Omega}$ and $\rho>0$, we will denote\\
\begin{equation}\label{6.01}
\mathfrak{m}^{\ast}(v;\rho):=\displaystyle\textrm{ess}\,{\sup_{^{\overline{\Omega}(x_0;\rho)}}}v(x)\indent\textrm{and}\indent
\mathfrak{m}_{\ast}(v;\rho):=\displaystyle\textrm{ess}{\inf_{^{\overline{\Omega}(x_0;\rho)}}}v(x),
\end{equation}
and
\begin{equation}\label{6.02}
\textrm{Osc}(v;\rho):=\mathfrak{m}^{\ast}(v;\rho)-\mathfrak{m}_{\ast}(v;\rho).
\end{equation}
Then, one has the following key oscillation estimate.

\begin{lemma}\label{oscilation}
For each $x_0\in\overline{\Omega}$ and $\rho\in(0,\rho^{\ast})$, if $v\in\mathcal{V}_2(\Omega)$ is a bounded function such that 
\begin{equation}\label{6.03}
\mathcal{E}(v,w)=0\,\,\,\,\,\,\,\,\,\,\,\textrm{for all}\,\,w\in\mathcal{V}_c(\Omega;D),
\end{equation}
where $D=B(x_0;\rho)$, then there exists a constant $\varsigma_0\in(0,1)$ such that
\begin{equation}
\label{6.04}\textrm{Osc}(v;\rho/4)\,\leq\,\varsigma_0\,\textrm{Osc}(v;\rho).
\end{equation}
\end{lemma}

\begin{proof}
Given $k_0\geq0$ and $\rho^{\ast}\in(0,1)$, fix $\rho\in(0,\rho^{\ast}]$,
and select a real number $k\geq k_0$ such that
$|\Omega_k(\rho)|\leq\frac{1}{2}|\Omega(x_0;\rho/4)|$. Then
$$|\Omega_k(\varrho_0)|\leq\frac{1}{2}|\Omega(x_0;\varrho_0)|\indent\indent\,\forall\,\varrho_0\in[\rho/4,\rho].$$
Let $h>k\geq k_0$, let $\varrho_0,\,\sigma_0\in[\rho/4,\rho]$ be such that $\varrho_0<\sigma_0$, and put
\begin{equation}
\label{6.05}\mathcal{T}(h,\varrho_0):=\|v_h\|^{r_0}_{_{2,\Omega_h(\varrho_0)}}
\,\,\,\,\,\,\,\,\textrm{and}\,\,\,\,\,\,\,\,
\Phi(h,\varrho_0):=\left|\left\|(\chi_{_{\Omega_h(\varrho_0)}},\chi_{_{\Gamma^{\infty}_h(\varrho_0)}})\right\|\right|^{^{\frac{r^2_0}{2-r_0}}}_{_{\mathbb{X\!}^{^{\,\frac{2r_0}{2-r_0},\frac{2s_0}{2-s_0}}}(\Omega;\Gamma)}}.
\end{equation}
for $r_0,\,s_0\in(1,2)$ arbitrarily fixed.
As $v\in\mathcal{V}_2(\Omega)$ solves (\ref{6.03}), we see that the conclusion of
Proposition \ref{gradient-local} is valid for $\Upsilon_p(f_0,f_1,f_2,g)=0$.
By Lemma \ref{consequences}(a) one obtains
\begin{equation}\label{6.07a}
\mathcal{T}(h,\varrho_0)\,\leq\,
c(\sigma_0-\varrho_0)^{-r_0}\,\mathcal{T}(k,\sigma_0)\,\Phi(k,\sigma_0)^{^{\frac{2}{r_0}}},
\end{equation}
where $c>0$ will denote a generic constant which may vary from line to line. Furthermore, an application of Lemma \ref{consequences}(b,c) entails that
\begin{equation}\label{6.09}
\Phi(h,\varrho_0)^{^{\frac{2-r_0}{2}}}\,\leq\,c(h-k)^{-r_0}(\sigma_0-\varrho_0)^{-r_0}\,\mathcal{T}(k,\sigma_0).
\end{equation}
Letting $\gamma$ be a positive solution of the equation
\begin{equation}
\label{6.07}r_0\gamma^2=(\gamma+1)(2-r_0).
\end{equation}
Replacing $\gamma$ in (\ref{6.07a}) and combining this with (\ref{6.09}), we have that
\begin{equation}\label{6.10}
\mathbf{\varpi}(h,\varrho_0)\,\leq\,c(h-k)^{-r_0}(\sigma_0-\varrho_0)^{-r_0-r_0\gamma}\,\mathbf{\varpi}(k,\sigma_0)^{\delta},
\end{equation}\indent\\
where $\mathbf{\varpi}(\cdot,\cdot):=\Phi(\cdot,\cdot)^{^{\frac{2-r_0}{2}}}\mathcal{T}(\cdot,\cdot)^{\gamma}$ and
$\delta:=1+\gamma^{-1}>1$. Letting $\sigma_0=\rho$ and $\varrho_0=\rho/2$, we apply Lemma \ref{lemma2} over
$\mathbf{\varpi}(\cdot,\cdot)$ to deduce that $\mathbf{\varpi}(k_0+\xi,\rho/2)=0$, for
$$\xi^2=c\rho^{-r_0-r_0\gamma}\mathbf{\varpi}_2(k_0,\rho)^{\delta-1}\,\leq\, c\rho^{-r_0-r_0\gamma}\rho^{^{\frac{2-r_0}{\gamma}}}
\mathcal{T}(k_0,\rho)=c\rho^{-2}\mathcal{T}(k_0,\rho).$$
Thus, as $v$ vanishes over $\Gamma\setminus\Gamma^{\infty}$, one deduces that $v(x)\leq k_0+\xi$ for almost all $x\in\overline{\Omega}(a,\rho/2)$, with
\begin{equation}
\label{6.11}\mathfrak{m}^{\ast}(v;\rho/2)\,\leq\,k_0+C^{\star}\,\rho^{-1}\left\|v-k
\right\|_{_{2,\Omega_{k_0}(\rho)}},
\end{equation}
for some constant $C^{\star}>0$. Next, for an arbitrary $\varrho_0\in[\rho/4,\rho]$ and for an integer $n\geq0$, put
\begin{equation}
\label{6.12}\omega_n(v;\varrho_0):=\mathfrak{m}^{\ast}(v;\varrho_0)-2^{-(n+1)}\textrm{Osc}(v;\varrho_0).
\end{equation}
Then $\omega_0(v;\varrho_0)=\frac{1}{2}\left(\mathfrak{m}^{\ast}(v;\varrho_0)+\mathfrak{m}_{\ast}(v;\varrho_0)\right)$ and $\omega_n(v;\varrho_0)\nearrow\omega_{\infty}(v;\varrho_0):=\mathfrak{m}^{\ast}(v;\varrho_0)$. By replacing $v$ by $-v$ in (\ref{6.03}) if necessary, we may assume that $\omega_0(v;\varrho)\geq0$, and one can find $m\in\mathbb{N\!}$\, large enough, such that
\begin{equation}\label{6.13}
\left|\Omega_{k_m}(\varrho_0)\right|\,\leq\,\frac{1}{2}|\Omega(x_0;\varrho_0)|,\indent\,\forall\,\varrho_0\in[\rho/4,\rho/2],
\,\,\,\,\,\,\,\,\textrm{where}\,\,\,\,k_m:=\omega_m(v;\varrho_0).
\end{equation}
Applying Corollary \ref{consequences}(d) for numbers $\varrho_0,\,\sigma_0\in[\rho/4,\rho]$ such that $\varrho_0<\sigma_0$, we have that
$$|\Omega_h(\varrho_0)|^{^{\frac{r_0}{2}}}
\,\leq\,c\,(h-k)^{-2}(\sigma_0-\varrho_0)^{-2}
(|\Omega_k(\varrho_0)|-|\Omega_h(\varrho_0)|)
\left\|v_k\right\|^2_{_{2,\Omega_{k}(\sigma_0)}},$$
whenever $h>k\geq0$, where $k\geq0$ is such that
$|\Omega_k(\sigma_0)|\,\leq\,\frac{1}{2}|\Omega(x_0;\sigma_0)|$.
Since by (\ref{6.13}), this condition holds for $k=k_m:=\omega_m(v;\varrho_0)$, taking $\varrho_0=\rho/2$ and $\sigma_0=\rho$, using the fact that $\Omega$ is a $2$-set, and observing that $\rho^{-2}\geq(\rho^{-2})^{r_0/2}$, we deduce that
\begin{equation}
\label{6.14}|\rho^{-2}\Omega_h(\rho/2)|^{^{\frac{r_0}{2}}}
\leq\,c\,(h-k_m)^{-2}(\mathfrak{m}^{\ast}(v;\rho)-k_m)^{2}
\left\{\rho^{-2}|\Omega_{k_m}(\rho/2)|-\rho^{-2}|\Omega_h(\rho/2)|\right\}.
\end{equation}
Using (\ref{6.14}), we apply Lemma \ref{lemma3} to the function $\varphi(t):=\rho^{-2}|\Omega_t(\rho/2)|$ to obtain that\\
\begin{equation}\label{6.15}
\rho^{-2}|\Omega_{k_n}(\rho/2)|\stackrel{n\rightarrow\infty}{\longrightarrow}0\,\,\,\,\,\,\,\,\,\,\,\,\textrm{as}\,\,\,\,\,k_n:=
\omega_n(v;\rho/2)\stackrel{n\rightarrow\infty}{\longrightarrow}\mathfrak{m}^{\ast}(v;\rho).
\end{equation}
Therefore, (\ref{6.15}) implies that we can select $n$ sufficiently large (and independent of $x_0\in\overline{\Omega}$), such that
\begin{equation}
\label{6.16}|\Omega_{k_n}(\rho/2)|\,\leq\,\frac{1}{2}|\Omega(x_0;\rho/2)|\,\,\,\,\,\,\,\,\,\,\textrm{and}\,\,\,\,\,\,\,\,\,\,
C^{\star}\rho^{-N}|\Omega_{k_n}(\rho/2)|\,\leq\,\frac{1}{4}\,,
\end{equation}
where $C^{\star}>0$ denotes the constant appearing in (\ref{6.11}). Combining (\ref{6.11}) and (\ref{6.16}) for $k_0:=\omega_n(v;\rho/2)$, we deduce that\\[2ex]
$\mathfrak{m}^{\ast}(v;\rho/4)\,\leq\,\omega_n(v;\rho/2)+\left(C^{\ast}(\mathfrak{m}^{\ast}(v;\rho/2)-\omega_n(v;\rho/2))^2
\rho^{-2}\left|\Omega_{_{\omega_n(v;\rho/2)}}(\rho/2)\right|\right)^{1/2}$\\
\begin{equation}\label{6.17}
\leq\,\omega_n(v)+\frac{1}{2}\left(\mathfrak{m}^{\ast}(v;\rho/2)-\omega_n(v;\rho/2)\right)
\,\leq\,\mathfrak{m}^{\ast}(v;\rho)-2^{-(n+2)}\textrm{Osc}(v;\rho).
\end{equation}
Inserting the property $\mathfrak{m}_{\ast}(v;\rho/4)\,\geq\,\mathfrak{m}_{\ast}(v;\rho)$ into (\ref{6.17}), we calculate that
\begin{equation}\label{6.18}
\textrm{Osc}(v;\rho/4)\,\leq\,\mathfrak{m}^{\ast}(v;\rho)-2^{-(n+2)}\textrm{Osc}(v;\rho)-\mathfrak{m}_{\ast}(v,\rho)
\,\leq\,(1-2^{-(n+2)})\,\textrm{Osc}(v;\rho).
\end{equation}
Letting $\varsigma_0:=(1-2^{-(n+2)})\in(0,1)$, and recalling that $\rho\in(0,\rho^{\ast})\subseteq(0,1)$ was chosen arbitrarily,
we are lead to the inequality (\ref{6.04}), completing the proof.
\end{proof}

\subsection{Proof of Theorem \ref{Holder-continuity}}

Given $u\in\mathcal{V}_2(\Omega)$ a weak solution of (\ref{E01}), take $x_0\in\overline{\Omega}$ and $0<\rho\leq\rho^{\ast}<1$, let $\hat{u}\in\mathcal{V}_c(\Omega;D)$ solve (\ref{c-Existence-Elliptic}) for
$D=B(x_0;\rho)$, and put $v:=u-\hat{u}\in\mathcal{V}_2(\Omega)$. Then $u,\,\hat{u},\,v$ are all globally bounded by virtue of Theorem \ref{Elliptic-Infinity-Case} with $\hat{u}\in\mathcal{V}_c(\Omega;D)$ fulfilling (\ref{Thm0-01b}), and moreover $v\in\mathcal{V}_2(\Omega)$ solves (\ref{6.03}). Therefore inequality (\ref{6.04}) holds for the function $v$, and consequently combining all these facts, we are lead to the inequality
\begin{equation}\label{Holder01}
\textrm{Osc}(u;\rho/4)\,\leq\,2|\|\hat{\mathbf{u}}\||_{_{\mathbb{X\!}^{\,\infty}(\Omega;\Gamma^{\infty})}}+\varsigma_0\,\textrm{Osc}(v;\rho)\,\leq\,4C^{\star}_1\rho^{\sigma^{\star}}\left(\displaystyle\sum^2_{i=0}\|f_i\|_{_{p_i,\Omega}}+\|g\|_{_{q,\Gamma^{\infty}}}\right)+\varsigma_0\,\textrm{Osc}(u;\rho),
\end{equation}
where the positive constants $C^{\star}_1$,\, $\sigma^{\star}$,\, $\varsigma_0$ are defined in Theorem \ref{Elliptic-Infinity-Case}(b) and Lemma \ref{oscilation}, respectively. Therefore, one applies Lemma \ref{lemma4} to the function $\Phi(t):=\textrm{Osc}(u;t)$ ($0<t\leq\rho^{\ast}$) to conclude that
\begin{equation}
\label{Holder02}\textrm{Osc}(u;\rho)\,\leq\,C
\left[\left(\frac{\rho}{\rho^{\ast}}\right)^{\vartheta}\textrm{Osc}(u;\rho^{\ast})+
\left(\displaystyle\sum^2_{i=0}\|f_i\|_{_{p_i,\Omega}}+\|g\|_{_{q,\Gamma^{\infty}}}\right)\rho^{\sigma^{\star}\delta}(\rho^{\ast})^{\sigma^{\star}(1-\delta)}\right],
\end{equation}
for all $\rho\in(0,\rho^{\ast}/4]$, for every $\delta\in(0,1)$, and for some constants $C,\,\vartheta>0$ (as in Lemma \ref{lemma4}). This shows that $u$ is H\"older continuous over $\overline{\Omega}(x_0;\rho/4)$, and the compactness of $\overline{\Omega}$ gives the existence of a constant $\delta_0\in(0,1)$ (independent of $u$) such that $u\in C^{0,\delta_0}(\overline{\Omega})$. Finally, a combination of (\ref{Holder02} with (\ref{Thm0-01}) gives (\ref{Holder-norm}), completing the proof. $\hfill\square$

\subsection{Inverse positivity}

\indent We conclude this section with the following result which will be of importance later on in the time-dependent case.

\begin{theorem}\label{Inverse-positivity}
Under the conditions of Theorem \ref{Holder-continuity}, let $u\in\mathcal{V}_2(\Omega)$ be a weak solution of (\ref{E01}). If $\Phi(w)\geq0$ for each $w\in\mathcal{V}_2(\Omega)^+$ (where the functional $\Phi(\cdot)$ is given by (\ref{Functional})), then $u(x)\geq0$ for all $x\in\overline{\Omega}$.
\end{theorem}

\begin{proof}
Since $\mathcal{E}(u^+,u^-)=0$, using the coercivity of the form we easily get that
$$c\,|\|\mathbf{u}^-\||^2_{_{\mathbb{X\!}^{\,2}(\Omega;\Gamma)}}\,\leq\,\mathcal{E}(u^-,u^-)=-\mathcal{E}(u,u^-)=\Phi(u^-)\,\leq\,0.$$
Thus $u^-=0$ over $\overline{\Omega}$, and this together with Theorem \ref{Holder-continuity} yield the desired result.
\end{proof}

\begin{remark}\label{pre-fractal-Holder}
In the case of the pre-fractal domains $\Omega^{(m)}$ ($m\in\mathbb{N\!}\,$), since these domains are Lipschitz, problem (\ref{E01}) can be defined in the classical way, and in particular from \cite{NITTKA}, one sees that the conclusions of Theorem \ref{Holder-continuity} and Theorem \ref{Inverse-positivity} hold for each set $\Omega^{(m)}$. Furthermore, adapting the proof in \cite[Theorem 4]{ACH10}, one can deduce that if $u^{(m)}\in\mathcal{V}_2(\Omega^{(m)})$ are weak solutions of problem (\ref{E01}) over the pre-fractal sets $\Omega^{(m)}$ for each $m\in\mathbb{N\!}_{\,0}$, and $u\in\mathcal{V}_2(\Omega)$ solves (\ref{E01}) over the ramified domain, then
$$\displaystyle\lim_{m\rightarrow\infty}\left\|u^{(m)}-u|_{_{\Omega^{(m)}}}\right\|_{_{H^1(\Omega^{(m)})}}\,=\,0.$$
\end{remark}

\section{The diffusion equation}\label{sec6}

\indent In this section we will consider the realization of the corresponding inhomogeneous heat equation
\begin{equation}
\label{D01}\left\{
\begin{array}{lcl}
\displaystyle\frac{\partial u}{\partial t}-\mathcal{A}u+\mathcal{B}u\,=\,f(t,x)\,\,\,\,\,\,\,\,\,\,\,\textrm{in}\,\,(0,\infty)\times\Omega;\\
\,\,\,\,u\,=\,0\,\,\,\,\,\,\,\,\,\,\,\,\,\,\,\,\,\,\,\,\,\,\,\,\,\,\,\,\,\,\,\,\,\,\,\,\,\,\,\,\,\textrm{on}\,\,(0,\infty)\times(\Gamma\setminus\Gamma^{\infty});\\
\frac{\partial u}{\partial\nu_{_{\mathcal{A}}}}+\beta u\,=\,g(t,x)\,\,\,\,\,\,\,\,\,\,\,\,\,\,\,\,\,\,\,\,\,\,\,\,\,\,\,\,\textrm{on}\,\,(0,\infty)\times\Gamma^{\infty},\\
\,\,\,\,u(0,x)\,=\,u_0(x)\,\,\,\,\,\,\,\,\,\,\,\,\,\,\,\,\,\,\,\,\,\,\,\,\,\,\,\textrm{in}\,\,\Omega,\\
\end{array}
\right.
\end{equation}
for given functions $f\in L^2(I;L^2(\Omega))$ and $g\in L^2(I;L^2_{\mu}(\Gamma^{\infty}))$, where $I$ denotes either $[0,\infty)$, or $[0,T]$ for some constant $T>0$. All the assumptions here are the same as in Theorem \ref{Holder-continuity}.

\subsection{Strong Feller resolvent and semigroup}\label{subsec6.1}

\indent In this part we apply the results in the previous section to derive a strong Feller resolvent and a corresponding Feller semigroup
associated with the bilinear form $\mathcal{E}(\cdot,\cdot)$ given by (\ref{E-form}). This form is in general non-symmetric nor sub-Markovian, unless $b_i=0$ for each $i\in\{1,2\}$ (Cf. \cite[Corollary 2.17]{OUHA05}).\\
\indent Let $A_{\mu}$ denote the unique operator associate with the form $\mathcal{E}(\cdot,\cdot)$, and
let $\{T_{\mu}(t)\}_{t\geq 0}$ denote the corresponding $C_0$-semigroup generated by $-A_{\mu}$ in $L^2(\Omega)$. Clearly $\{T_{\mu}(t)\}_{t\geq 0}$ is analytic  and compact (Cf. \cite[Chapter 17, Section 6]{DautrayLions}). Consequently, for each $u_0\in L^2(\Omega)$, the function $u(t):=T_{\mu}(t)u_0$ is the unique mild solution of the abstract Cauchy problem
\begin{equation}\label{parabolicpbomog}
\left\{
     \begin{array}{ll}
       u_t= A_{\mu}u\,\,\,\,\,\,\,\,\,\,\,\,\,\,\,\,\,\,\,\,\,\,\textrm{for}\,\,t\in(0,\infty);\\
       u(0)=u_0\,\,\,\,\,\,\,\,\,\,\,\,\,\,\,\,\,\,\,\,\,\,\textrm{in}\,\,\Omega
     \end{array}
   \right.
\end{equation}
Furthermore, proceeding in an analogous way as in the derivation of \cite[Theorem 6.8]{OUHA05}, one has that
\begin{equation}\label{ultracontractivity}
\|T_{\mu}(t)\|_{_{\mathcal{L}(L^2(\Omega);L^\infty(\Omega))}}\,\leq\, c_1 t^{-1/2}e^{c_2t}\,\,\,\,\,\,\textrm{for all}\,\,t>0,
\end{equation}
and thus the semigroup $\{T_{\mu}(t)\}_{t\geq 0}$ is ultracontractive in the sense that $T_{\mu}(t)$ maps $L^2(\Omega)$ into $L^{\infty}(\Omega)$.
In addition, by working with the adjoint semigroup and performing a duality argument, we infer that
$$\|T_{\mu}(t)\|_{\mathcal{L}(L^1(\Omega);L^\infty(\Omega))}\leq c'_1 t^{-1}e^{c'_2t}\,\,\,\,\,\,\textrm{for each}\,\, t>0.$$
Therefore, by Dunford Pettis's Theorem, the semigroup $\{T_{\mu}(t)\}_{t\geq 0}$
admits an integral representation with kernel $K_{\mu}(x,y,t)\in L^\infty(\Omega\times \Omega)$ (e.g. \cite[Theorem 1.3]{Arendt-Bukl94}), which satisfies Gaussian estimates, and consequently in views of \cite{AR-ELST97},
it follows that $\{T_{\mu}(t)\}_{t\geq 0}$ extrapolates
to a family of compact analytic semigroups over $L^p(\Omega)$ for all $p\in[1,\infty]$, which are strongly continuous whenever $p\in[1,\infty)$,
and all have the same angle of analyticity. Furthermore, one has that $\{T_{\mu}(z)\}_{_{\textrm{Re}(z)>0}}$ are kernel operators satisfying Gaussian estimates.
In particular, for each $u_0\in L^p(\Omega)$ with $1\leq p<\infty$, the function $u(t):=T_{\mu}(t)u_0\in L^p(\Omega)$ is the unique mild solution
to problem (\ref{parabolicpbomog}).
Moreover, following the approach as in \cite[Proposition 7.1]{Daners}, we find the existence of constants $M\geq 0$ and $\omega\in \mathbb{R\!}\,$ such that
\begin{equation}\label{infinity-semigroup}
\|u(t)\|_{_{\infty,\Omega}}\,\leq\, Me^{\omega t}\|u_0\|_{_{\infty,\Omega}}
\end{equation}
for each $u_0\in L^{\infty}(\Omega)$.\\
\indent Next, given $\gamma>0$, let $G_{\gamma}$ denote the resolvent associated with the form $\mathcal{E}(\cdot,\cdot)$. Recall that here and in the rest of the paper, we are assuming that the condition (\ref{sharp-coercivity}) holds. Then, we first establish the following result.

\begin{theorem}\label{feller-res}
If $f\in L^p(\Omega)$ for $p>1$, then for each $\gamma>0$, one has $G_{\gamma}(L^p(\Omega))\subseteq C^{0,\delta_0}(\overline{\Omega})$
for some $\delta_0\in(0,1)$. In particular, $G_{\gamma}$ is a strong Feller resolvent.
\end{theorem}

\begin{proof}
Given $f\in L^p(\Omega)$ for $p>1$, put $u_{\gamma}:=G_{\gamma}(f)$. Then $u_{\gamma}$ is a weak solution of (\ref{E01}), and thus from here, Theorem \ref{Holder-continuity} implies that $u_{\gamma}\in C^{0,\delta_0}(\overline{\Omega})$ for some constant $\delta_0\in(0,1)$. Hence $G_{\gamma}$ is a Feller resolvent.
\end{proof}

\begin{remark}\label{Kernel}
From Theorem \ref{feller-res} together with the arguments in \cite[proof of Theorem 4.1]{NITTKA}, one clearly see that the semigroup $\{T_{\mu}(t)\}_{t\geq 0}$ has the Feller property, and the corresponding integral kernel $K_{\mu}(\cdot,\cdot,\cdot)\in C^{0,\delta_0}([a,b]\times\overline{\Omega}\times\overline{\Omega})$ whenever $0<a\leq b<\infty$.
\end{remark}

\indent Given $A_{\mu}$ the operator associated with the form $\mathcal{E}(\cdot,\cdot)$, denote by
$A^c_{\mu}$ the part of the operator $A_{\mu}$ in $C(\overline{\Omega})$, in the sense that
$$D(A^c_{\mu}):=\{w\in D(A_{\mu})\cap C(\overline{\Omega})\mid A_{\mu}w\in C(\overline{\Omega})\},
\,\,\,\,\,\,\,\,\,\,\,\,A^c_{\mu}w=A_{\mu}w.$$ By Theorem \ref{feller-res}, $A^c_{\mu}$ is well-defined. Moreover, we now establish a key density result which cannot be proven in the classical way as in \cite{NITTKA,WAR06}, since in our case one cannot give a concrete description of the notion of generalized derivative (as in Definition \ref{Def-gen-normal}), and we do not know if it is possible to approximate a given function in $L^s_{\mu}(\Gamma^{\infty})$ by the generalized co-normal derivative of suitable functions. Therefore, motivated by some results in \cite{NITTKA2010}, we attack this problem through another direction.

\begin{theorem}\label{feller-density}
The part $A^c_{\mu}$ of the operator $A_{\mu}$ in $C(\overline{\Omega})$ is densely defined over $C(\overline{\Omega})$.
\end{theorem}

\begin{proof}
We first recall that if $f_i\in L^{p_i}(\Omega)$ for $p_0>1$ and $\displaystyle\min_{1\leq i\leq2}\{p_i\}>2$, if one consider the elliptic equation (\ref{E01}) for $g\equiv0$, recalling condition (\ref{sharp-coercivity}), from Theorem \ref{exist-uniq-weak-soln} and Theorem \ref{Holder-continuity} we get that there exists a unique weak solution $u_0\in\mathcal{V}_2(\Omega)\cap C^{0,\delta_0}(\overline{\Omega})$ of problem (\ref{E01}) (for $g\equiv0$), fulfilling the Schauder norm (\ref{Holder-norm}). In particular, one sees that $u_0$ depends continuously in $C(\overline{\Omega})$ on the data $f_i$ in $L^{p_i}(\Omega)$ ($i\in\{0,1,2\}$).\\
\indent Now, since $C^{\infty}(\overline{\Omega})$ is dense in $C(\overline{\Omega})$ (due to the Stone-Weierstrass theorem), it suffices to show that $C^{\infty}(\overline{\Omega})$ lies in the closure of $D(A^c_{\mu})$. Take $v\in C^{\infty}(\overline{\Omega})$, and notice that $A_{\mu}v\in H^1(\Omega)^{\ast}$. Recalling that $\langle A_{\mu}v,w\rangle=\mathcal{E}(v,w)$ and setting $\vartheta:=\min\{\zeta,r,s,2-\epsilon\}\in (1,2)$ for $\epsilon\in(0,1)$, given $w\in H^1(\Omega)$, a closely examination of the structure of $\mathcal{E}(\cdot,\cdot)$ gives that
$$|\langle A_{\mu}v,w\rangle|\,\leq\,\|\bar{\alpha}\|_{_{\infty,\Omega}}\|\nabla v\|_{_{\infty,\Omega}}\|\nabla w\|_{_{1,\Omega}}+\|\bar{\eta}\|_{_{\vartheta,\Omega}}\|\nabla v\|_{_{\infty,\Omega}}\|w\|_{_{\frac{\vartheta}{\vartheta-1},\Omega}}+\|\lambda\|_{_{\vartheta,\Omega}}\|v\|_{_{\infty,\Omega}}\|w\|_{_{\frac{\vartheta}{\vartheta-1},\Omega}}+\|\beta\|_{_{\vartheta,\Gamma^{\infty}}}\|v\|_{_{\infty,\Gamma^{\infty}}}\|w\|_{_{\frac{\vartheta}{\vartheta-1},\Gamma^{\infty}}}$$
$$\indent\indent\indent\indent\leq\,C_0\max\left\{\|\bar{\alpha}\|_{_{\infty,\Omega}}\|\nabla v\|_{_{\infty,\Omega}},\,\|\bar{\eta}\|_{_{\vartheta,\Omega}}\|\nabla v\|_{_{\infty,\Omega}},\,\|\lambda\|_{_{\vartheta,\Omega}}\|v\|_{_{\infty,\Omega}},\,\|\beta\|_{_{\vartheta,\Gamma^{\infty}}}\|v\|_{_{\infty,\Gamma^{\infty}}}\right\}\|w\|_{_{W^{1,\varsigma}(\Omega)}},$$ for some constant $C_0>0$, where $\varsigma:=2\vartheta'(d+\vartheta')^{-1}\in[1,2)$ (recall that $d:=-\log(2)/\log(\tau)\in[1,2)$ and $\vartheta'$ denotes the conjugate of $\vartheta$). Also, in the last inequality, we have recalled Theorem \ref{prop-ramified}(c,d) and have used the corresponding linear continuous interior and trace Sobolev maps following an approach analogous as in the proof of Theorem \ref{embeddings-trace}. Thus, we have that the functional $\Psi_v(\cdot):=\langle A_{\mu}v,\cdot\rangle$ extends to a bounded linear functional on $W^{1,\varsigma}(\Omega)$, so $A_{\mu}v\in W^{1,\varsigma}(\Omega)^{\ast}$. From here, an application of \cite[Theorem 4.3.3]{ZIE89} implies the existence of functions $\tilde{f}_i\in L^{^{\frac{\varsigma}{\varsigma-1}}}(\Omega)$ ($i\in\{0,1,2\}$), such that
$$\langle A_{\mu}v,w\rangle=\displaystyle\int_{\Omega}\tilde{f}_0w\,dx+\displaystyle\int_{\Omega}\displaystyle\sum^2_{i=1}\tilde{f}_i\partial_{x_i}w\,dx,\,\,\,\,\,\,\,\,\,\,\,\,\,\,\,\,\textrm{for each}\,\,\,w\in H^1(\Omega).$$ Hence $v$ solves (\ref{E01}) for $g\equiv0$ and $f_i$ replaced by $\tilde{f}_i$. Now, for $i\in\{0,1,2\}$, choose sequences $\{\tilde{f}_{n,i}\}_{_{n\in\mathbb{N\!}}}\subseteq C^{\infty}(\overline{\Omega})$ fulfilling $\tilde{f}_{n,i}\stackrel{n\rightarrow\infty}{\longrightarrow}\tilde{f}_{i}$ in $L^{^{\frac{\varsigma}{\varsigma-1}}}(\Omega)$, and let $\{v_n\}_{_{n\in\mathbb{N\!}}}$ denote the corresponding sequence of weak solutions to problem (\ref{E01}) with homogeneous boundary conditions related to the data sequences $\{\tilde{f}_{n,i}\}_{_{n\in\mathbb{N\!}}}$. Since $\varsigma(\varsigma-1)^{-1}>2$, Theorem \ref{Holder-continuity} shows that $\{v_{n}\}_{_{n\in\mathbb{N\!}}}\subseteq D(A^c_{\mu})$. By the continuity of the solution operator together with the linear dependence of the function with respect to the data, it follows that $v_{n}\stackrel{n\rightarrow\infty}{\longrightarrow}v$ in $C(\overline{\Omega})$. Thus, $v\in\overline{D(A^c_{\mu})}$, proving that $C^{\infty}(\overline{\Omega})\subseteq\overline{D(A^c_{\mu})}$, as desired.
\end{proof}

\subsection{Proof of Theorem \ref{feller-semi}}\label{proof-main-2}

By virtue of Theorem \ref{feller-res}, we see that $\mathcal{A}^c_{\mu}$ generates a $C_0$-semigroup
$T^c_{\mu}(t)$ over $C(\overline{\Omega})$, which in views of the preceding arguments, is analytic of angle $\pi/2$. Furthermore, $T^c_{\mu}(t)$
is a compact operator, since one can write
$$T^c_{\mu}(t)=T_4\circ T_3\circ T_2\circ T_1,$$
for $T_1:=Id:C(\overline{\Omega})\rightarrow L^2(\Omega)$ is bounded, $T_2:=T_{\mu}(t/3):L^2(\Omega)\rightarrow L^2(\Omega)$ is compact,
$T_3:=T_{\mu}(t/3):L^2(\Omega)\rightarrow L^{\infty}(\Omega)$ is bounded (since the semigroup is ultracontractive), and
$T_4:=T_{\mu}(t/3):L^{\infty}(\Omega)\rightarrow C(\overline{\Omega})$ is bounded (since the semigroup has the Feller property). This completes the proof. $\hfill\square$

\subsection{Solvability of the inhomogeneous diffusion equation}\label{subsec6.2}
\indent The aim of this subsection is to establish the well-posedness of the inhomogeneous Robin-type diffusion equation (\ref{D01}). We will adopt most of the tools employed by Nittka \cite[section 2]{NITTKA2014}, modifying the required ones in order to include the case in this paper over our highly irregular domain.\\
\indent We begin our discussion by posing the general definition for a weak solution of problem (\ref{D01}).

\begin{definition}\label{weak-sol-parabolic}
 Given $T>0$, we say that a function $u\in C([0,T];L^2(\Omega))\cap L^2((0,T); \mathcal{V}_2(\Omega))$ is a \textbf{weak solution of the
 parabolic equation (\ref{D01}) on $[0,T]$}, if\\[2ex]
\,\,\,\,\,\,$-\displaystyle\int_{0}^{T}\displaystyle\int_{\Omega} u(\xi)\psi_t(\xi)\, dx d\xi +\int_{0}^{T}\mathcal{E}(u(\xi),\psi(\xi))\, d\xi$\\
\begin{equation}
=\displaystyle\int_{\Omega} u_0 \psi(0) dx+ \int_{0}^{T}\int_{\Omega} f(\xi,x)\psi(\xi)\,dx d\xi +\displaystyle\int_{0}^{T}\int_{\Gamma} g(\xi,x) \psi(\xi) d\mu d\xi,
\end{equation}
 \noindent for all $\psi\in H^1(0,T;\mathcal{V}_2(\Omega))$ such that $\psi(T)=0$, where $u_0:=u(0,x)$ and $\mathcal{E}(\cdot,\cdot)$
 denotes the bilinear form given by (\ref{E-form}).
 Similarly, we say that a function $u:[0, \infty)\to L^2(\Omega)$ is a weak solution of (\ref{D01}) on $[0,\infty)$ if for every $T>0$
 its restriction to $[0,T]$ is a weak solution on $[0,T].$
 \end{definition}

\indent Transitioning from the homogeneous problem into the inhomogeneous case via quadratic form's methods is not a immediate fact, and requires a clever operator technique which we present below (as in \cite{NITTKA2014}).

\begin{definition}\label{a-normal-a}
Given $u\in \mathcal{V}_2(\Omega)$, let $\mathcal{A}$ and $\mathcal{B}$ be the linear operators defined by (\ref{E02}).
\begin{enumerate}
\item[(a)]\,\, We say that $\mathcal{L}u:=(\mathcal{A}+\mathcal{B})u\in L^2(\Omega)$ if there exists a $f\in L^2(\Omega)$ such that
$$\mathcal{E}_0(u,v)=-\langle f,v\rangle_{_{L^2(\Omega)}},\,\,\,\,\,\,\,\,\,\,\,\textrm{for all}\,\,v\in H^1_0(\Omega),$$
where $\mathcal{E}_0(u,v):=\mathcal{E}(u,v)-\displaystyle\int_{\Gamma^{\infty}}\beta uv\,d\mu$.
\item[(b)]\,\, Assume that $\mathcal{L}u\in L^2(\Omega)$. Then we say that $\frac{\partial u}{\partial\nu_{_{\mathcal{A}}}}\in L^2_{\mu}(\Gamma^{\infty})$, if there exists a unique function $g\in L^2_{\mu}(\Gamma^{\infty})$ such that
$$\langle-\mathcal{L}u, v\rangle_{_{L^2(\Omega)}}+\langle g, v\rangle_{_{L^2_{\mu}(\Gamma^{\infty})}}=\mathcal{E}_0(u,v),\,\,\,\,\,\,\,\,\,\,\textrm{for every}\,\,v\in \mathcal{V}_2(\Omega).$$
\item[(c)]\,\, We define the operator $\mathbf{A}_{\mu}:\mathbb{X\!}^{\,2}(\Omega;\Gamma)\rightarrow L^2(\Omega)\times L^2_{\mu}(\Gamma^{\infty})$ by:
$$\left\{
\begin{array}{lcl}
D(\mathbf{A}_{\mu}):=\left\{(u,0)\mid u \in \mathcal{V}_2(\Omega),\,\,\,\mathcal{L}u\in L^{2}(\Omega),\,\,\,\frac{\partial u}{\partial\nu_{_{\mathcal{A}}}}\in L^{2}_{\mu}(\Gamma^{\infty})\right\};\\[1ex]
\mathbf{A}_{\mu}(u,0):=\left(\mathcal{L}u,\,-\frac{\partial u}{\partial\nu_{_{\mathcal{A}}}}-\beta u\right).
\end{array}
\right.$$
 \end{enumerate}
 \end{definition}

 \begin{remark}\label{caratt- D-A2}
 In views of Definition \ref{a-normal-a}(c), it is easily verified that $(u,0)\in D(\mathbf{A}_{\mu})$ with  $-\mathbf{A}_{\mu}(u,0)=(f,g)$ if and only if
 $$\mathcal{E}(u,v) =\langle f,v\rangle_{_{L^2(\Omega)}}+ \langle g,v\rangle_{_{L^2_{\mu}(\Gamma^{\infty})}}\,\,\,\,\,\,\,\,\,\,\textrm{for every}\,\,\,v\in
 \mathcal{V}_2(\Omega).$$
 Therefore if $(f,g)\in\mathbb{X\!}^{\,2}(\Omega;\Gamma^{\infty})$, then
$(u,0)\in D(\mathbf{A}_{\mu})$ with  $-\mathbf{A}_{\mu}(u,0)=(f,g)$ if and only if $u\in\mathcal{V}_2(\Omega)$ is a weak solution of the boundary value problem (\ref{E01}).
 \end{remark}

 We now have some basic properties for the operator $\mathbf{A}_{\mu}$.

 \begin{prop}\label{A-operator-prop}
 The operator $\mathbf{A}_{\mu}$ is resolvent positive. More precisely, the operator
 $\gamma-\mathbf{A}_{\mu}$ on $D(\mathbf{A}_{\mu})$ is invertible for all $\gamma>0$, and if $\mathbf{A}_{\mu}(u,0)=(f,g)$ for $(f,g)\in L^{2}(\Omega)\times L^{2}_{\mu}(\Gamma)$ fulfilling
 $\langle f,w\rangle_{_{L^2(\Omega)}}+\langle g,w\rangle_{_{L^2_{\mu}(\Gamma^{\infty})}}\geq0$ for each $w\in\mathcal{V}_2(\Omega)^+$, then $u\geq 0$ a.e. Moreover,
if $D(\mathbf{A}_{\mu})$ is equipped with the graph norm, then $D(\mathbf{A}_{\mu})$ is continuously embedded into $\mathcal{V}_2(\Omega)\times \{0\}.$
 \end{prop}

 \begin{proof}
 If  $\gamma>0$, then by virtue of (\ref{sharp-coercivity}) together with Lax-Milgram's Lemma, for each $(f,g)\in L^{2}(\Omega)\times L^{2}_{\mu}(\Gamma^{\infty})$ , there exists a unique function $u\in {\mathcal{V}_2(\Omega)}$ such that
$$\gamma\langle u,v\rangle_{_{L^2(\Omega)}}+\mathcal{E}(u,v)=\langle f,v\rangle_{_{L^2(\Omega)}} +\langle g,v\rangle_{_{L^2_{\mu}(\Gamma^{\infty})}},
\,\,\,\,\,\,\,\,\,\,\,\,\,\,\,\textrm{for all}\,\,v\in \mathcal{V}_2(\Omega).$$
This together with Remark \ref{caratt- D-A2} we have that $(\gamma-\mathbf{A}_{\mu})(u,0)=(f,g)$, which implies
that the operator $\gamma-\mathbf{A}_{\mu}: D(\mathbf{A}_{\mu})\rightarrow L^{2}(\Omega)\times L^{2}_{\mu}(\Gamma^{\infty})$ is a bijection for each
$\gamma>0$. Furthermore, if $\langle f,w\rangle_{_{L^2(\Omega)}}+\langle g,w\rangle_{_{L^2_{\mu}(\Gamma^{\infty})}}\geq0$ for each $w\in\mathcal{V}_2(\Omega)^+$, then from the proof of Theorem \ref{Inverse-positivity} one clearly deduces that $u\geq 0$ a.e. in $\Omega$. Hence the resolvent $(\gamma-\mathbf{A}_{\mu})^{-1}$ is a positive operator, and since every positive operator is continuous, we get that $\gamma-\mathbf{A}_{\mu}$ is invertible. Moreover, the closedness of $\mathbf{A}_{\mu}$ implies that $D(\mathbf{A}_{\mu})$ is a Banach space with respect to the graph norm. Since the embedding $D(\mathbf{A}_{\mu})\hookrightarrow\mathbb{X\!}^{\,2}(\Omega;\Gamma)$ and the embedding $\mathcal{V}_2(\Omega)\times \{0\}\hookrightarrow\mathbb{X\!}^{\,2}(\Omega;\Gamma)$ are both continuous, an application of the Closed Graph Theorem implies that the embedding
$D(\mathbf{A}_{\mu})\hookrightarrow\mathcal{V}_2(\Omega)\times \{0\}$ is continuous, as asserted.
 \end{proof}

\indent Although the operator $\mathbf{A}_{\mu}$ defined above exhibits great properties as one sees in Proposition \ref{A-operator-prop}, it has the disadvantage that its domain is not dense in $\mathbb{X\!}^{\,2}(\Omega;\Gamma)$, and consequently it is not a generator of a $C_0$-semigroup. However, one may avoid this issue by considering integrated semigroups.\\
\indent The remaining of the results in this section can be established in the exact way as in \cite[section 2]{NITTKA2014} (with some even simpler steps and minor modifications), so we will omit all the proofs of the results, and will only state the relevant definitions and results. We first introduce the notion of mild and classical solution of problem (\ref{D01}), recalling that from now on we will assume that $f\in L^2(0,T; L^2(\Omega))$ and $g\in L^2(0,T,L^2_{\mu}(\Gamma))$.\\

 \begin{definition}\label{classical sol}
Let $I=[0,T]$ or $I=[0,\infty)$.
\begin{enumerate}
\item[(a)]\,\, We say that function $u$ is a \textbf{classical solution  of problem (\ref{D01})}, if
$u\in C^1((0,T], L^2(\Omega))$ with  $u(0)=u_0,$  the mapping $t\to (u(t),0)$ lies in $C(I,D(\mathbf{A}_{\mu}))$, and
$$(u_t(\xi),0)- (\mathbf{A}_{\mu}u(\xi),0)=(f(\xi,x), g(\xi,x)),\,\,\,\,\,\,\,\,\textrm{for all}\,\,\xi\in I.$$
\item[(b)]\,\,We call $u$ is a \textbf{mild solution of problem (\ref{D01})},
if $u\in  C(I, L^2(\Omega)))$ with  $(\int_{0}^{t}u(\xi)\,d\xi, 0) \in D(\mathbf{A}_{\mu})$ for all $t\geq 0$, and
\begin{center}
$(u(t), 0)-\mathbf{A}_{\mu}(\int_{0}^{t}u(\xi)\,d\xi, 0)=(\int_{0}^{t} f(\xi,x)\,d\xi,\int_{0}^{t} g(\xi,x)\,d\xi)$.
\end{center}
\end{enumerate}
 \end{definition}

\indent Given $I=[0,T]$ or $I=[0,\infty)$, it follows that every classical solution of problem (\ref{D01}) on $I$ is a weak solution of problem (\ref{D01} on $I$, and
every weak solution of equation (\ref{D01}) on $I$ is a mild solution of the heat equation (\ref{D01}) on $I$. Moreover, if $(u_0,0)\in D(\mathbf{A}_{\mu}),$\,  $\mathbf{A}_{\mu}(u_0,0)+(f(0,x),g(0,x)) \in D(\mathbf{A}_{\mu})$,\,
$f\in C^2([0,T], L^2(\Omega))$, and $g\in C^2([0,T], L^2_{\mu}(\Gamma^{\infty}))$, then equation (\ref{D01}) has a classical solution (Cf. \cite[proof of Theorem 2.6 and Proposition 2.7]{NITTKA2014}). We also have the following important inequality.

\begin{prop}\label{apriori estimates}
Given $T>0$, if $u$ is a classical solution of (\ref{D01}) on $[0,T]$, then there exists a constant $C>0$ (depending on $T$, \,$|\Omega|$, and the measurable coefficients of the form $\mathcal{E}(\cdot,\cdot)$) such that
\begin{equation}\label{aux-estimate}
\displaystyle\sup_{t\in [0,T]}\displaystyle\int_{\Omega} |u(t)|^2 dx +\displaystyle\int_{0}^{T}\|\nabla u(\xi)\|^2_{_{2,\Omega}}\,d\xi\,
\leq\,C\left(\|u_0\|^2_{_{2,\Omega}}+\displaystyle\int_{0}^{T}\|f(\xi)\|^2_{_{2,\Omega}}d\xi+\displaystyle\int_{0}^{T}\|g(\xi)\|^2_{_{2,\Gamma^{\infty}}}d\xi\right).
\end{equation}
Furthermore, by approximating a weak solution by classical solutions, it follows that (\ref{aux-estimate}) is valid for $u\in C([0,T];L^2(\Omega))\cap L^2((0,T); \mathcal{V}_2(\Omega))$ a weak solution of problem (\ref{D01}).
\end{prop}

\indent From here, we have the existence and uniqueness result for weak solutions to problem (\ref{D01}).

\begin{theorem}\label{exist-uniq-weak-soln}
Given $T>0$ arbitrary fixed, let $u_0\in L^2(\Omega)$,\, $f\in L^2(0,T; L^2(\Omega))$ and $g\in L^2(0,T,L^2_{\mu}(\Gamma^{\infty}))$.
Then there exists a weak solution of (\ref{D01}) in $[0,T]$ which is unique, even in the class of mild solutions. Moreover, the notions of weak and mild solution coincide.
\end{theorem}

\subsection{Fine regularity for weak solutions}\label{subsec6.3}

\indent In this part we are devoted to establish that weak solutions of the diffusion equation (\ref{D01}) are globally bounded under minimal assumptions. Our approach is highly motivated and follows the refined procedures developed by Nittka \cite{NITTKA2014}, whose motivations are based on results in \cite{LAD-SOL-URAL68,LAD-URAL68}.\\
\indent To begin, we first comment that for convenience (and without loss of generality), we will work with negative time. Then, following the same approach as in subsection \ref{subsec5.2}, given $u\in L^{\infty}(-T,0;L^2(\Omega))\cap L^2(-T,0;H^1(\Omega))$ and $k\geq0$ a fixed real number, we put
\begin{equation}\label{2.16a}
u_{k}(t):=(u(t)-k)^+.
\end{equation}
Clearly $u_{k}(t)\in L^{\infty}(-T,0;L^2(\Omega))\cap L^2(-T,0;H^1(\Omega))$,
with $\partial_{x_i}u_{k}(t)=(\partial_{x_i}u(t))\chi_{_{\Omega_k(t)}}$ for each $i\in\{1,2\}$, where
\begin{equation}\label{2.17a}
E_k(t):=\{x\in\overline{\Omega}\mid u(t)\geq k\}\,\,\,\,\,\,\,\,\textrm{and}\,\,\,\,\,\,\,\,D_k(t):=D\cap E_k(t),
\end{equation}
for each $D\subseteq\mathbb{R}^N$ such that $D\cap\overline{\Omega}\neq\emptyset$. Then for convenience, given $\varsigma>0$, we write
\begin{equation}\label{Special-Norm}
\Upsilon^2_{\varsigma}(u):=\displaystyle\sup_{-\varsigma\leq t\leq0}\displaystyle\int_{\Omega}|u(t)|^2\,dx+
\displaystyle\int^0_{-\varsigma}\displaystyle\int_{\Omega}|\nabla u(t)|^2\,dx\,dt.
\end{equation}
By virtue of (\ref{2.03}) and (\ref{2.04}), we get the following parabolic Sobolev inequality:
\begin{equation}\label{Sobolev-parabolic}
\|u\|_{_{L^{\kappa_p}(-\varsigma,0;L^{p}(\Omega))}}+\|u\|_{_{L^{\kappa_q}(-\varsigma,0;L^{q}_{\mu}(\Gamma^{\infty}))}}\,\leq\,C_0\,\Upsilon^2_{\varsigma}(u)
\end{equation}
for all $u(t)\in L^{\infty}(-T,0;L^2(\Omega))\cap L^2(-T,0;H^1(\Omega))$ and for some constant $C_0>0$, where
$\kappa_{p},\,\kappa_q\in[2,\infty)$ and $p,\,q\in[2,\infty)$ are such that
$$\displaystyle\frac{1}{\kappa_p}+\displaystyle\frac{1}{p}=\displaystyle\frac{1}{2}\,\,\,\,\,\textrm{and}
\,\,\,\,\,\displaystyle\frac{1}{\kappa_q}+\displaystyle\frac{1}{q}=\displaystyle\frac{1}{2}.$$
\indent\\
\begin{lemma}\label{Lem1}
Given $T>0$ fixed, let $u\in L^{\infty}(-T,0;L^2(\Omega))\cap L^2(-T,0;H^1(\Omega))$. Given $K_1,\,K_2\geq1$ integers and
$\xi\in\{1,\ldots,K_1\}$, and $\zeta\in\{1,\ldots,K_2\}$, let $r_{1,\xi},\,r_{2,\zeta},\,s_{1,\xi},\,s_{2,\zeta}\in[2,\infty)$ be such that
\begin{equation}\label{Lem1-01}
\displaystyle\frac{1}{r_{1,\xi}}+\displaystyle\frac{1}{s_{1,\xi}}=\displaystyle\frac{1}{2}\,\,\,\,\,\textrm{and}
\,\,\,\,\,\displaystyle\frac{1}{r_{2,\zeta}}+\displaystyle\frac{1}{s_{2,\zeta}}=\displaystyle\frac{1}{2}.
\end{equation}
Assume that there exists constants $k_0,\,\gamma_0\geq0$ and positive numbers $\theta_{1,\xi},\,\theta_{2,\zeta}$, such that
for all $\varsigma\in(0,T]$,\, $\varepsilon\in(0,1/2)$, and $k\geq k_0$, we have
\begin{equation}\label{Lem1-02}
\Upsilon^2_{(1-\varepsilon)\varsigma}(u_{k})\,\leq\,\displaystyle\frac{\gamma_0}{\varepsilon\varsigma}
\displaystyle\int^0_{-\varsigma}\displaystyle\int_{\Omega_k}|u_{k}(t)|^2\,dx\,dt+
\gamma_0k^2\displaystyle\sum^{K_1}_{\xi=1}\left(\displaystyle\int^0_{-\varsigma}|\Omega_k(t)|
^{^{\frac{r_{1,\xi}}{s_{1,\xi}}}}\,dt\right)^{^{\frac{2(1+\theta_{1,\xi})}{r_{1,\xi}}}}
+\gamma_0k^2\displaystyle\sum^{K_2}_{\zeta=1}\left(\displaystyle\int^0_{-\varsigma}\mu(\Gamma^{\infty}_k(t))
^{^{\frac{r_{2,\zeta}}{s_{2,\zeta}}}}\,dt\right)^{^{\frac{2(1+\theta_{2,\zeta})}{r_{2,\zeta}}}}.
\end{equation}
Then there exists a constant $C^{\ast}_0>0$ (independent of $u$ and $k_0$) such that
\begin{equation}\label{Lem1-03}
\textrm{ess}\displaystyle\sup_{[-T/2,0]\times\overline{\Omega}}\{u\}
\,\leq\,C^{\ast}_0\left(\|u\|^2_{_{L^{2}(-T,0;L^{2}(\Omega))}}+k^2_0\right)^{1/2}.
\end{equation}
Moreover, if instead of (\ref{Lem1-02}) one has
\begin{equation}\label{Lem1-04}
\Upsilon^2_{_T}(u_{k})\,\leq\,\gamma_0\displaystyle\int^0_{-T}\displaystyle\int_{\Omega_k}|u_{k}(t)|^2\,dx\,dt+
\gamma_0k^2\displaystyle\sum^{K_1}_{\xi=1}\left(\displaystyle\int^0_{-T}|\Omega_k(t)|
^{^{\frac{r_{1,\xi}}{s_{1,\xi}}}}\,dt\right)^{^{\frac{2(1+\theta_{1,\xi})}{r_{1,\xi}}}}
+\gamma_0k^2\displaystyle\sum^{K_2}_{\zeta=1}\left(\displaystyle\int^0_{-T}\mu(\Gamma^{\infty}_k(t))
^{^{\frac{r_{2,\zeta}}{s_{2,\zeta}}}}\,dt\right)^{^{\frac{2(1+\theta_{2,\zeta})}{r_{2,\zeta}}}},
\end{equation}
Then there exists a constant $C^{\star}>0$ (independent of $u$ and $k_0$) such that
\begin{equation}\label{Lem1-05}
\textrm{ess}\displaystyle\sup_{[-T,0]\times\overline{\Omega}}\{u\}
\,\leq\,C^{\ast}_0\left(\|u\|^2_{_{L^{2}(-T,0;L^{2}(\Omega))}}+k^2_0\right)^{1/2}.
\end{equation}
\end{lemma}

\begin{proof}
Since $|\Omega_k(t)|\leq|\Omega|$ and $\mu(\Gamma^{\infty}_k(t))\leq\mu(\Gamma^{\infty})$, it suffices to consider
$$\theta_{1,\xi}=\theta_{2,\zeta}=\theta:=\min\left\{\min\{\theta_{1,\xi}\mid1\leq\xi\leq K_1\}\cup\min\{\theta_{2,\zeta}\mid1\leq\zeta\leq K_2\}\right\}$$
for each $\xi\in\{1,\ldots,K_1\}$ and $\zeta\in\{1,\ldots,K_2\}$, provided that $\gamma_0$ is replaced by a larger constant $\gamma'_0$
(depending on $r_{1,\xi}$,\, $s_{1,\xi}$,\, $\theta_{1,\xi}$,\, $r_{2,\zeta}$,\, $s_{2,\zeta}$,\, $\theta_{2,\zeta}$,\, $T$,\, $\gamma_0$,\, $|\Omega|$, y
$\mu(\Gamma^{\infty})$). Given $\hat{k}>\hat{k}_0\geq0$ fixed real numbers,
define the sequences $\{\varsigma_n\}_{n\geq0}$,\, $\{k_n\}_{n\geq0}$,\, $\{y_n\}_{n\geq0}$,\, $\{z_n\}_{n\geq0}$ of positive real numbers, by:
\begin{equation}\label{Lem1-06}
\varsigma_n:=\displaystyle\frac{T}{2}(1+2^{-(n+1)}),\,\,\,\,\,\,k_n:=\hat{k}_0\chi_{_{\{n=0\}}}+(2-2^{-(n-1)})\hat{k}\chi_{_{\{n\geq1\}}},\,\,\,\,\,\,
y_n:=\displaystyle\frac{1}{\hat{k}^2}\displaystyle\int^0_{-\varsigma_n}\displaystyle\int_{\Omega_{k_n}}|u_{k_n}(t)|^2\,dx\,dt,
\end{equation}
and
\begin{equation}\label{Lem1-07}
z_n:=\displaystyle\sum^{K_1}_{\xi=1}\left(\displaystyle\int^0_{-\varsigma_n}|\Omega_{k_n}(t)|
^{^{\frac{r_{1,\xi}}{s_{1,\xi}}}}\,dt\right)^{^{\frac{2}{r_{1,\xi}}}}
+\displaystyle\sum^{K_2}_{\zeta=1}\left(\displaystyle\int^0_{-\varsigma_n}\mu(\Gamma^{\infty}_{k_n}(t))
^{^{\frac{r_{2,\zeta}}{s_{2,\zeta}}}}\,dt\right)^{^{\frac{2}{r_{2,\zeta}}}}.
\end{equation}
Clearly $\varsigma_n\in[T/2,T]$ and $k_n\geq\hat{k}_0$ for each nonnegative integer $n$. Moreover, it is clear that
$$|u_{k_n}(t)|\geq(k_{n+1}-k_n)^2\chi_{_{\overline{\Omega}_{k_{n+1}}(t)}},$$ so taking into account
the definitions of the sequences defined above together with H\"older's inequality and (\ref{Sobolev-parabolic}), we calculate and get that\\[2ex]
$\hat{k}^2y_{n+1}\,\leq\,c\left(\displaystyle\int^0_{-\varsigma_{n+1}}|\Omega_{k_{n+1}}(t)|\,dt\right)^{^{1-\frac{2}{\kappa_p}}}
\left(\displaystyle\int^0_{-\varsigma_{n+1}}\left\|u_{k_{n+1}}(t)
\right\|^{\kappa_p}_{_{p,\Omega_{k_{n+1}}}}\,dt\right)^{2/\kappa_p}\,
\leq\,c\left(\displaystyle\int^0_{-\varsigma_{n+1}}|\Omega_{k_{n+1}}(t)|\,dt\right)^{^{1-\frac{2}{\kappa_p}}}\Upsilon^2_{\varsigma_{n+1}}(u_{k_{n+1}})$\\
\begin{equation}\label{Lem1-08}
\indent\indent\leq\,
c\left\{(k_{n+1}-k_n)^{-2}\hat{k}^2y_n\right\}^{^{1-\frac{2}{\kappa_p}}}\Upsilon^2_{\varsigma_{n+1}}(u_{k_{n+1}})
\,\leq\,c\,4^{n+1}(y_n)^{^{1-\frac{2}{\kappa_p}}}\Upsilon^2_{\varsigma_{n+1}}(u_{k_{n+1}}),
\end{equation}
where $c>0$ is a constant that varies from line to line (this will be assumed for the remaining of the proof), and $p,\,\kappa_{p}\in(2,\infty)$ are such that
$$\displaystyle\frac{1}{\kappa_p}+\displaystyle\frac{1}{p}=\displaystyle\frac{1}{2},$$ and we are choosing $p>\kappa_p$. In the same way, we get that\\[2ex]
$4^{-(n+1)}\hat{k}^2z_{n+1}\,\leq\,(k_{n+1}-k_n)^2z_{n+1}$\\
\begin{equation}\label{Lem1-09}
\leq\,\displaystyle\sum^{K_1}_{\xi=1}\left(\displaystyle\int^0_{-\varsigma_{n+1}}\left\{\displaystyle\int_{\Omega_{k_n}}|u_{k_n}(t)|^{s_{1,\xi}}\,dx\right\}
^{^{\frac{r_{1,\xi}}{s_{1,\xi}}}}\,dt\right)^{^{\frac{2}{r_{1,\xi}}}}+
\displaystyle\sum^{K_2}_{\zeta=1}\left(\displaystyle\int^0_{-\varsigma_{n+1}}\left\{\displaystyle\int_{\Gamma^{\infty}_{k_n}}|u_{k_n}(t)|^{s_{2,\zeta}}\,d\mu\right\}
^{^{\frac{r_{2,\zeta}}{s_{2,\zeta}}}}\,dt\right)^{^{\frac{2}{r_{2,\zeta}}}}\,\leq\,c\,\Upsilon^2_{\varsigma_{n+1}}(u_{k_{n}}).
\end{equation}
Taking $\varsigma=\varsigma_n$ and $\varepsilon=1-\varsigma_{n+1}/\varsigma_n\geq2^{-(n+3)}$, we apply the assumption (\ref{Lem1-02})
to obtain that
\begin{equation}\label{Lem1-10}\Upsilon^2_{\varsigma_{n+1}}(u_{k_{n+1}})\,\leq\,
\frac{\gamma_0}{\varepsilon\varsigma_n}\hat{k}^2y_n+\gamma_0k^2_nz^{1+\theta}_n
\,\leq\,\gamma_02^{n+4}\hat{k}^2(T^{-1}+1)(y_n+z^{1+\theta}_n),
\end{equation}
and then combining (\ref{Lem1-08}), (\ref{Lem1-09}), and (\ref{Lem1-10}), we arrive at
\begin{equation}\label{Lem1-11}
y_{n+1}\,\leq\,c'\,8^n(y^{1+\delta}_n+z^{1+\theta}_ny^{\delta}_n)\,\,\textrm{and}\,\,
z_{n+1}\,\leq\,c'\,8^n(y_n+z^{1+\theta}_n),\,\,\textrm{for each integer}\,\,n\geq0,
\end{equation}
for some constant $c'>0$, where $\delta:=1-2/\kappa_p$. To complete the arguments, we seek estimations for $y_0$ and $z_0$. In fact, we easy see that
\begin{equation}\label{Lem1-12}
y_0=\displaystyle\frac{1}{\hat{k}^2}\displaystyle\int^0_{-\frac{3T}{4}}\displaystyle\int_{\Omega_{k_0}}|u_{k_0}(t)|^2\,dx\,dt\,\leq\,
\displaystyle\frac{1}{\hat{k}^2}\displaystyle\int^0_{-T}\displaystyle\int_{\Omega}|u(t)|^2\,dx\,dt.
\end{equation}
On the other hand, using (\ref{Lem1-02}) and following the approach as in (\ref{Lem1-08}) and (\ref{Lem1-09}), we have that\\[2ex]
$(\hat{k}-k_0)^2z_0\,\leq\,\displaystyle\sum^{K_1}_{\xi=1}\left(\displaystyle\int^0_{-\varsigma_{0}}\left\{\displaystyle\int_{\Omega_{k_0}}
|u_{k_0}(t)|^{s_{1,\xi}}\,dx\right\}^{^{\frac{r_{1,\xi}}{s_{1,\xi}}}}\,dt\right)^{^{\frac{2}{r_{1,\xi}}}}+
\displaystyle\sum^{K_2}_{\zeta=1}\left(\displaystyle\int^0_{-\varsigma_{0}}\left\{\displaystyle\int_{\Gamma^{\infty}_{k_0}}|u_{k_0}(t)|^{s_{2,\zeta}}\,d\mu\right\}
^{^{\frac{r_{2,\zeta}}{s_{2,\zeta}}}}\,dt\right)^{^{\frac{2}{r_{2,\zeta}}}}$\\
$$\leq\,c\,\Upsilon^2_{\varsigma_{0}}(u_{k_{0}})
\leq\,\displaystyle\frac{4c\gamma_0}{T}\displaystyle\int^0_{-T}\displaystyle\int_{\Omega}|u_{k_0}(t)|^2\,dx\,dt
+c\gamma_0k^2_0\displaystyle\sum^{K_1}_{\xi=1}\left[T|\Omega|^{^{\frac{r_{1,\xi}}{s_{1,\xi}}}}\right]^{^{\frac{2(1+\theta)}{r_{1,\xi}}}}
+\,c\gamma_0k^2_0\displaystyle\sum^{K_2}_{\zeta=1}\left[T\mu(\Gamma^{\infty})^{^{\frac{r_{2,\zeta}}{s_{2,\zeta}}}}\right]^{^{\frac{2(1+\theta)}{r_{2,\zeta}}}}.$$
Thus letting
$$c'':=c\gamma_0\max\left\{\frac{4}{T}\,,\,2\displaystyle\sum^{K_1}_{\xi=1}\left[T|\Omega|^{^{\frac{r_{1,\xi}}{s_{1,\xi}}}}\right]^{^{\frac{2(1+\theta)}{r_{1,\xi}}}}
+2\displaystyle\sum^{K_2}_{\zeta=1}\left[T\mu(\Gamma^{\infty})^{^{\frac{r_{2,\zeta}}{s_{2,\zeta}}}}\right]^{^{\frac{2(1+\theta)}{r_{2,\zeta}}}}\right\},$$
we arrive at
\begin{equation}\label{Lem1-13}
z_0\,\leq\,\displaystyle\frac{c''}{(\hat{k}-k_0)^2}\left(\displaystyle\int^0_{-T}\displaystyle\int_{\Omega_{k_0}}|u_{k_0}(t)|^2\,dx\,dt+k^2_0\right),
\,\,\,\,\,\,\,\,\,\,\,\,\textrm{for all}\,\,\,\hat{k}\geq k_0.
\end{equation}
Then, we define $$l:=\min\left\{\delta,\,\frac{\theta}{1+\theta}\right\}\,\,\,\,\,\textrm{and}\,\,\,\,\,
\omega:=\min\left\{\frac{1}{(2c')^{^{\frac{1}{\delta}}}8^{^{\frac{1}{\delta l}}}},\,\frac{1}{(2c')^{^{\frac{1+\theta}{\theta}}}\,8^{^{\frac{1}{\theta l}}}}\right\},$$
and then choose
\begin{equation}\label{Lem1-14}
\hat{k}:=\max\left\{\displaystyle\frac{2}{\sqrt{\eta}}\left(\displaystyle\int^0_{-T}\displaystyle\int_{\Omega}|u(t)|^2\,dx\,dt\right)^{1/2},\,
\displaystyle\frac{2\sqrt{c''}}{\eta^{^{\frac{1}{2(1+\theta)}}}}
\left(\displaystyle\int^0_{-T}\displaystyle\int_{\Omega}|u(t)|^2\,dx\,dt+k^2_0\right)^{1/2}+k_0\right\}.
\end{equation}
Then one sees that
\begin{equation}\label{Lem1-15}
\hat{k}\,\leq\,c^{\star}\left(\displaystyle\int^0_{-T}\displaystyle\int_{\Omega}|u(t)|^2\,dx\,dt+k^2_0\right)^{1/2}
\end{equation}
for some constant $c^{\star}>0$, and moreover the selection implies that
\begin{equation}\label{Lem1-16}
y_0\,\leq\,\omega\,\,\,\,\,\,\,\,\textrm{and}\,\,\,\,\,\,\,\,z_0\,\leq\,\omega^{^{\frac{1}{1+\theta}}}.
\end{equation}
In views of (\ref{Lem1-11}) and (\ref{Lem1-16}), we see that the sequences $\{y_n\}$,\, $\{z_n\}$ satisfy the conditions of Lemma \ref{lemma1}.
Thus, applying Lemma \ref{lemma1} for $\hat{k}$ given by (\ref{Lem1-14})
shows that $\displaystyle\lim_{n\rightarrow\infty}z_n=0$, and consequently
\begin{equation}\label{Lem1-17}
u(t)\,\leq\,\displaystyle\lim_{n\rightarrow\infty}k_n=2\hat{k}\,\,\,\,\textrm{a.e. in}\,\,\overline{\Omega},\,\,\,\,\,\textrm{for a.e.}\,\,
t\in\displaystyle\bigcap_{n\in\mathbb{N\!}}[-\varsigma_n,0]=[-T/2,0].
\end{equation}
Combining (\ref{Lem1-17}) with (\ref{Lem1-15}), we get lead to the fulfillment of (\ref{Lem1-03}). To complete the proof,
assume now that (\ref{Lem1-04}) is valid. Then by selecting $\varsigma_n:=T$ in the previous arguments, the proof runs in the exact way as in the previous case.
In particular, one has that (\ref{Lem1-17}) is valid for almost all $t\in\displaystyle\bigcap_{n\in\mathbb{N\!}}[-\varsigma_n,0]=[-T,0]$, leading to
the fulfillment of (\ref{Lem1-05}), as desired.
\end{proof}

\indent We now establish a central $L^{\infty}$-estimate which is global in space, but in general local in time. However, if the initial condition in problem (\ref{D01}) is zero, then the result can be shown to be also global in time.

\begin{theorem}\label{Thm-Partial-Bounded}
Given $T>0$ fixed, let $\kappa_p,\,\kappa_q,\,p,\,q\in[2,\infty)$ be such that
\begin{equation}\label{Thm1-00}
\displaystyle\frac{1}{\kappa_p}+\displaystyle\frac{1}{p}<1\,\,\,\,\,\,\,\,\textrm{and}
\,\,\,\,\,\,\,\,\displaystyle\frac{1}{\kappa_q}+\displaystyle\frac{1}{2q}<\frac{1}{2}.
\end{equation}
Given $u_0\in L^2(\Omega)$, let $f\in L^{\kappa_p}(0,T;L^p(\Omega))$ and $g\in L^{\kappa_q}(0,T;L^q_{\mu}(\Gamma^{\infty}))$, and
let $u\in C(0,T;L^2(\Omega))\cap L^2(0,T;\mathcal{V}_2(\Omega))$ be a weak solution of problem (\ref{D01}).
Then there exists a constant $C^{\ast}_1=C^{\ast}_1(T,N,\kappa_p,\kappa_q,p,q,|\Omega|,\mu(\Gamma^{\infty}))>0$ such that
\begin{equation}\label{Thm1-01}
|\|\mathbf{u}\||_{_{L^{\infty}(T/2,T;\mathbb{X\!}^{\,\infty}(\Omega;\Gamma))}}
\,\leq\,C^{\ast}_1\left(\|u\|_{_{L^{2}(0,T;L^{2}(\Omega))}}+\|f\|_{_{L^{\kappa_p}(0,T;L^{p}(\Omega))}}+
\|g\|_{_{L^{\kappa_q}(0,T;L^{q}_{\mu}(\Gamma^{\infty}))}}\right).
\end{equation}
If in addition $u_0=0$, then we get the global estimate
\begin{equation}\label{Thm1-02}
|\|\mathbf{u}\||_{_{L^{\infty}(0,T;\mathbb{X\!}^{\,\infty}(\Omega;\Gamma))}}
\,\leq\,C^{\ast}_1\left(\|u\|_{_{L^{2}(0,T;L^{2}(\Omega))}}+\|f\|_{_{L^{\kappa_p}(0,T;L^{p}(\Omega))}}+
\|g\|_{_{L^{\kappa_q}(0,T;L^{q}_{\mu}(\Gamma^{\infty}))}}\right).
\end{equation}
\end{theorem}

\begin{proof}
Let $u\in C(0,T;L^2(\Omega))\cap L^2(0,T;\mathcal{V}_2(\Omega))$ be a weak solution of problem (\ref{D01}). We divide this technical result into three major steps.\\
$\bullet$\,\, {\it \underline{Step 1}}. Assume in addition that $u\in C(0,T;L^2(\Omega))\cap L^2(0,T;\mathcal{V}_2(\Omega))$
is a classical solution of (\ref{D01}), and that $T\leq T_0$ for some sufficiently small constant $T_0>0$ (depending only in $\kappa_p$,
$\kappa_q$, $p$, $q$, $|\Omega|$, $\mu(\Gamma^{\infty})$, and the coefficients of the form $\mathcal{E}(\cdot,\cdot)$ given by (\ref{E-form}); this connstant will be specified later in the proof). By a linear
translation in time, we may assume that problem (\ref{D01}) is solved over $[-T,0]$, with initial value $u_0=u(-T)$.
Let $u_k(t)$ be the function given by (\ref{2.16a}), fix $\varsigma\in[-T,0]$, and select a function $\psi\in H^1(-\varsigma,0)$
such that $0\leq\psi(t)\leq1$ for all $t\in[-\varsigma,0]$. Assume in addition that either $\psi(-\varsigma)=0$, or $\varsigma=T$ and
$u_k(-T)=0$. Then for each $t\in[-\varsigma,0]$, observe that\\[2ex]
$\displaystyle\frac{\psi(t)^2}{2}\displaystyle\int_{\Omega_k(t)}|u_k(t)|^2\,dx=\displaystyle\int^t_{-\varsigma}
\displaystyle\frac{d}{d\xi}\left(\displaystyle\frac{\psi(\xi)^2}{2}\displaystyle\int_{\Omega_k(\xi)}|u_k(\xi)|^2\,dx\right)d\xi$\\
\begin{equation}\label{thm1-02}
=\displaystyle\int^t_{-\varsigma}\psi(\xi)\psi'(\xi)\left(\displaystyle\int_{\Omega_k(\xi)}|u_k(\xi)|^2\,dx\right)d\xi+
\displaystyle\int^t_{-\varsigma}\psi(\xi)^2\left(\displaystyle\int_{\Omega_k(\xi)}u_k(\xi)\frac{ du_k(\xi)}{d\xi}\,dx\right)d\xi.
\end{equation}
Using the fact that $u\in C(0,T;L^2(\Omega))\cap L^2(0,T;\mathcal{V}_2(\Omega))$
is a classical solution of (\ref{D01}), we see that
\begin{equation}\label{thm1-03}
\displaystyle\int_{\Omega_k(\xi)}u_k(\xi)\frac{ du_k(\xi)}{d\xi}\,dx=\displaystyle\int_{\Omega_k(\xi)}u_k(\xi)\frac{ du(\xi)}{d\xi}\,dx
=\displaystyle\int_{\Omega_k(\xi)}f(\xi,x)u_k(\xi)\,dx+\displaystyle\int_{\Gamma^{\infty}_k(\xi)}g(\xi,x)u_k(\xi)\,d\mu-\mathcal{E}(u(\xi),u_k(\xi)),
\end{equation}
Estimating the right hand side in (\ref{thm1-03}) and recalling (\ref{sharp-coercivity}), we clearly see that
$\mathcal{E}(u(\xi),u_k(\xi))\,\geq\,\mathcal{E}(u_k(\xi),u_k(\xi))$, and moreover\\[2ex]
\indent$\displaystyle\int_{\Omega_k(\xi)}f(\xi,x)u_k(\xi)\,dx+\displaystyle\int_{\Gamma^{\infty}_k(\xi)}g(\xi,x)u_k(\xi)\,d\mu$\\
\begin{equation}\label{thm1-09}
\leq\,\displaystyle\frac{1}{k}\displaystyle\int_{\Omega_k(\xi)}|f(\xi,x)|(|u_k(\xi)|^2+k^2)\,dx+
\displaystyle\frac{1}{k}\displaystyle\int_{\Gamma^{\infty}_k(\xi)}|g(\xi,x)|(|u_k(\xi)|^2+k^2)\,d\mu.
\end{equation}
From here, recalling the condition (\ref{E03}) and combining all the above into (\ref{thm1-03}), we deduce that
\begin{equation}\label{thm1-11}
\displaystyle\int_{\Omega_k(\xi)}u_k(\xi)\frac{ du_k(\xi)}{d\xi}\,dx\,\leq\,-c_0\|\nabla u_k(\xi)\|^2_{_{2,\Omega_k(\xi)}}
+\displaystyle\frac{1}{k}\displaystyle\int_{\Omega_k(\xi)}|f(\xi,x)|(|u_k(\xi)|^2+k^2)\,dx+\displaystyle\frac{1}{k}
\displaystyle\int_{\Gamma^{\infty}_k(\xi)}|g(\xi,x)|(|u_k(\xi)|^2+k^2)\,d\mu.
\end{equation}
From here, we insert (\ref{thm1-11}) into (\ref{thm1-02}) to get that\\[2ex]
$\min\left\{\displaystyle\frac{1}{2},c_0\right\}\Upsilon^2_{\varsigma}(\psi u_k)
\,\leq\,\displaystyle\sup_{-\varsigma\leq t\leq0}\left(\displaystyle\frac{\psi(t)^2}{2}\|u_k(t)\|^2_{_{2,\Omega_k(t)}}\right)+
c_0\displaystyle\int^0_{-\varsigma}\psi(t)^2\|\nabla u_k(t)\|^2_{_{2,\Omega_k(t)}}\,dt$\\
\begin{equation}\label{thm1-12}
\leq\,L_{\psi'}\displaystyle\int^0_{-\varsigma}\|u_k(t)\|^2_{_{2,\Omega_k(t)}}\,dt
+\displaystyle\frac{1}{k}\displaystyle\int^0_{-\varsigma}\displaystyle\int_{\Omega_k(t)}|f(t,x)|(\psi(t)^2|u_k(t)|^2+k^2)\,dx\,dt+\displaystyle\frac{1}{k}
\displaystyle\int^0_{-\varsigma}\displaystyle\int_{\Gamma^{\infty}_k(t)}|g(t,x)|(\psi(t)^2|u_k(t)|^2+k^2)\,d\mu\,dt,
\end{equation}
where $L_{\psi'}:=\displaystyle\sup_{-\varsigma\leq t\leq0}|\psi'(t)|$. To estimate the last two terms in (\ref{thm1-12}), let
$\theta_b:=(2\tilde{b}+2\tilde{\kappa}_b-\tilde{b}\tilde{\kappa}_b)(\tilde{b}\tilde{\kappa}_b)^{-1}>0$
for $b\in\{p,q\}$ and $\tilde{\xi}:=2\xi(\xi-1)^{-1}$ ($\xi\in(1,\infty)$), and notice that
$$\frac{1}{\tilde{\kappa}_p(1+\theta_p)}+\frac{1}{\tilde{p}(1+\theta_p)}=\frac{1}{\tilde{\kappa}_q(1+\theta_q)}+\frac{1}{\tilde{q}(1+\theta_q)}=\frac{1}{2}.$$
Letting
\begin{equation}\label{thm1-13}
k^2_0:=\|f\|^2_{_{L^{\kappa_p}(0,T;L^{p}(\Omega))}}+\|g\|^2_{_{L^{\kappa_q}(0,T;L^{q}_{\mu}(\Gamma^{\infty}))}},
\end{equation}
we apply H\"older's inequality to deduce that
\begin{equation}\label{thm1-14}
\displaystyle\int^0_{-\varsigma}\displaystyle\int_{\Omega_k(t)}\frac{1}{k}|f(t,x)|\psi(t)^2|u_k(t)|^2\,dx\,dt\,
\leq\,\displaystyle\frac{k_0}{k}\left(\displaystyle\int^0_{-\varsigma}|\Omega_k(t)|^{\tilde{\kappa}_p/\tilde{p}}dt\right)^{^{\frac{2\theta_p}{\tilde{\kappa}_p(1+\theta_p)}}}
\|\psi u_k\|^2_{_{L^{\tilde{\kappa}_p(1+\theta_p)}(-\varsigma,0;L^{\tilde{p}(1+\theta_p)}(\Omega))}}.
\end{equation}
Now observe that $$\displaystyle\lim_{\varsigma\rightarrow0^+}\left\{\left(\displaystyle\int^0_{-\varsigma}|\Omega_k(t)|^{\tilde{\kappa}_p/\tilde{p}}dt\right)
^{^{\frac{2\theta_p}{\tilde{\kappa}_p(1+\theta_p)}}}\|\psi u_k\|^2_{_{L^{\tilde{\kappa}_p(1+\theta_p)}(-\varsigma,0;L^{\tilde{p}(1+\theta_p)}(\Omega))}}\right\}=0,$$
so using this together with the selection of $\theta_p$, (\ref{thm1-14}), and (\ref{Sobolev-parabolic}), entail the existence of a constant
$T_1=T_1(\kappa_p,\theta_p,p,|\Omega|)>0$ such that
\begin{equation}\label{thm1-15}
\displaystyle\int^0_{-\varsigma}\displaystyle\int_{\Omega_k(t)}\frac{1}{k}|f(t,x)|\psi(t)^2|u_k(t)|^2\,dx\,dt
\,\leq\,\displaystyle\frac{k_0}{4k}\min\left\{\displaystyle\frac{1}{2},c_0\right\}\Upsilon^2_{\varsigma}(\psi u_k)
\end{equation}
whenever $\varsigma\leq T_1$. In the same way, recalling the selection of $\theta_q$, we deduce the existence
of positive constants $T_2=T_2(\kappa_q,\theta_q,q,\mu(\Gamma^{\infty}))$ such that
\begin{equation}\label{thm1-17}
\displaystyle\int^0_{-\varsigma}\displaystyle\int_{\Gamma^{\infty}_k(t)}\frac{1}{k}|g(t,x)|\psi(t)^2|u_k(t)|^2\,d\mu\,dt
\,\leq\,\displaystyle\frac{k_0}{4k}\min\left\{\displaystyle\frac{1}{2},c_0\right\}\Upsilon^2_{\varsigma}(\psi u_k),
\end{equation}
if $\varsigma\leq T_2$. Letting $T_0:=\min\{T_1,T_2\}$, and taking $\varsigma\in(0,T_0]$ and $k\geq k_0$, we insert
(\ref{thm1-15}) and (\ref{thm1-17}) into (\ref{thm1-12}), to arrive at
\begin{equation}\label{thm1-19}
\Upsilon^2_{\varsigma}(\psi u_k)\,\leq\,\displaystyle\frac{2}{\min\{1/2,c_0\}}\left\{
L_{\psi'}\displaystyle\int^0_{-\varsigma}\|u_k(t)\|^2_{_{2,\Omega_k(t)}}\,dt
+k^2\displaystyle\int^0_{-\varsigma}\displaystyle\int_{\Omega_k(t)}\frac{1}{k}|f(t,x)|\,dx\,dt+
k^2\displaystyle\int^0_{-\varsigma}\displaystyle\int_{\Gamma^{\infty}_k(t)}\frac{1}{k}|g(t,x)|\,d\mu\,dt\right\}.
\end{equation}
Furthermore,
\begin{equation}\label{thm1-20}
\displaystyle\int^0_{-\varsigma}\displaystyle\int_{\Omega_k(t)}\frac{1}{k}|f(t,x)|\,dx\,dt\,\leq\,\left(\displaystyle\int^0_{-\varsigma}|\Omega_k(t)|^{^{\frac{\kappa_p(p-1)}{p(\kappa_p-1)}}}dt\right)^{^{\frac{\kappa_p-1}{\kappa_p}}}\,\,\,\textrm{and}\,\,\,
\displaystyle\int^0_{-\varsigma}\displaystyle\int_{\Gamma^{\infty}_k(t)}\frac{1}{k}|g(t,x)|\,d\mu\,dt\,\leq\,
\left(\displaystyle\int^0_{-\varsigma}\mu(\Gamma^{\infty}_k(t))^{^{\frac{\kappa_q(q-1)}{q(\kappa_q-1)}}}dt\right)
^{^{\frac{\kappa_q-1}{\kappa_q}}}.
\end{equation}
Write $\gamma_0:=2(\min\{1/2,c_0\})^{-1}$ and choose
$\theta_{1}:=\theta_p$,\, $\theta_2:=\theta_q$,\,
$r_{1}:=(1+\theta_{1})\tilde{\kappa}_p$,\, $r_2:=(1+\theta_{2})\tilde{\kappa}_q$,\, $s_1:=(1+\theta_{1})\tilde{p},$ and $s_2:=(1+\theta_{2})\tilde{q}$. A straightforward calculation shows that the parameters $\theta_{i}$,\, $r_{i}$,\, $s_{i}$ ($i\in\{1,2\}$) fulfill (\ref{Lem1-01}). Then using all these notations, and substituting (\ref{thm1-20}) into (\ref{thm1-19}), we arrive at
\begin{equation}\label{thm1-24}
\Upsilon^2_{\varsigma}(\psi u)\,\leq\,\gamma_0L_{\psi'}\displaystyle\int^0_{-\varsigma}\|u_k(t)\|^2_{_{2,\Omega_k(t)}}\,dt
+\gamma_0k^2\left(\displaystyle\int^0_{-\varsigma}|\Omega_k(t)|
^{^{\frac{r_{1}}{s_{1}}}}\,dt\right)^{^{\frac{2(1+\theta_{1})}{r_{1}}}}
+\gamma_0k^2\left(\displaystyle\int^0_{-\varsigma}\mu(\Gamma^{\infty}_k(t))
^{^{\frac{r_{2}}{s_{2}}}}\,dt\right)^{^{\frac{2(1+\theta_{2})}{r_{2}}}}.
\end{equation}
Next, given $\varepsilon\in(0,1/2)$ fixed, select the function $\psi_{\varepsilon}\in H^1(-\varsigma,0)$ defined by: $$\psi_{\varepsilon}(t):=
\displaystyle\frac{t+\varsigma}{\varepsilon\varsigma}\chi_{_{\{-\varsigma\leq t<-(1-\varepsilon)\varsigma\}}}+\chi_{_{\{-(1-\varepsilon)\leq t\leq\varsigma,0\}}}.$$
Clearly $\psi_{\varepsilon}(-\varsigma)=0$,\,
$L_{\psi'_{\varepsilon}}:=\displaystyle\sup_{-\varsigma\leq t\leq0}|\psi'_{\varepsilon}(t)|\leq(\varepsilon\varsigma)^{-1}$, and $\Upsilon^2_{_{(1-\varepsilon)\varsigma}}(u)\,\leq\,\Upsilon^2_{\varsigma}(\psi_{\varepsilon}u)$. Thus, replacing $\psi=\psi_{\varepsilon}$ in (\ref{thm1-24}), we obtain that
\begin{equation}\label{thm1-25}
\Upsilon^2_{_{(1-\varepsilon)\varsigma}}(u)\,\leq\,\displaystyle\frac{\gamma_0}{\varepsilon\varsigma}\displaystyle\int^0_{-\varsigma}\|u_k(t)\|^2_{_{2,\Omega_k(t)}}\,dt
+\gamma_0k^2\left(\displaystyle\int^0_{-\varsigma}|\Omega_k(t)|
^{^{\frac{r_{1}}{s_{1}}}}\,dt\right)^{^{\frac{2(1+\theta_{1})}{r_{1}}}}
+\gamma_0k^2\left(\displaystyle\int^0_{-\varsigma}\mu(\Gamma^{\infty}_k(t))
^{^{\frac{r_{2}}{s_{2}}}}\,dt\right)^{^{\frac{2(1+\theta_{2})}{r_{2}}}}.
\end{equation}
Therefore, (\ref{thm1-25}) shows that (\ref{Lem1-02}) and all the conditions of Lemma \ref{Lem1} are fulfilled, and thus
applying Lemma \ref{Lem1} to both $u$ and $-u$ (the latter one being a weak solution of (\ref{D01}) for $f,\,g$ replaced by
$-f,\,-g$, respectively), and using the fact that $u|_{_{\Gamma\setminus\Gamma^{\infty}}}=0$, we establish (\ref{Thm1-01}) in the case $T\leq T_0$.\\
$\bullet$\,\, {\it \underline{Step 2}}. Assume now in addition that $u(-T)=u_0=0$. Then by choosing $\varsigma:=T$ and
and $\psi(t):=1$ for al $t\in[-T,0]$, equality (\ref{thm1-02}) becomes
\begin{equation}\label{thm1-26}
\displaystyle\frac{1}{2}\displaystyle\int_{\Omega_k(t)}|u_k(t)|^2\,dx
=\displaystyle\int^t_{-T}\left(\displaystyle\int_{\Omega_k(\xi)}u_k(\xi)\frac{ du_k(\xi)}{d\xi}\,dx\right)d\xi.
\end{equation}
From here, the calculations follow in the exact way as in the previous step, and in particular (\ref{thm1-12}) is transformed into\\[2ex]
$\min\left\{\displaystyle\frac{1}{2},c_0\right\}\Upsilon^2_{_T}(u)
\,\leq\,\displaystyle\sup_{-T\leq t\leq0}\left(\displaystyle\frac{1}{2}\|u_k(t)\|^2_{_{2,\Omega_k(t)}}\right)+
c_0\displaystyle\int^0_{-T}\|\nabla u_k(t)\|^2_{_{2,\Omega_k(t)}}\,dt$\\
\begin{equation}\label{thm1-27}
\leq\,\displaystyle\int^0_{-T}\left(\|u_k(t)\|^2_{_{2,\Omega_k(t)}}
+\frac{1}{k}\displaystyle\int_{\Omega_k(t)}|f(t,x)|(|u_k(t)|^2+k^2)\,dx
+\frac{1}{k}\displaystyle\int_{\Gamma^{\infty}_k(t)}|g(t,x)|(|u_k(t)|^2+k^2)\,d\mu\right)\,dt.
\end{equation}
Continuing analogously as in the previous case, we arrive a the inequality
\begin{equation}\label{thm1-28}
\Upsilon^2_{_T}(u)\,\leq\,\gamma_0\displaystyle\int^0_{-T}\|u_k(t)\|^2_{_{2,\Omega_k(t)}}\,dt
+\gamma_0k^2\left(\displaystyle\int^0_{-\varsigma}|\Omega_k(t)|
^{^{\frac{r_{1}}{s_{1}}}}\,dt\right)^{^{\frac{2(1+\theta_{1})}{r_{1}}}}
+\gamma_0k^2\left(\displaystyle\int^0_{-\varsigma}\mu(\Gamma^{\infty}_k(t))
^{^{\frac{r_{2}}{s_{2}}}}\,dt\right)^{^{\frac{2(1+\theta_{2})}{r_{2}}}},
\end{equation}
where all the constants and parameters appearing in (\ref{thm1-28}) are defined as in the previous step. Consequently,
(\ref{Lem1-04}) is valid for this case, and another application of Lemma \ref{Lem1} together with the fact that $u|_{_{\Gamma\setminus\Gamma^{\infty}}}=0$ gives (\ref{Thm1-02}).
This completes the proof of Theorem \ref{Thm-Partial-Bounded} in the case when the weak solution $u\in C(0,T;L^2(\Omega))\cap L^2(0,T;\mathcal{V}_2(\Omega))$ of (\ref{D01}) is also a classical solution of (\ref{D01}), and for $T\leq T_0$.\\
$\bullet$\,\, {\it \underline{Step 3}}. We now prove the general case. Since the proof runs in a similar way as in
\cite[proof of Proposition 3.1]{NITTKA2014}, we will only sketch the main ingredients of the proof. Let
$u\in C(0,T;L^2(\Omega))\cap L^2(0,T;\mathcal{V}_2(\Omega))$ be a weak solution of problem (\ref{D01}). Since $D(A_{\mu})$ is dense
in $L^2(\Omega)$, we pick $\{u_{0,n}\}\subseteq D(A^2_{\mu})$ such that $u_{0,n}\stackrel{n\rightarrow\infty}{\longrightarrow}u_0$ in $L^2(\Omega)$.
We also take sequences $\{f_n\}\subseteq C^2([0,T];L^{\infty}(\Omega))$ and $\{g_n\}\subseteq C^2([0,T];L^{\infty}_{\mu}(\Gamma^{\infty}))$
fulfilling $f_{n}\stackrel{n\rightarrow\infty}{\longrightarrow}f$ in $L^{\kappa_p}(0,T;L^p(\Omega))$ and
$g_{n}\stackrel{n\rightarrow\infty}{\longrightarrow}g$ in $L^{\kappa_q}(0,T;L^q_{\mu}(\Gamma^{\infty}))$, and $f_n(0,x)=g_n(0,x)=0$ for every $n\in\mathbb{N\!}\,$.
Using $f_n,\,g_n$ as the data functions in problem (\ref{D01}), from the arguments in subsection \ref{subsec6.2}, we get that problem (\ref{D01}) (with data $f_n,\,g_n$) admits a unique classical solution
$u_n$, and $u_{n}\stackrel{n\rightarrow\infty}{\longrightarrow}u$ in $C([0,T];L^2(\Omega))\cap L^2(0,T;\mathcal{V}_2(\Omega))$.
Then, given $T>0$ arbitrary, put $T'_0:=\min\{T_0,T\}>0$, for $T_0>0$ the parameter given in the previous two steps, and let $I'\subseteq[T'_0/2,T'_0]$ be an interval.
Then the inequality (\ref{Thm1-01}) is valid (over $I'$) for the function $u_n$ (with data $f_n,\,g_n$), and moreover an application of (\ref{Thm1-01}) to
the function $u_n-u_m$ ($n,\,m\in\mathbb{N\!}\,$) shows that $\{u_n\}$ is a Cauchy sequence in $L^{\infty}(I';L^{\infty}(\Omega))$. Thus
$u_{n}\stackrel{n\rightarrow\infty}{\longrightarrow}u$ over $L^{\infty}(I';L^{\infty}(\Omega))$. Therefore, passing to the limit and
using the fact that $[T/2,T]\subseteq\displaystyle\bigcup^m_{j=1}I'_i$ (for some $m\in\mathbb{N\!}\,$) with $\ell(I'_i)\leq L'_0/2$, we are lead into the fulfillment
of the inequality (\ref{Thm1-01}). If in addition $u_0=0$, then the same procedure is repeated over the inequality (\ref{Thm1-02}), and this
completes the proof of the theorem.
\end{proof}

\indent Now we provide a priori estimates depending in the initial condition of problem (\ref{D01}). We begin
with an estimate global in time. In this case, since the initial condition is only assumed to be bounded in $\Omega$, if $u_0\neq0$,
then it is still an open problem to obtain global $L^{\infty}$-estimates in the sense that in this case weak solutions to equation (\ref{D01}) are only known to be bounded
over $\Omega$ (in fact, there are globally bounded whenever $t>0$, but $u_0|_{_{\Gamma^{\infty}}}$ may fail to be bounded unless specified otherwise).
To be more precise, we have the following result whose proof runs as in \cite[proof of Theorem 3.2]{NITTKA2014}.

\begin{theorem}\label{infitiny-u_0-infinity}
Given $T>0$ arbitrary fixed, let $\kappa_p,\,\kappa_q,\,p,\,q\in[2,\infty)$ be as in Theorem \ref{Thm-Partial-Bounded}.
Given $u_0\in L^{\infty}(\Omega)$, let $f\in L^{\kappa_p}(0,T;L^p(\Omega))$ and $g\in L^{\kappa_q}(0,T;L^q_{\mu}(\Gamma^{\infty}))$, and
let $u\in C(0,T;L^2(\Omega))\cap L^2(0,T;\mathcal{V}_2(\Omega))$ be a weak solution of problem (\ref{D01}).
Then there exists a constant $C^{\ast}_2=C^{\ast}_2(T,\kappa_p,\kappa_q,p,q,|\Omega|,\mu(\Gamma^{\infty}))>0$ such that
\begin{equation}\label{Thm2-01}
\|u\|_{_{L^{\infty}(0,T;L^{\infty}(\Omega))}}
\,\leq\,C^{\ast}_2\left(\|u_0\|_{_{\infty,\Omega}}+\|f\|_{_{L^{\kappa_p}(0,T;L^{p}(\Omega))}}+
\|g\|_{_{L^{\kappa_q}(0,T;L^{q}_{\mu}(\Gamma^{\infty}))}}\right).
\end{equation}
\end{theorem}

\begin{proof}
If $u\in C(0,T;L^2(\Omega))\cap L^2(0,T;\mathcal{V}_2(\Omega))$ solves (\ref{D01}), then by linearity one can write $u(t):=T_{\mu}(t)u_0+\tilde{u}(t)$,
where $T_{\mu}(t)u_0$ denotes the weak solution of the homogeneous problem (\ref{D01}) for $f=g=0$, and $\tilde{u}(t)$ solves (\ref{D01}) for $u_0=0$. Then an application of (\ref{infinity-semigroup}) together with (\ref{Thm1-02}) entail that\\[2ex]
$\|u\|^2_{_{L^{\infty}(0,T;L^{\infty}(\Omega))}}\,\leq\,2\left(\displaystyle\sup_{t\in[0,T]}\|T_{\mu}(t)u_0\|^2_{_{\infty,\Omega}}
+\|\tilde{u}\|^2_{_{L^{\infty}(0,T;L^{\infty}(\Omega))}}\right)$\\
\begin{equation}\label{Thm2-02}
\leq\,2\left\{M^2e^{2|\omega|T}\|u_0\|^2_{_{\infty,\Omega}}+C^{\ast}_1\left(\|\tilde{u}\|_{_{L^{2}(0,T;L^{2}(\Omega))}}+
\|f\|_{_{L^{\kappa_p}(0,T;L^{p}(\Omega))}}+\|g\|_{_{L^{\kappa_q}(0,T;L^{q}_{\mu}(\Gamma^{\infty}))}}\right)\right\}.
\end{equation}
Furthermore, an application of Proposition \ref{apriori estimates} and H\"older's inequality gives that
\begin{equation}\label{Thm2-03}
\|\tilde{u}\|_{_{L^{2}(0,T;L^{2}(\Omega))}}
\leq\,cT\left(\|u_0\|^2_{_{\infty,\Omega}}+
\|f\|_{_{L^{\kappa_p}(0,T;L^{p}(\Omega))}}+\|g\|_{_{L^{\kappa_q}(0,T;L^{q}_{\mu}(\Gamma^{\infty}))}}\right),
\end{equation}
for some constant $c>0$.
Combining (\ref{Thm2-02}) and (\ref{Thm2-03}) yield the inequality (\ref{Thm2-01}), as asserted.
\end{proof}

\indent If problem (\ref{D01}) is local in time, we can deduce global boundedness in space.

\begin{theorem}\label{infitiny-u_0-L2}
Given $T_2>T_1>0$ arbitrary fixed, let $\kappa_p,\,\kappa_q,\,p,\,q\in[2,\infty)$ be as in Theorem \ref{infitiny-u_0-infinity}.
Given $u_0\in L^{2}(\Omega)$, let $f\in L^{\kappa_p}(0,T;L^p(\Omega))$ and $g\in L^{\kappa_q}(0,T;L^q_{\mu}(\Gamma^{\infty}))$, and
let $u\in C(0,T;L^2(\Omega))\cap L^2(0,T;\mathcal{V}_2(\Omega))$ be a weak solution to equation (\ref{D01}).
Then there exists a constant $C^{\ast}_3=C^{\ast}_3(T_1,T_2,\kappa_p,\kappa_q,p,q,|\Omega|,\mu(\Gamma^{\infty}))>0$ such that
\begin{equation}\label{Thm3-01}
|\|\mathbf{u}\||_{_{L^{\infty}(T_1,T_2;\mathbb{X\!}^{\,\infty}(\Omega;\Gamma))}}
\,\leq\,C^{\ast}_3\left(\|u_0\|_{_{2,\Omega}}+\|f\|_{_{L^{\kappa_p}(0,T;L^{p}(\Omega))}}+
\|g\|_{_{L^{\kappa_q}(0,T;L^{q}_{\mu}(\Gamma^{\infty}))}}\right).
\end{equation}
\end{theorem}

\begin{proof}
If $u\in C(0,T;L^2(\Omega))\cap L^2(0,T;\mathcal{V}_2(\Omega))$ is a weak solution of (\ref{D01}), then by using the fact that $[T_1,T_2]\subseteq\displaystyle\bigcup^m_{j=1}[T^{\ast}_j/2,T^{\ast}_j]$
with $T^{\ast}_1<T^{\ast}_2<\ldots<T^{\ast}_m$ (for some $m\in\mathbb{N\!}\,$), and $T^{\ast}_1=T_1$, \,$T^{\ast}_m=T_2$,
one can use the full strength of Theorem \ref{Thm-Partial-Bounded} to get that
\begin{equation}\label{Thm3-02}
|\|\mathbf{u}\||_{_{L^{\infty}(T_1,T_2;\mathbb{X\!}^{\,\infty}(\Omega;\Gamma))}}
\,\leq\,C^{\ast}_1\left(\|u\|_{_{L^{2}(0,T;L^{2}(\Omega))}}+\|f\|_{_{L^{\kappa_1}(0,T;L^{p}(\Omega))}}+
\|g\|_{_{L^{\kappa_2}(0,T;L^{q}_{\mu}(\Gamma))}}\right).
\end{equation}
Moreover, an application of Proposition \ref{apriori estimates} leads to the inequality (\ref{Thm3-01}), completing the proof.
\end{proof}

\subsection{Proof of Theorem \ref{Main-T2}}\label{subsec6.4}

The existence and uniqueness of a weak solution is given by Theorem \ref{exist-uniq-weak-soln}. To prove the global uniform continuity of the weak solution,
in views of Definition \ref{a-normal-a}, denote by $\mathbf{A}_{_{\mathbb{X\!}}}:=\mathbf{A}_{\mu}|_{_{\mathbb{X\!}^{\,p,q}(\Omega;\Gamma)}}$ the realization of the operator $\mathbf{A}_{\mu}$ in $\mathbb{X\!}^{\,p,q}(\Omega;\Gamma)$, with $D(\mathbf{A}_{_{\mathbb{X\!}}}):=\{(u,0)\in D(\mathbf{A}_{\mu})\mid \mathcal{L}u\in L^p(\Omega)\}$, and $\mathbf{A}_{_{\mathbb{X\!}}}(u,0)=\mathbf{A}_{\mu}(u,0)$. Then if $(f,g)\in\mathbb{X\!}^{\,p,q}(\Omega;\Gamma^{\infty})$, then $(u,0)\in D(\mathbf{A}_{_{\mathbb{X\!}}})$ with $\mathbf{A}_{_{\mathbb{X\!}}}(u,0)=(f,g)$ if and only if $u\in\mathcal{V}_2(\Omega)$ is a weak solution of the elliptic problem (\ref{E01}). By Theorem \ref{Holder-continuity}, we get that $u\in C(\overline{\Omega})$, and thus $D(\mathbf{A}_{_{\mathbb{X\!}}})\subseteq C(\overline{\Omega})\times\{0\}\subseteq\mathbb{X\!}^{\,p,q}(\Omega;\Gamma)$. Consequently, $\mathbf{A}_{_{\mathbb{X\!}}}$ is resolvent positive, and thus from \cite[Theorem 3.11.7]{AR-BA-HIE-NEU} we have that $\mathbf{A}_{_{\mathbb{X\!}}}$ generates an once integrated semigroup on $\mathbb{X\!}^{\,p,q}(\Omega;\Gamma)$. Now select sequences $\{f_n\}\subseteq C^2([0,T];L^{\infty}(\Omega))$ and $\{g_n\}\subseteq C^2([0,T];L^{\infty}_{\mu}(\Gamma^{\infty}))$ such that $f_n(0)=g_n(0)=0$, \,$f_n\stackrel{n\rightarrow\infty}{\longrightarrow}f$ in $L^{\kappa_p}(0,T;L^p(\Omega))$, and $g_n\stackrel{n\rightarrow\infty}{\longrightarrow}g$ in $L^{\kappa_q}(0,T;L^q_{\mu}(\Gamma^{\infty}))$, and let $v_n\in C(0,T;L^2(\Omega))\cap L^2(0,T;\mathcal{V}_2(\Omega))$ be the unique weak solution to problem (\ref{D01}) with data $(f_n,g_n)$ and initial condition $v(0)=0$. On the other hand, consider the abstract Cauchy problem
\begin{equation}\label{Abstract-Cauchy-Continuity}
\left\{
     \begin{array}{ll}
       \displaystyle\frac{\partial\mathfrak{w}_n}{\partial t}-\mathbf{A}_{_{\mathbb{X\!}}}\mathfrak{w}_n=(f_n(t),g_n(t))\,\,\,\,\,\,\,\,\,\,\,\,\,\,\,\,\textrm{for}\,\,t\in(0,\infty);\\
       \mathfrak{w}_n(0)=(0,0).
     \end{array}
   \right.
\end{equation}
By virtue of \cite[Corollary 3.2.11]{AR-BA-HIE-NEU}, it follows that problem (\ref{Abstract-Cauchy-Continuity}) has a unique solution $\mathfrak{w}_n:=(w_n,0)\in C^1([0,T];\mathbb{X\!}^{\,p,q}(\Omega;\Gamma))\cap C([0,T];D(\mathbf{A}_{_{\mathbb{X\!}}}))$, and since $D(\mathbf{A}_{_{\mathbb{X\!}}})\subseteq C(\overline{\Omega})\times\{0\}$, we see that $w_n\in C([0,T];C(\overline{\Omega}))$. Moreover, one gets that $w_n$ is a classical solution of problem (\ref{D01}) (with respect to the data $(f_n,g_n)$), and thus by uniqueness one gets $w_n=v_n$ for all $n\in\mathbb{N\!}\,$. Therefore $v_n(t)\in C([0,T];C(\overline{\Omega}))$, and by recalling Theorem \ref{infitiny-u_0-infinity}, we have that $v_n\stackrel{n\rightarrow\infty}{\longrightarrow}v$ uniformly over $[0,T]\times\overline{\Omega}$, and thus $v(t)\in C([0,T];C(\overline{\Omega}))$ is a weak solution of equation (\ref{D01}) with $v(0)=0$. Additionally, if $u_0\in C(\overline{\Omega})$, then in views of Theorem \ref{feller-semi} we have that $T_{\mu}(t)u_0\in C([0,T];C(\overline{\Omega}))$. Therefore, $u(t)=T_{\mu}(t)u_0+v(t)\in C([0,T];C(\overline{\Omega}))$. Finally, inequality (\ref{Main-Eq2}) follows at once from the above conclusion together with the estimate (\ref{Thm2-01}) in Theorem \ref{infitiny-u_0-infinity}, and we have completed the proof of the result.\,\,\,\,\,\,$\Box$

\begin{remark}
In the case of positive time, from Theorem \ref{feller-res} we see that $\left(T_{\mu}(t)\right)(L^p(\Omega))\subseteq C(\overline{\Omega})$ for all $p\in(1,\infty)$ and for all $t>0$. This together with Theorem \ref{infitiny-u_0-L2} and the arguments above imply that if $u_0\in L^2(\Omega)$, then for each $0<T_1<T_2$ fixed, one has that $u(t)=T_{\mu}(t)u_0+v(t)\in C([T_1,T_2];C(\overline{\Omega}))$, with
\begin{equation}\label{Remark-L2-Continuity}
\|u\|_{_{C([T_1,T_2];C(\overline{\Omega}))}}
\,\leq\,C^{\ast}_3\left(\|u_0\|_{_{2,\Omega}}+\|f\|_{_{L^{\kappa_p}(0,T;L^{p}(\Omega))}}+
\|g\|_{_{L^{\kappa_q}(0,T;L^{q}_{\mu}(\Gamma^{\infty}))}}\right).
\end{equation}
\end{remark}

\begin{remark}
The results obtained in this paper can be extended to diffusion equations involving the non-symmetric operator
$$\mathcal{A}_0u:=\mathcal{A}-\displaystyle\sum^2_{i=1}\partial_{x_i}(\theta_iu),$$
for $\theta_i\in W^{1,\infty}_0(\Omega)$, provided that the coefficient $\lambda$ is ``large enough." In the case of Lipschitz domains, Nittka \cite{NITTKA2014,NITTKA} obtained global regularity results under lower assumptions on the coefficients $\theta_i$, but the methods employed in \cite{NITTKA2014,NITTKA} are only valid for Lipschitz domains, and cannot be replicated nor extended to the case of non-Lipschitz domains. We are optimistic that the results can be obtained under less regularity in the coefficients $\theta_i$ ($i\in\{1,2\}$, but so far there is no proof obtained achieving this.\\
\end{remark}

\begin{center}
    \textsc{A. Appendix}
\end{center}

\indent In this part we develop a mechanism to approximate the Hausdorff measure of the ramified boundary $\Gamma^{\infty}$, which gives us a correct notion of the length of the ramified generation of bronchial trees as the ones considered in this paper.\\

\noindent{\it A1. Language and tools}
\indent\\

\indent We briefly introduce the language and tools needed in order to establish the first main result of the paper. We begin with the following definition.

\begin{adefinicion}
    
    Let $K$ be the unique non-empty compact self-similar invariant set under an iterated function system (IFS) $\mathfrak{G}=\{G_{j}\}_{j=1}^{M}$ satisfying the open set condition (OSC), where $G_{j}$ has ratio $0<r_{j}$. Let $\mathbf{M}:=\{1, 2, \ldots, M\}$ and $n\geq 1$. We define the word space associated to $K$ as $\Theta:=\mathbf{M}^{\mathbb{N}}$ and $\Theta_{n}:=\mathbf{M}^{n}$ with the $n$-truncation map $[\cdot]_{n}: \Theta \to \Theta_{n}$ define for a word $\omega=\omega_1 \omega_2 \cdots \in \Theta$ by $[\omega]_{n}:=\omega_{1}\cdots \omega_{n}$
\label{A1}
    \end{adefinicion}

\indent We will not deal with the trivial case when $M=1$. Observe that there is a relation between the word space $\Theta$ and the attractor $K$ of an IFS with $M$ maps, where we identify points in $K$ with infinite words and regions with finite words. Then, for $\omega \in \Theta_{n}$, we define $K_{(\omega)}:=G_{\omega}(K)$ where $G_{\omega_1 \omega_2 \cdots \omega_n}$ is given inductively as $F_{\omega_2 \cdots \omega_n}\circ F_{\omega_1}$. Moreover for $\omega \in \Theta$, we define the point $K_{(\omega)}$ as the unique point in $\displaystyle\bigcap_{n\in\mathbb{N}}K_{[\omega]_{n}}$. Given $\mu_p$ the natural probability measure on $K$, then for $\omega=\omega_1 \cdots \omega_n \in \Theta_n$ we get that $\mu(K_{(\omega)})=\displaystyle\prod_{j=1}^{n}r_{\omega_{j}}^{d}=\displaystyle\sum_{j=0}^{M}r_{j}^{d}\mu(K_{(\omega)})$ (since $\mathfrak{G}$ fulfills OSC).\\

\begin{adefinicion}
     Let $K(n):=\{K_{(\omega)}| \omega\in \Theta_{n}\}$ be the set of $n$-cells of $K$, where we reserve the notation $\Delta_{j}^{(n)}$ for elements of $K(n)$, which we call $n$-cells. We also define $K(0):=K$.

     \label{A2}
\end{adefinicion}

\indent We now collect several facts given by Jia \cite{JIA07-1, JIA07-2} into one main result which will be very useful in the proof of the first main result.\\

\begin{ateorema}
     (See \cite{JIA07-1, JIA07-2})\,\,
Suppose that $K$ is a self-similar set that satisfies the open set condition. For $n\geq 1$, $1\leq k\leq M^{n}$, let $\Delta_{1}^{(n)}, \Delta_{2}^{(n)}, \ldots, \Delta_{k}^{(n)} \in K(n)$ and $\mu$ be the common self-similar probability measure on the $K$, defined by $\mu(K(n))=\prod_{j=1}^{n}r_{j}^{d}$. Let
$$q_{k}:=\min_{\{\Delta_{j}^{(n)}\}_{j=1}^{k} \subseteq K(n)}\left\{\frac{\left|\bigcup_{i =1}^{k}\Delta_{j}\right|^{d}}{\mu\left(\bigcup_{j=1}^{k}\Delta_{j}\right)}\right\}, \ j=1, \ldots, k$$
where the minimum is taken for every possible union of $k$ elements of $K(n)$ and $a_{n}=\min_{1\leq k\leq M^{n}}\{q_{k}\}$. If there exists a constant $\alpha> 0$ such that $\alpha \leq a_{n}$ for all $n\in\mathbb{N}$, then $\mathcal{H}^{s}(K)\geq \alpha$. Furthermore, $\{a_{n}\}$ is a decreasing sequence with $\displaystyle\lim_{n\to\infty}a_{n}=\mathcal{H}^{d}(K)$. Moreover, if $r_{j}=r$ for $j\in\{1,\ldots, m\}$, then $q_{k}=\min_{\{\Delta_{j}^{(n)}\}_{j=1}^{k} \subseteq K(n)}\left\{\frac{\left|\bigcup_{i =1}^{k}\Delta_{i}\right|^{d}}{kr^{nd}}\right\}$, where the minimum is taken for all possible union of $k$-elements of $K(n)$ and $a_{n}=\min_{1\leq k\leq M^{n}}\{q_{k}\}$. Also in the latter case, $\{a_{n}\}_{_{n\in\mathbb{N\!}}}$ is a decreasing sequence and $\displaystyle\lim_{n\to\infty}a_{n}=\mathcal{H}^{d}(K)$.
\label{TA1}
\end{ateorema} 

\indent We think of $\left | \bigcup_{j=1}^{k}\Delta_{j}^{(n)} \right |^{d}$ as an analog of the $d$-dimensional volume, and since $\mu$ is the natural probability measure in $K$, we can think of $q_{k}$ as the search for the collection of $k$ $n$-cells with the highest \textit{density}. By means of this analogy, one can think of $a_n$ as finding the collection of $n$-cells with the highest \textit{density} for each $n$.  In views of Theorem A.1, the main goal of this section consists in constructing an increasing sequence converging towards $\mathcal{H}^{d}(\Gamma^\infty)$. This will be achieved using case-by-case analysis. To proceed, we add some additional definitions and notations, as in \cite{FERR-VELEZ18-1}.

\begin{adefinicion}
    We say a proposition $P$ on the subsets of $K$ is a {\bf (valid) case}, if (\textbf{I}) for every $n$, there is a family $\{\Delta_{j}^{(n)}\}_{j=1}^{k} \subseteq K(n)$ such that $\bigcup_{j=1}^{k}\Delta_{j}^{(n)}$ satisfies $P$ and (\textbf{II}) for $n>1$ and $\{\Delta_{j}^{(n)}\}_{j=1}^{k} \subseteq K(n)$ such that $\bigcup_{j=1}^{k}\Delta_{j}^{(n)}$ satisfies case $P$, there is a family $\{\Delta_{i}^{(n-1)}\}_{i} \subseteq K(n-1)$ such that $\bigcup_{j=1}^{k}\Delta_{j}^{(n)} \subseteq \bigcup_{j}\Delta_{j}^{(n-1)}$ and $\bigcup_{j}\Delta_{j}^{(n-1)}$ satisfies case $P$. 

    \label{A3}
\end{adefinicion}

\indent   We close this subsection by considering the notion of height measure of $\Omega$ that will depend on the parameter $\tau$ with $\tau \leq \tau^{\ast}$, as in \cite{ACH10}.

\begin{adefinicion}
    We define
    \begin{apequation}
        \label{H}
    \mathfrak{H}:=\sup_{(x_1, x_2)\in \Omega}x_{2}=\frac{1+\frac{3\tau}{\sqrt{2}}}{1-\tau^2}.
\end{apequation}
\label{A4}
\end{adefinicion}

\noindent{\it A2. Proof of Theorem \ref{Main-T1}}
\indent\\

\indent To begin, we observe that the measure of the base hexagon of our domain is approximately $\text{diam}(\overline{\Omega^{(0)}})=\text{diam}(\overline{V})=1+\sqrt{2}\tau$; see table (\ref{tabla1}) below.

\begin{table}[h!]
\begin{center}
\begin{tabular}{|c|c|}
\hline
\textbf{Measurements of the diameters} & \textbf{\begin{tabular}[c]{@{}c@{}}Number of Hexagons\\ {[}Newly generated{]}\end{tabular}} \\ \hline \hline
           $\text{diam}(\overline{\Omega^{(1)}}\setminus \Omega^{(0)})=\frac{1}{2}+\frac{\tau}{\sqrt{2}}$                   &    2                         \\ \hline
                  
    $\text{diam}(\overline{\Omega^{(2)}}\setminus \Omega^{(1)})=\frac{1}{4}+\frac{\tau}{2\sqrt{2}}$         &       4                                                                     \\ \hline
  $\text{diam}(\overline{\Omega^{(3)}}\setminus \Omega^{(2)})=\frac{1}{8}+\frac{\tau}{4\sqrt{2}}$                              &        8                                      \\ \hline
        $\text{diam}(\overline{\Omega^{(4)}}\setminus \Omega^{(3)})=\frac{1}{16}+\frac{\tau}{8\sqrt{2}}$                              &        16                       \\ \hline
                $\vdots$              &   $\vdots$                                           \\ \hline
            $\text{diam}(\overline{\Omega^{(k)}}\setminus \Omega^{(k-1)})=\frac{1}{2^{k}}+\frac{\tau}{2^{k-1}\sqrt{2}}$                              &        $2^{k}$                                              \\ \hline
\end{tabular}
\caption{Measurements of the diameters of  $2^{k}$-hexagons.}
\label{tabla1}
\end{center}

\end{table}

Note that $\text{diam}(\overline{\Omega^{(k)}}\setminus \Omega^{(k-1)})\leq \text{diam}(\overline{\Omega^{(k-1)}}\setminus \Omega^{(k-2)})$ for all $k\in \mathbb{Z}^{+}$. We can express $\overline{\Omega}$ as a disjoint union, as follows: $\overline{\Omega}:=\bigcup_{k=1}^{\infty}(\overline{\Omega^{(k)}}\setminus \Omega^{(k-1)})$. In views of table \ref{tabla1} we get $\text{diam}(\overline{\Omega})=\sum_{k=1}^{\infty}\text{diam}(\overline{\Omega^{(k)}}\setminus \Omega^{(k-1)})=\sum_{k=1}^{\infty}\left(\frac{1}{2^{k}}+\frac{\tau}{2^{k-1}\sqrt{2}}\right)=1+\frac{2\tau}{\sqrt{2}}$. Simple geometric constructions were used to calculate the diagonal of the first hexagon, assuming that the base length is 2. Initially, we construct a rectangle with a base of length 2 and position it on a plane centered at the origin, so that the base aligns with the interval $[-1, 1]$. The height of the rectangle can vary. Then, we take transversal cuts at the upper corners with a common angle, for example, 45 degrees, and these cuts have a size of $2\tau$. Next, we make two copies of this hexagon, each with a base of length $2\tau$, and place them at the corners. This process repeats indefinitely, so that for the final terminations of each $k$-prefractal, there are $2^{k}$ hexagons. \\
\indent We now wish to find a way to estimate level below similar to the one used in \cite{JIA07-1, JIA07-2}. First, to apply this method there must be no overlap, and thus we take the ramified domain $\Omega$ with parameter $\tau\leq \tau^{*}$. Also, it may be simpler to use the fractal construction that depends on a previous iteration and not on the initial hexagonal domain. It is worth noting that due to the \textit{mass distribution principle} (Cf. \cite[Theorem 4.2]{FAL90}), one immediately sees that $\mathcal{H}^{d}(\Gamma^{\infty})$ is strictly positive.\\

\noindent{\bf Remark A1.}\, {\it
$\Gamma_{n}^{\infty}:=G_{i_1}\circ \cdots G_{i_n}(\Gamma^{\infty})$ with $1\leq i_{1}, \ldots, i_{n}\leq 2$}.\\

Let $\{\Delta_{j}^{(n)}\}_{j}$ be a collection of $n$-cells and let $s=\textrm{diam}\left(\bigcup_{j=1}^{k}\Delta_{j}^{(n)}\right)$ the diameter of such a collection. We organize the proof by cases based on how many lower 1-cells $\Gamma^\infty_1, \Gamma^\infty_2$ the family $\bigcup_{j=1}^{k}\Delta_{i}^{(n)}$ intersects. 
Recall that the Hausdorff dimension of $\Gamma^{\infty}$ is $d=-\frac{\log 2}{\log \tau}$. Now, Let $n\geq 1$, $1\leq k\leq 2^{n}$, $\Delta_{1}, \Delta_{2}, \ldots, \Delta_{k} \in \Gamma_{n}^{\infty}$. Note that $[\Gamma_1^\infty]$ and $[\Gamma_2^\infty]$ can be regarded as equivalence classes of the elements of $\{\Gamma_k^\infty\}_{k=1}^{n}$, which $(i_{l}, i_{m})=(1, 2)$ and $(i_{m}, i_{l})=(2, 1)$ with $i_k< i_j$. Strategically, it is sufficient to work with these equivalence classes, and we can exclude the other scenarios (see Definition \ref{A3}), and therefore one can exclude the other cases, and simply work with these classes.

\thinspace

\begin{figure}[h!]
    \centering
  \begin{tikzpicture}[thick]
\draw[fill=gray!15] (0,1.1) -- (0,0) -- (3,0) -- (3,1.1) -- (2.15,2) -- (0.85,2) -- cycle; 
\draw[shift={(0,1.1)}, rotate=46.2,scale=0.415,semithick,fill=gray!15] (0,1.1) -- (0,0) -- (3,0) -- (3,1.1) -- (2.15,2) -- (0.85,2) -- cycle;
\draw[shift={(2.15,2)}, rotate=-46.2,scale=0.415,semithick,fill=gray!15] (0,1.1) -- (0,0) -- (3,0) -- (3,1.1) -- (2.15,2) -- (0.85,2) -- cycle;
\draw[thin] (2.15,1.9) .. controls (2.1,1.7) and (2.4,1.55) .. (2.3,1.35);
\draw[thin] (2.3,1.35) .. controls (2.5,1.45) and (2.85,0.97) .. (2.9,1.1);
\node[rotate=-32] at (2.2,1.2) {\scriptsize $2\tau$};
\draw[shift={(3.37,1.93)}, rotate=-93,scale=0.172,thin,fill=gray!15] (0,1.1) -- (0,0) -- (3,0) -- (3,1.1) -- (2.15,2) -- (0.85,2) -- cycle;
\draw[shift={(2.476,2.3165)}, rotate=0.5,scale=0.1715,thin,fill=gray!15] (0,1.1) -- (0,0) -- (3,0) -- (3,1.1) -- (2.15,2) -- (0.85,2) -- cycle;
\draw[shift={(0.02,2.318)}, rotate=-0.65, scale=0.1725,thin,fill=gray!15] (0,1.1) -- (0,0) -- (3,0) -- (3,1.1) -- (2.15,2) -- (0.85,2) -- cycle;
\draw[shift={(-0.3285,1.41)}, rotate=93,scale=0.172,thin,fill=gray!15] (0,1.1) -- (0,0) -- (3,0) -- (3,1.1) -- (2.15,2) -- (0.85,2) -- cycle;
\draw[shift={(3.56195,1.919)}, rotate=-46.2,scale=0.06915,very thin,gray,fill=gray!15] (0,1.1) -- (0,0) -- (3,0) -- (3,1.1) -- (2.15,2) -- (0.85,2) -- cycle;
\draw[shift={(3.693,1.54)}, rotate=-139.8,scale=0.06915,very thin,gray,fill=gray!15] (0,1.1) -- (0,0) -- (3,0) -- (3,1.1) -- (2.15,2) -- (0.85,2) -- cycle;

\draw[shift={(2.843,2.662)}, rotate=-46,scale=0.06918,very thin,gray,fill=gray!15] (0,1.1) -- (0,0) -- (3,0) -- (3,1.1) -- (2.15,2) -- (0.85,2) -- cycle;
\draw[shift={(2.4776,2.508)}, rotate=46.98,scale=0.06921,very thin,gray,fill=gray!15] (0,1.1) -- (0,0) -- (3,0) -- (3,1.1) -- (2.15,2) -- (0.85,2) -- cycle;

\draw[shift={(0.024,2.509)}, rotate=46.7,scale=0.06925,very thin,gray,fill=gray!15] (0,1.1) -- (0,0) -- (3,0) -- (3,1.1) -- (2.15,2) -- (0.85,2) -- cycle;
\draw[shift={(0.3954,2.659)}, rotate=-47.2,scale=0.06927,very thin,gray,fill=gray!15] (0,1.1) -- (0,0) -- (3,0) -- (3,1.1) -- (2.15,2) -- (0.85,2) -- cycle;

\draw[shift={(-0.688,1.765)}, rotate=46.6,scale=0.06927,very thin,gray,fill=gray!15] (0,1.1) -- (0,0) -- (3,0) -- (3,1.1) -- (2.15,2) -- (0.85,2) -- cycle;
\draw[shift={(-0.5188,1.403)}, rotate=140,scale=0.06927,very thin,gray,fill=gray!15] (0,1.1) -- (0,0) -- (3,0) -- (3,1.1) -- (2.15,2) -- (0.85,2) -- cycle;
\end{tikzpicture}
\caption{The pre-fractal $\Omega^{(3)}$.}
\label{omega2}
\end{figure}
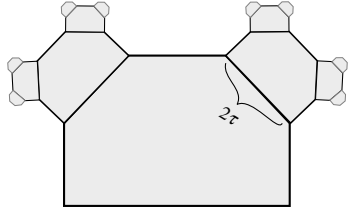

\noindent$\bullet$\,\,\underline{\textit{Case 1.}} Assume that $\bigcup_{i=1}^{k}\Delta_{i}^{(n)}$ intersects one of $G_{i}(\Gamma^{\infty})$,  for $i\in\{1,2\}$. Without loss of generality, we will simply work with the case $i=1$, since the case $i=2$ is similar since it is a symmetry (see Figure \ref{omega2}).

\noindent As mentioned before, we will proceed in a similar way to \cite{ JIA07-1, JIA07-2}. Given that $\bigcup_{i=1}^{k}{\Delta_{i}^{(n)}}$ intersects $G_1(\Gamma^\infty)$, we can assume that the diameters of each of them are smaller than the pre-fractals in which they are contained. It now holds that every diameter of the family $\{\Delta_i\}_{i}$ is minimized by the diameter of one of the $2^N$-hexagons, that is,\\[2ex] 
\begin{align*}
\textrm{diam}(\Delta_i)  &\leq \textrm{diam}(\overline{\Omega^{(0)}})+\displaystyle\sum_{m=1}^{\infty}\textrm{diam}(\overline{\Omega^{(m)}}\setminus \Omega^{(m-1)})\\ 
 &= 1+\tau\sqrt{2}+\displaystyle\sum_{m=1}^{\infty}\left(\frac{1}{2^{m}}+\frac{\tau}{2^{m-1}\sqrt{2}}\right)=(1+\tau\sqrt{2})+1+\frac{2\tau}{\sqrt{2}}=2+2\sqrt{2}\tau
\end{align*}whenever $k<N$ with $i=1, \ldots, k$, and therefore $\textrm{diam}\left(\bigcup_{i=1}^{k}\Delta_i\right)\leq 2+2\sqrt{2}\tau$ (see Table \ref{tabla1}). On the other hand, the diameter of each $\Delta_i$ will exceed the diameter of one of the $2^N$ hexagons (see table \ref{tabla1}), i.e., $\textrm{diam}(\Delta_i)\geq \frac{1}{2^{N}}+\frac{\tau}{2^{N-1}\sqrt{2}}$ for each $i=1, \ldots, k$. All the more reason, $\textrm{diam}\left(\bigcup_{i=1}^{k}\Delta_{i} \right)\geq \frac{1}{2^{N}}+\frac{\tau}{2^{N-1}\sqrt{2}}$, with fixed $N$. Thus, 
\begin{apequation}
    \frac{1}{2^{N}}+\frac{\tau}{2^{N-1}\sqrt{2}}\leq \textrm{diam}\left(\bigcup_{i=1}^{k}\Delta_{i} \right) \leq 2+2\sqrt{2}\tau
    \label{deltas}
\end{apequation}

Let $s=\textrm{diam}\left(\bigcup_{i=1}^{k}\Delta_{i} \right)$. By virtue of (\ref{deltas}) it follows that $s\geq \frac{1}{2^{N}}+\frac{\tau}{2^{N-1}\sqrt{2}}$, with $\tau\leq \tau^{*}$ and $N\in \mathbb{N}$. For each $\Delta_{i}$, there exists $\Delta_{j}^{n+1}\in \Gamma_{n-1}^{\infty}$ s.t. $\Delta_{i}\subset \Delta_{j}^{n-1}$ $(j=1, 2, \ldots, k_{n-1})$ and $\Delta_{1}^{n-1}, \Delta_{2}^{n-1}, \ldots, \Delta_{k_{n-1}}^{n-1}$ are different from each other. It can also be seen that $\tau k\leq k_{n-1}$ and
$$\textrm{diam}\left(\bigcup_{j=1}^{k_{n-1}}\Delta_{j}^{n-1} \right)\leq s+2\tau^{n},$$ also
$$\frac{\textrm{diam}\left( \bigcup_{i=1}^{k}\Delta_{i} \right)}{\textrm{diam}\left ( \bigcup_{j=1}^{k_{n-1}}\Delta_{j}^{n-1} \right )}\geq \frac{s}{s+2\tau^{n}}\geq \frac{\frac{1}{2^{N}}+\frac{\tau}{2^{N-1}\sqrt{2}}}{\frac{1}{2^{N}}+\frac{\tau}{2^{N-1}\sqrt{2}}+2\tau^n}=\frac{\frac{1+\sqrt{2}\tau}{2^{N}}}{\frac{1+\sqrt{2}\tau}{2^{N}}+2\tau^{n}}=\frac{1}{1+\tau^{n}\left(2^{1-N}(1+\sqrt{2}\tau)\right)}$$ Therefore, 
\begin{align*}
\frac{\left(\textrm{diam}\left ( \bigcup_{i=1}^{k}\Delta_{i} \right )\right)^{d}}{k\tau^{nd}} &\geq \left (\frac{1}{1+\tau^{n}\left(2^{1-N}(1+\sqrt{2}\tau)\right)}\right)^{d}\frac{\left(\textrm{diam}\left(\bigcup_{j=1}^{k_{n-1}}\Delta_{j}^{n-1} \right )\right)^{d}}{(\tau k)\tau^{d(n-1)}}\\ 
 &\geq \left (\frac{1}{1+\tau^{n}\left(2^{1-N}(1+\sqrt{2}\tau)\right)} \right )^{d}\frac{\left(\textrm{diam}\left( \bigcup_{j=1}^{k_{n-1}}\Delta_{j}^{n-1} \right )\right)^{d}}{k_{n-1}\tau ^{d(n-1)}}
\end{align*}
Having in mind the construction and properties of $\{a_{n}\}$ in Theorem \ref{TA1} (second part), we put
$$a_{n}^{(1)}=\min_{1\leq k\leq 2^{n}}\min_{\Delta_{i} \in \Gamma_{n}^{\infty}}\min_{i=1, \ldots, k}\left\{\frac{\textrm{diam}\left(\bigcup_{i =1}^{k}\Delta_{i}\right)^{d}}{k\tau^{nd}}, \, \bigcup_{i=1}^{k}\Delta_{i}\,\,\,\,\textrm{satisfies Case 1}\,\right\}$$
It is obvious that $a_{n}\leq a_{n}^{(1)}$. From the above inequality:
$$a_{n}^{(1)}\geq \left ( \frac{1}{1+\tau^{n}\left(2^{1-N}(1+\sqrt{2}\tau)\right)} \right )^{d} a_{n-1}^{(1)}$$
Thus, for any $j\geq 1$, define $$\hat{a}_{j}:=\left(\frac{1}{1+\tau^{n}\left(2^{1-N}(1+\sqrt{2}\tau)\right)}\right)^{d},$$
and notice that
$$a_{j+n}^{(1)}\,\geq\,\hat{a}_{j+n}a_{j+n-1}^{(1)}\,\geq\,\hat{a}_{l+n}\cdot \hat{a}_{j+n-1}\cdots \hat{a}_{n+1}\cdot a_{n}^{(1)}.$$
Thus taking logarithms we deduce that
\begin{apequation}\label{B02}
\log(a_{j+n}^{(1)})\,\geq \,\log \left ( \prod_{j\geq 1}\hat{a}_{j+n} \right )+\log(a_{n}^{(1)})\,\geq\, d\sum_{j\geq 1}\left(\frac{1}{1+\tau^{n+j}\left(2^{1-N}(1+\sqrt{2}\tau)\right)}\right)+\log(a_{n}^{(1)}).
\end{apequation}
Clearly $\{a_{n}^{(1)}\}$ is a decreasing sequence. If $\displaystyle\lim_{n\to\infty}a_{n}^{(1)}=\alpha_{1}$, then $\displaystyle\lim_{j\to\infty}\log(a_{j+n}^{(1)})=\log(\alpha_{1})$. Use inequality, $\log(1+x)<x$, $x>0$, we get

\begin{align*}
\log \alpha_1 &\geq -d \sum_{j=0}^{n}\left(\log\left(1+\tau^{j}\left(2^{1-N}(1+\sqrt{2}\tau)\right)\right)+\log(a_{n}^{(1)}) \right )\\ 
 &\geq -d\sum_{j=0}^{n}\tau^{j}\left(2^{N-1}(1+\sqrt{2}\tau)\right))+\log(a_{n}^{(1)}) \\ 
 &=-d\frac{\left(2^{1-N}(1+\sqrt{2}\tau)\right)\tau^{n}}{1-\tau}+\log(a_{n}^{(1)})\\
 &= \log\left(a_{n}^{(1)}\exp\left(-d\frac{2^{1-N}(1+\sqrt{2}\tau)}{1-\tau}\tau^{n}\right)\right)
\end{align*}

Since $a_{n}^{(1)}\geq a_{n}$, it follows that
\begin{apequation}\label{B03}
    \alpha_{1}\,\geq\, a_{n}^{(1)}\exp\left(-d\frac{2^{1-N}(1+\sqrt{2}\tau)}{1-\tau}\tau^{n}\right)
\end{apequation}
with $\tau \leq \tau ^{\ast}$, $n\in\mathbb{N_0\!}\,$.

\begin{figure}[h!]
    \centering
    \includegraphics[scale=0.6]{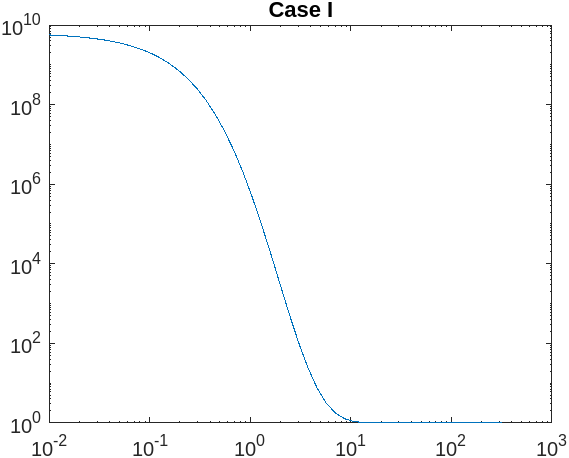}
    \caption{$\left\{\exp\left(-d\frac{2^{1-N}(1+\sqrt{2}\tau)}{1-\tau}\tau^{n}\right)\right\}$ in time and logarithmic scale}.
    \label{Esc-II}
\end{figure}

Now, taking the notion of $\Omega$ height in the definition A.4, we study a new case of the above such that the diameters of each $\Delta_i$ exceed $\delta\mathfrak{H}$ with $\delta$ sufficiently small $(\delta \ll 1)$, i.e., $\textrm{diam}(\Delta_i)\geq \delta\mathfrak{H}$ (recall the definition of $\mathfrak{H}$ given in (\ref{H})).\\

\noindent$\bullet$\,\,\underline{\textit{Case 2}}  As in Case 1, take now $\textrm{diam}\left(\bigcup_{i=1}^{k}\Delta_{i} \right )\geq \delta\mathfrak{H}$ with $\mathfrak{H}$ defined in Definition A.4, and $\delta$ is sufficiently small $(\delta \ll 1)$. Take $s\geq \delta\mathfrak{H}$ with $\tau\leq \tau^{\ast}$. Following the same procedure as above, we obtain that
$$a_{n}^{(2)}=\min_{1\leq k\leq 2^{n}}\min_{\Delta_{i} \in \Gamma_{n}^{\infty}}\min_{i=1, \ldots, k}\left\{\frac{\textrm{diam}\left(\bigcup_{i =1}^{k}\Delta_{i}\right)^{d}}{k\tau^{nd}}, \, \bigcup_{i=1}^{k}\Delta_{i}\,\,\,\,\textrm{satisfies   Case 2}\,\right\}.$$
Thus, as previously done, we have that
$$a_{n}^{(2)}\,\geq\,\left(\frac{1}{1+\frac{2\sqrt{2}\delta}{(\sqrt{2}+3\tau)}\tau^{n}} \right )^{d} a_{n-1}^{(2)}\,\,\,\,\,\,\,\,\textrm{and}\,\,\,\,\,\,\,\,\displaystyle\lim_{n\rightarrow\infty}a_{n}^{(2)}=\alpha_{2},$$
and thus
\begin{apequation}\label{B04}
    \alpha_{2}\,\geq\,a_{n}\exp\left(-\frac{2\sqrt{2}d}{\delta(1-\tau)(\sqrt{2}+3\tau)}\tau^{n}\right).
\end{apequation}
 Then comparing Case 1 and Case 2, we know that it is a better estimate if for each time $n\to\infty$ the argument of the exponential is larger. Define
\begin{apequation}
    \label{Delta_1}\xi_{1}^{(N)}:=\frac{2^{1-N}(1+\sqrt{2}\tau)}{1-\tau}
\,\,\,\,\,\,\,\,\textrm{and}\,\,\,\,\,\,\,\,
    \xi_{2}^{(\delta)}:=\frac{2\sqrt{2}}{\delta(1-\tau)(\sqrt{2}+3\tau)}
\end{apequation}
Note that the sequences $\{\xi_{1}^{(N)}\tau^{n}\}_{_{n\in\mathbb{N\!}_{\,0}}}$ and $\{\xi_{2}^{(\delta)}\tau^{n}\}_{_{n\in\mathbb{N\!}_{\,0}}}$ tends to zero eventually, since $\tau \leq \tau^\ast$, but $\xi_{1}^{(N)} \to 0$ as $N\to \infty$ and $\xi_{2}^{(\delta)} \to \infty$ as $\delta \to 0^+$.

\begin{figure}[h!]
    \centering
    \includegraphics[scale=0.5]{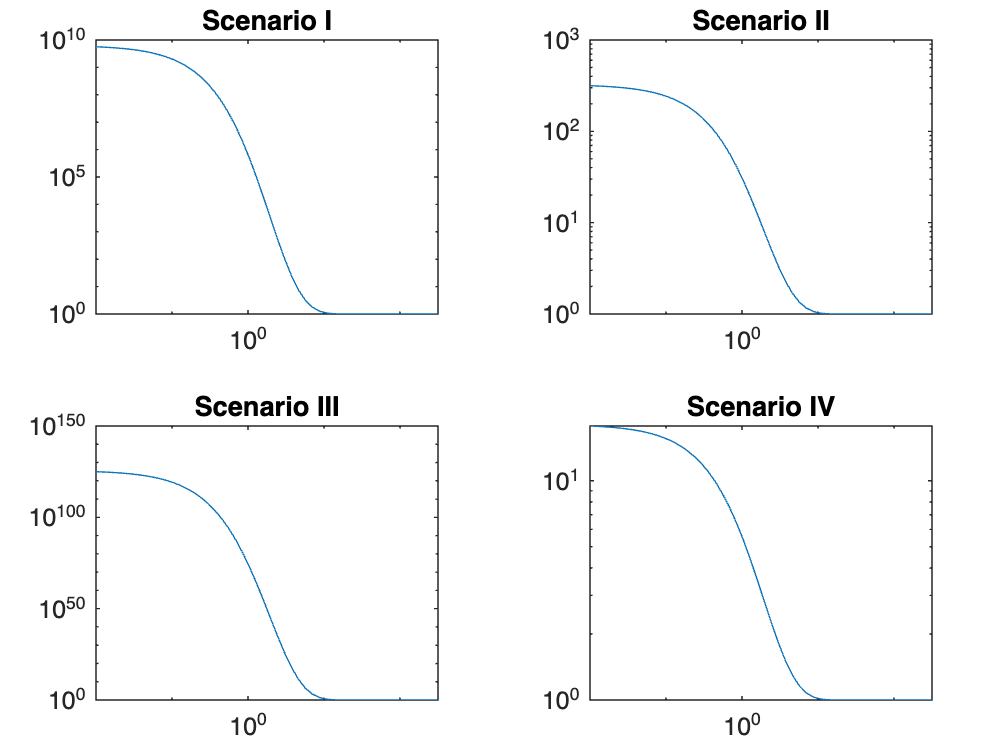}
    \caption{Different scenarios for case II. \textbf{Scenario I:} $\delta=\hat{\delta}_{\tau^*}$. \textbf{Scenario II:}  $\delta=0.5$. \textbf{Scenario III:}  $\delta=.01$. \textbf{Scenario IV:} $\delta=1$}.
    \label{Esc-IIb}
\end{figure}

In the figure \ref{Esc-IIb}, it can be seen that a smaller error is demonstrated for the scenario III, for values of
 $\delta\in [0, \hat{\delta}_{\tau^*}]$, the closer this $\delta$ is to $0$, then $a_{n}\xi_{2}^{(\delta)}\tau^{n}\nearrow\mathcal{H}^d(\Gamma^{\infty})$.\\

 \thinspace

 \noindent{\bf Remark A2.}\, {\it
    For $\tau_{1}, \tau_{2} \in [1/2, \tau^*]$, with $\tau_{1}>\tau_{2}$, it follows that $a_{n}\xi_{1}^{(N)}\tau_{1}^{n}$ and $a_{n}\xi_{2}^{(\delta)}\tau_{1}^{n}$  converge faster to $\mathcal{H}^d(\Gamma^{\infty})$, i.e., for the case $\tau=\tau^*$ one has the optimum for approximating this measure}.\\

\indent \underline{Case 1 and Case 2.} 

Let $\Delta := \{\delta \in (0, 1): \text{diam}(\bigcup_{i=1}^{k}\Delta_{i}) \geq \delta \mathfrak{H} \ \text{in the sense of case 1}\}$ be defined, where $\mathfrak{H}$ denotes a certain quantity. Take $\hat{\delta}_{\tau} := \min_{\delta \in \Delta}\{\delta\}$. Consequently, it is observed that $\xi_{2}^{\hat{\delta}_{\tau}} \leq \xi_{2}^{\delta}$ for all $\delta \in (0, 1)$. \\

 Recall that $a_n\nearrow\mathcal{H}^d(\Gamma^{\infty})$. Consequently,

\begin{apequation}
    \min\left\{a_{n}\exp\left(-d\xi^{(N)}_{1}\tau^{n}\right), a_{n}\exp\left(-d\xi^{(\delta)}_{2}\tau^{n}\right)\right\}:=\left\{\begin{matrix} 
a_{n}\exp\left(-d\xi^{(N)}_{1}\tau^{n}\right),  & \text{if} \ \delta>\hat{\delta}_{\tau}\\
a_{n}\exp\left(-d\xi^{(\delta)}_{2}\tau^{n}\right),  & \text{if} \ \delta\leq \hat{\delta}_{\tau}
\end{matrix}\right.
\label{min1}
\end{apequation}

\noindent We call $\hat{\delta}$ a \textit{control estimate} of $\delta$. From the above assuming that $\delta\leq \hat{\delta}_{\tau}$, then 

$$\lim_{n\to\infty}\{a_{n}^{(1)}, a_{n}^{(2)}\}\geq \min\left\{a_{n}\exp\left(-d\xi^{(\infty)}_{1}\tau^{n}\right), a_{n}\exp\left(-d\xi^{(\delta)}_{2}\tau^{n}\right)\right\}=a_{n}\exp\left(-d\xi^{(\delta)}_{2}\tau^{n}\right)$$

Let us now analyze the case where the family $\{\Delta_{i}\}_{i}$ is contained in $G_{i}\circ G_{i} \circ \cdots G_{i} (\Gamma^\infty)= G_{i}^{(N)}(\Gamma^\infty)$ with $i=1, 2$. The case is similar for both similarities since they are symmetries of each other, so it would be sufficient to study one of them. This means that each $\Delta_i$ will be contained in one of the symmetries (right or left) of the pre-fractals. In each of these cases we know that the diameter of each $\Delta_i$ will be greater than $2\tau$ because it is the measure of the base of the following hexagon. This is true, since the base always exceeds the measure of the other sides of the hexagon, and in this way we can say that $\textrm{diam}(\Delta_i)\geq 2\tau$, and therefore $\textrm{diam}\left(\bigcup_{i=1}^{k}\Delta_i\right)\geq 2\tau$.

\noindent$\bullet$\,\,\underline{\textit{Case 3.}}  If $\bigcup_{i=1}^{k}\Delta_{i} \subset G_{1}^{(N)}(\Gamma^{\infty})$. At that time, $s=\textrm{diam}\left(\bigcup_{i=1}^{k}\Delta_i\right)\geq 2\tau$ with $\tau \leq \tau^{\ast}$. In a similar way to the definition of $\{a_{n}^{(1)}\}$ (in Case 1), let
$$a_{n}^{(3)}=\min_{1\leq k\leq 2^{n}}\min_{\Delta_{i} \in \Gamma_{n}^{\infty}}\min_{i=1, \ldots, k}\left\{\frac{\textrm{diam}\left(\bigcup_{i =1}^{k}\Delta_{i}\right)^{d}}{k\tau^{nd}}, \, \bigcup_{i=1}^{k}\Delta_{i}\,\,\,\,\textrm{satisfies Case 3}\,\right\}$$
As before, it is easy to see that $\{a_{n}^{(2)}\}$ decreases. Letting $\displaystyle\lim_{n\to\infty}a_{n}^{(3)}=\alpha_3$, we get that
$$\frac{\textrm{diam}\left( \bigcup_{i=1}^{k}\Delta_{i} \right)}{\textrm{diam}\left ( \bigcup_{j=1}^{k_{n-1}}\Delta_{j}^{n-1} \right )}\,\geq\, \frac{s}{s+2\tau^{n}}\,\geq\, \frac{2\tau}{2\tau +2\tau^{n}}= \frac{1}{1+\tau^{n-1}}=\frac{1}{1+\frac{1}{\tau}\tau^{n}}$$
Proceeding in the same way as in the previous case, we arrive at
$$a_{n}^{(3)}\geq \left ( \frac{1}{1+\frac{1}{\tau}\tau^{n}} \right )^{d} a_{n-1}^{(3)}$$
Again we proceed in a similar way to obtain that
\begin{apequation}\label{B05}
    \alpha_{3}\geq a_{n}\exp\left(-\frac{d}{1-\tau}\tau^{n}\right).
\end{apequation}
In views of the definition of pre-fractals $\Omega^{(m)}$ $(m\in\mathbb{N\!}_{\,0})$ given in (\ref{Omega-ramified}), then we can consider the possibility that for some $m_0\in\mathbb{N\!}_{\,0}$, one gets that $\bigcup_{i=1}^{k}\Delta_{i}$ intercepts with $\overline{\Omega^{
(m_0)}}$. Moreover, there may exist a family of closed sets $\{\overline{\Omega^{(m)}}\}_{_{m\in\mathbb{N\!}_{\,0}}}$ such that $\bigcup_{i=1}^{k}\Delta_{i}$ intercepts with $\overline{\Omega^{(m)}}$ for each $m\in \mathfrak{A} \subset \mathbb{N\!}_{\,0}$ but in other cases they do not intercept. This leads to the following:\\

\noindent$\bullet$\,\,\underline{\textit{Case 4.}}  If $\bigcup_{i=1}^{k}\Delta_{i} \cap \overline{\Omega^{(m_0)}}\neq \emptyset$ and $\bigcup_{i=1}^{k}\Delta_{i} \cap \overline{\Omega^{(m'_0)}}= \emptyset$ with $m_0,\,m'_0\in\mathbb{N}$, \,$m_0 \neq m'_0$, then at that time, $s=\textrm{diam}\left(\bigcup_{i=1}^{k}\Delta_i\right)\geq (2\tau)^{\mathfrak{m}+1}$ for $\mathfrak{m}:=\max\{m_0,m'_0\}$ and $\tau \leq \tau^{\ast}$. Following analogously as in Cases 1 and 2, define the decreasing sequence $\{a^{(3)}_n\}$ by:
$$a_{n}^{(4)}=\min_{1\leq k\leq 2^{n}}\min_{\Delta_{i} \in \Gamma_{n}^{\infty}}\min_{i=1, \ldots, k}\left\{\frac{\textrm{diam}\left(\bigcup_{i =1}^{k}\Delta_{i}\right)^{d}}{k\tau^{nd}}, \, \bigcup_{i=1}^{k}\Delta_{i}\,\,\,\,\textrm{satisfies Case 4}\,\right\}$$
Given that $\displaystyle\lim_{n\to\infty}a_{n}^{(4)}=\alpha_4$, proceeding similarly as before we are lead into

$$a_{n}^{(4)}\geq \left ( \frac{1}{1+\frac{2}{(2\tau)^{\mathfrak{m}+1}}\tau^{n}} \right )^{d} a_{n-1}^{(4)}$$

\noindent and therefore, 

\begin{apequation}\label{B06}
    \alpha_{4}\,\geq\,a_{n}\exp\left(-\frac{d}{(1-\tau)\tau^{\mathfrak{m}}}\tau^{n}\right) 
\end{apequation}
Combining (\ref{B04}), (\ref{B05}) and (\ref{B06}), we conclude that

\begin{align*}
 \mathcal{H}^{d}(\Gamma^{\infty}) &= \lim_{n\to\infty}\min\{a_{n}^{(2)}, a_{n}^{(3)}, a_{n}^{(4)}\}\\ 
 &\geq \min \left \{ a_{n}\exp\left(-\frac{2\sqrt{2}d}{\delta(1-\tau)(\sqrt{2}+3\tau)}\tau^{n}\right), a_{n}\exp\left(-\frac{d}{1-\tau}\tau^{n}\right), a_{n}\exp\left(-\frac{d}{(1-\tau)\tau^{\mathfrak{m}}}\tau^{n}\right) \right \}
\end{align*}

\noindent Now, it is clear that if we take $\mathfrak{m}=0$, Cases 3 and 4 are equivalent, so we will examine Cases 2 and Case 3. For this, we define
\begin{apequation}\label{B07}
    \xi^{(\mathfrak{m})}:=\frac{1}{(1-\tau)\tau^{\mathfrak{m}}}
\end{apequation}

Note that $\xi^{(0)}$ coincides with the case 3. It is clear that $\xi_{1}^{(N)}, \xi^{(\mathfrak{m})} > 0$ (see (\ref{Delta_1})). With this in mind, we now check what eventually happens with $\xi_{1}^{(N)}$ and $\xi^{(\mathfrak{m})}$. 

\thinspace

\noindent$\bullet$\,\,\underline{Case 2 and Case 4 when $\tau=\tau^{\ast}$.}  It can be seen that the value that $\mathfrak{m}$ must take to coincide with $\xi^{(\hat{\delta})}_{2}$is the integer part of

$$\mathfrak{m}=\frac{\log\left(\frac{\hat{\delta}_{\tau^*}(\sqrt{2}+3\tau)}{2\sqrt{2}}\right)}{\log(\tau)} \approx 3.70431$$

\noindent If $\tau=\tau^{\ast}$, then $\lfloor \mathfrak{m} \rfloor=3$. Note that for $\mathfrak{m}\leq 3$ we have that $\xi^{(\mathfrak{m})}<\xi_{2}
^{(\delta)}$, and for the case where $\mathfrak{m}>3$ the opposite case. Similarly as discussed in (\ref{min1}), we obtain that

\begin{apequation}
    \min\left\{a_{n}\exp\left(-d\xi^{(\delta)}_{2}\tau^{n}\right), a_{n}\exp\left(-d\xi^{(\mathfrak{m})}\tau^{n}\right)\right\}:=\left\{\begin{matrix}
a_{n}\exp\left(-d\xi^{(\delta)}_{2}\tau^{n}\right),  & \text{if} \ \mathfrak{m}\leq 3, \ \ n>\mathfrak{m}\\ 
a_{n}\exp\left(-d\xi^{(\mathfrak{m})}\tau^{n}\right),  & \text{if} \ \mathfrak{m} >3, \ \ n\leq \mathfrak{m}
\end{matrix}\right.
\end{apequation}

\begin{figure}[h!]
    \centering
    \includegraphics[scale=0.5]{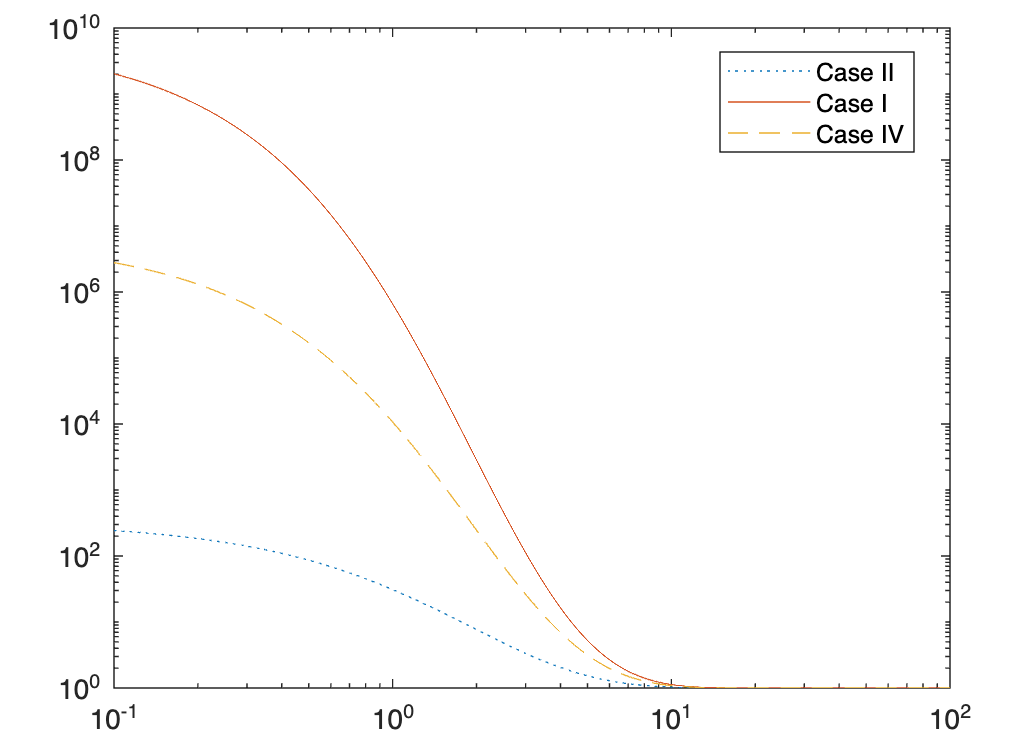}
    \caption{Convergence analysis of Cases 1, 2 and 4.}.
    \label{Esc-VARIOS}
\end{figure}

 \noindent Thus, one can rule out the case $n \leq 3$, so to try to find the terms of $\{b_{n}\}$ it is necessary to take the case when $\mathfrak{m}\leq 3$, accordingly

$$\mathcal{H}^{d}(\Gamma^{\infty})\geq \min\left\{a_{n}\exp\left(-d\xi^{(\delta)}_{2}\tau^{n}\right), a_{n}\exp\left(-d\xi^{(\mathfrak{m})}\tau^{n}\right)\right\} =a_{n}\exp\left(-d\xi^{(\delta)}_{2}\tau^{n}\right)=a_n\,\exp\left(-\frac{2\sqrt{2}d\delta}{(1-\tau)(\sqrt{2}+3\tau)}\tau^{n}\right)$$

Note that in the figure \ref{Esc-VARIOS} the sequence associated to case II is better compared to the other two cases, because it starts from a lower error, considering that they all converge to zero using the same number of iterations. Moreover, in case the program stops at a finite number of iterations, the error of the sequence in case II is smaller. Therefore, 

\begin{apequation}  a_{n}\xi_{2}^{(\delta)}\tau^{n}\,\leq\,\mathcal{H}^d(\Gamma^{\infty})\,\leq\,a_n\,\,\,\,\,\,\,\textrm{and}\,\,\,\,\,\,a_{n}\xi_{2}^{(\delta)}\tau^{n}\nearrow\mathcal{H}^d(\Gamma^{\infty})\swarrow a_{n}, \ \ \delta\in[0, \hat{\delta}_{\tau}]
    \label{estimation}
\end{apequation}

i.e., 

\begin{apequation}
    a_n\,\exp\left(-\frac{2\sqrt{2}d}{\delta(1-\tau)(\sqrt{2}+3\tau)}\tau^{n}\right)\leq\mathcal{H}^{d}(\Gamma^{\infty})\leq a_{n}, \ \ \delta\in[0, \hat{\delta}_{\tau}],
\end{apequation}

Note also that $\exp\left(\xi_{2}^{(\delta)}\tau^{n}\right)\leq 1$ for $\delta\leq \hat{\delta}_{\tau}, $ then $b_{n}:=a_{n}\exp\left(\xi_{2}^{(\delta)}\tau^{n}\right)\leq 1$ with $0\leq a_n\leq 1$, for all $n\in\mathbb{N}_0$. Thus, $\{b_n\}\subseteq [0, 1]$, which concludes the proof of Theorem \ref{Main-T1}. $\hfill\square$\\


\noindent{\bf Remark A3.}\, {\it
One can calculate computationally that $a_{1}\approx 0.894684280597$.
Consequently, if one considers the critical case $\tau=\tau^{\ast}$ and $\delta=\hat{\delta}_{\tau^*}$, we obtain the approximation of $\mathcal{H}^d(\Gamma^{\infty})$ mentioned in (\ref{a_1-computation})}.\\

\noindent{\bf Acknowledgements.}\, We would like to thank Dr. Yves Achdou and Dr. Nicoletta Tchou for their help in providing us with suitable images of ramified domains.\\

\noindent The second author is supported by:
The Puerto Rico Science, Technology and Research Trust, under agreement number 2022-00014.\\

\indent\\

\noindent{\bf Disclaimer.}  {\it This content is only the responsibility of the authors and does not necessarily represent the official views of The Puerto Rico Science, Technology and Research Trust}.\\

\end{document}